\documentclass[11pt]{article}

\usepackage{amsmath,amsfonts,amssymb,graphicx}

\newtheorem{lemma}{Lemma}[section]

\newtheorem{prop}{Proposition}[section]
\newtheorem{thm}{Theorem}[section]
\newtheorem{cor}{Corollary}[section]
\newtheorem{defn}{Definition}[section]
\newtheorem{claim}{Claim}[section]

\numberwithin{equation}{section}

\newcommand{\be}{\begin}
\newcommand{\refp}[1]{(\ref{#1})}

\newcommand\ca{\mathcal A}
\newcommand{\cb}{\mathcal B}

\newcommand\ce{\mathcal E}
\newcommand{\cg}{\mathcal G}

\newcommand\ck{\mathcal K}
\newcommand\cl{\mathcal L}
\newcommand{\cm}{\mathcal M}

\newcommand\fl{\mathcal R}

\newcommand\ct{\mathcal T}

\newcommand{\al}{\alpha}
\newcommand{\del}{\delta}
\newcommand\Del{\Delta}
\newcommand{\eps}{\epsilon}
\newcommand\ga{\gamma}
\newcommand\Ga{\Gamma}
\newcommand\ka{\kappa}
\newcommand\lla{\lambda}

\newcommand\ups{\upsilon}
\newcommand\om{\omega}

\newcommand{\si}{\sigma}
\newcommand\Si{\Sigma}

\newcommand\bbc{\mathbb C}

\newcommand\bbe{\mathbb E}

\newcommand\q{\mathbb Q}
\newcommand{\R}{\mathbb R}

\newcommand\bbu{\mathbb U}

\newcommand{\z}{\mathbb Z}

\newcommand\frg{\mathfrak g}

\newcommand\frp{\mathfrak p}

\newcommand\frP{\mathfrak P}
\newcommand\frR{\mathfrak R}

\newcommand\bfh{\mathbf h}

\newcommand\RP[1]{{\mathbb{RP}^{#1}}}
\newcommand\CP[1]{{\mathbb{CP}^{#1}}}

\newcommand{\uset}{\underset}
\newcommand{\oset}{\overset}
\newcommand{\uline}{\underline}
\newcommand{\oline}{\overline}

\newcommand{\la}{\langle}
\newcommand{\ra}{\rangle}
\newcommand{\st}{\,|\,}
\newcommand{\ti}{\tilde}

\newcommand{\prtl}{\partial}
\renewcommand\square{\kern20pt{\vbox{\hrule height.4pt
        \hbox{\vrule width.4pt height 6pt\kern6pt
                \vrule width.4pt}
        \hrule height.4pt}}}
\newcommand\im{\text{im}}
\renewcommand\Im{\text{Im}}
\renewcommand\Re{\text{Re}}
\newcommand\rank{\text{rank}}

\newcommand\inv{^{-1}}
\newcommand\hol{\text{Hol}}

\newcommand\const{\text{const}}

\newcommand\proof{{\em Proof.}\ }
\newcommand\loc{{\text{loc}}}

\newcommand{\ry}{\R\times Y}

\newcommand\U[1]{\text{U}(#1)}
\newcommand\SU[1]{\text{SU}(#1)}

\newcommand\SO[1]{\text{SO}(#1)}

\newcommand\cs{\vartheta}
\newcommand\ind{\text{ind}}

\newcommand\red{^{\text{red}}}

\newcommand\torsion{\text{torsion}}

\newcommand\lloc[2]{L^#1_{#2,\loc}}
\newcommand\llw[3]{{L^{#1,#2}_#3}}

\newcommand\id{\text{Id}}

\newcommand\chm{\check M}
\newcommand\fct{\frac1{32\pi^2}}
\newcommand\cbs{\cb^*(Y[0])}
\newcommand\mab{M(\al,\beta)}

\newcommand\mtb{M(\theta,\beta)}
\newcommand\ufl{\uline\fl^*}
\newcommand\ual{{\uline\alpha}}
\newcommand\ubeta{{\uline\beta}}
\newcommand\uga{{\uline\gamma}}
\newcommand\yb[1]{Y[#1]}
\newcommand\qt{q_2}

\newcommand\ix[2]{I^{#1}_{#2}}

\newcommand\mx[2]{M_{#1,#2}}
\newcommand\mxl[2]{M_{#1,#2}^L}
\newcommand\tmx[2]{\ti M_{#1,#2}}

\newcommand\energy[2]{\ce_{#1}(#2)}
\newcommand\ttp{\tilde\tau{{\vphantom\tau}^+}}
\newcommand\trp[1]{\tilde\rho{{\vphantom\rho}^+_{#1}}}
\newcommand\tta{\tilde\tau_a}
\newcommand\ttd{\tilde\tau_d}

\newcommand\tn[1]{t_{n,#1}}
\newcommand\hrp{\hat R{\vphantom R}^+}

\newcommand\omn{\om_n}
\newcommand\flag{\fl_{\al'\ga}}
\newcommand\bbei[2]{\bbe(#1,#2)}
\newcommand\ubbe{\uline{\bbe}}
\newcommand\uell{\uline{\ell}}
\newcommand\cmsi{\cm_\Si}
\newcommand\cmcyl{\cm_{\text{cyl}}}
\newcommand\tcm{\ti\cm}
\newcommand\tcmsi{\ti\cm_\Si}
\newcommand\tcmcyl{\ti\cm_{\text{cyl}}}
\newcommand\tcmint{\ti\cm_{\text{int}}}
\newcommand\cmlsi{\cm^L_\Si}
\newcommand\cmlcyl{\cm^L_{\text{cyl}}}
\newcommand\wsi{W_\Si}
\newcommand\bbesi{\bbe_\Si}
\newcommand\wcyl{W_{\text{cyl}}}
\newcommand\nsi{N_\Si}
\newcommand\ncyl{N_{\text{cyl}}}
\newcommand\cnc{{\mathcal N}_C}
\newcommand\sco[1]{_{:#1}}
\newcommand\va{\vec\al}
\newcommand\uz[3]{{\uline z}_{#1#2#3}}
\newcommand\oz[3]{{\oline z}_{#1#2#3}}
\newcommand\zz[3]{z_{#1#2#3}}
\newcommand\oc[1]{{\oline C}^#1}
\newcommand\hsp{\hspace{2cm}}
\newcommand\ath{A^\theta}
\newcommand\ev{_{\text{ev}}}

\newcommand\npart{/\kern-.19cm\partial}
\newcommand\bev{\cb_{\text{ev}}}
\newcommand\xend{X_{\text{end}}}
\newcommand\twod[2]{D_{#1,#2}}
\newcommand\hcob{\theta^3_{\z}}
\newcommand\pt[1]{T_#1}
\newcommand\uEue{\uline E\otimes\uline\ell}

\begin{document}

\title{Mod~$2$ instanton homology and $4$-manifolds\\with boundary}

\author{Kim A.\ Fr\o yshov}

\date{}

\maketitle

\begin{abstract}
Using instanton homology with coefficients in $\z/2$ we construct a
homomorphism $\qt$ from the homology cobordism group $\hcob$ to the integers
which is not a rational linear combination of
the instanton $h$--invariant and the Heegaard Floer
correction term $d$. If an oriented
homology $3$--sphere $Y$ bounds a smooth, compact,
negative definite $4$--manifold without $2$--torsion in its homology
then $\qt(Y)\ge0$, with strict inequality if the intersection form 
is non-standard.
\end{abstract}

\thispagestyle{empty}

\bibliographystyle{plain}

\tableofcontents

\section{Introduction}

This paper will introduce an integer invariant $\qt(Y)$ of oriented
integral homology $3$--spheres $Y$. This invariant is defined in terms of
instanton cohomology with coefficients in $\z/2$ and may be regarded as a
mod~$2$ analogue of the $h$--invariant \cite{Fr3}, which was defined with
rational coefficients. Both invariants grew out of efforts to extend
Donaldson's diagonalization theorem \cite{D1,D2} to $4$--manifolds with
boundary.

We will use the instanton (co)homology originally introduced by Floer
\cite{F1}, an exposition of which can be found in \cite{D5}. With
coefficients in $\z/2$, instanton cohomology $I(Y;\z/2)$ comes equipped with
some extra structure, namely two ``cup products'' $u_2$ and $u_3$ of degrees
$2$ and $3$, respectively, and homomorphisms
\[I^4(Y;\z/2)\oset{\del_0}\longrightarrow \z/2\oset{\del_0'}\longrightarrow I^1(Y;\z/2)\]
counting index~$1$ trajectories running into and from the trivial flat
$\SU2$ connection, respectively.
This extra structure enters in the definition of the invariant $q_2$.
Reversing the r\^oles of the cup products $u_2,u_3$ in the definition
yields another invariant $q_3$. However, the present paper will focus on
$\qt$.

It would be interesting to try to express the invariants $h,q_2,q_3$ in terms of
the equivariant instanton homology groups recently introduced by Miller Eismeier
\cite{Miller-Eismeier1}. 

\subsection{Statement of main results}

\begin{thm}[Additivity]
  For any oriented homology $3$--spheres $Y_0$ and $Y_1$ one has
  \[\qt(Y_0\#Y_1)=\qt(Y_0)+\qt(Y_1).\]
\end{thm}

The proof of additivity is not quite straightforward and occupies more
than half the paper.

\be{thm}[Monotonicity]\label{thm:monotone}
Let $W$ be a smooth compact oriented $4$-manifold with boundary
$\prtl W=(-Y_0)\cup Y_1$, where $Y_0$ and $Y_1$ are oriented homology
$3$--spheres. Suppose the intersection form of $W$ is negative
definite and $H^2(W;\z)$ contains no element of order $4$. Then
\[\qt(Y_0)\le\qt(Y_1).\]
\end{thm}

If the manifold $W$ in the theorem actually satisfies $b_2(W)=0$ then one can
apply the theorem to $-W$ as well so as to obtain $\qt(Y_0)=\qt(Y_1)$.
This shows that $\qt$ descends to a group homomorphism
$\hcob\to\z$, where $\hcob$ is the integral homology cobordism group.

We observe that the properties of $\qt$ described so
far also hold for the instanton
$h$--invariant, the negative of its monopole analogue \cite{Fr4,KM5}, and
the Heegaard Floer correction term $d$ \cite{OS2}. Note that the latter three
invariants are monotone
with respect to any negative definite cobordism, without any assumption on the
torsion in the cohomology.


\be{thm}[Lower bounds]\label{thm:defub}
Let $X$ be a smooth compact oriented $4$-manifold whose boundary
is a homology sphere $Y$. Suppose the intersection form of $X$ is negative
definite and $H^2(X;\z)$ contains no $2$-torsion. Let
\[J_X:=H^2(X;\z)/torsion,\]
and let $w$ be an element of $J_X$ which is not divisible by $2$.
Let $k$ be the minimal square norm (with
respect to the intersection form) of any element
of $w+2J_X$. Let $n$ be the number of elements of $w+2J_X$ of square norm $k$.
If $k\ge2$ and $n/2$ is odd then
\be{equation*}
\qt(Y)\ge k-1.
\end{equation*}
\end{thm}


By an {\em integral lattice} we mean a free abelian group of finite rank
equipped with
a symmetric bilinear integer-valued form. Such a lattice is called
{\em odd} if it contains an element of odd
square; otherwise it is called {\em even}.

\be{cor}\label{cor:JtiJ}
Let $X$ be as in Theorem~\ref{thm:defub}.
Let $\ti J_X\subset J_X$ be the orthogonal complement of the sublattice of $J_X$
spanned by all vectors of square $-1$, so that $J_X$ is an orthogonal sum
\[J_X=m\la-1\ra\oplus\ti J_X\]
for some non-negative integer $m$.
\be{description}
\item[(i)]If $\ti J_X\neq0$, i.e.\ if $J_X$ is not diagonal, then $\qt(Y)\ge1$.
\item[(ii)]If $\ti J_X$ is odd then $\qt(Y)\ge2$.
\end{description}
\end{cor}

To deduce (i) from the theorem, take $C:=v+2J_X$ where $v$ is any non-trivial
element of
$\ti J_X$ of minimal square norm. To prove (ii), choose a $v$ with minimal odd
square norm.

\be{thm}\label{thm:surgery-bound}
Let $Y$ be the result of $(-1)$ surgery on a knot $K$ in $S^3$. If changing
$n^-$ negative crossings in a diagram for $K$ produces a positive knot then
\[0\le\qt(Y)\le n^-.\]
\end{thm}

For $k\ge2$ the Brieskorn sphere $\Si(2,2k-1,4k-3)$ is the boundary of a
plumbing manifold with intersection form $-\Ga_{4k}$ (see
Section~\ref{sec:further-props}), and it is also
the result of $(-1)$
surgery on the $(2,2k-1)$ torus knot. In these examples
the upper bound on $\qt$ given by
Theorem~\ref{thm:surgery-bound} turns out to coincide with the lower bound
provided by Theorem~\ref{thm:defub}, and one obtains the following.

\begin{prop}\label{prop:Brieskorn-calc}
   For $k\ge2$ one has
\[\qt(\Si(2,2k-1,4k-3))=k-1.\]
\end{prop}

On the other hand, by \cite[Proposition~1]{Fr7} one has
\[h(\Si(2,2k-1,4k-3))=\lfloor k/2\rfloor,\]
and in these examples the correction term $d$ satisfies $d=h/2$, as follows
from \cite[Corollary~1.5]{OS6}. This shows:

\begin{prop}
  The invariant $\qt$ is not a rational linear combination of the
  $h$--invariant and the correction term $d$.\square
\end{prop}
In particular,
\[h,\qt:\hcob\to\z\]
are linearly independent homomorphisms, and the same is true for $d,\qt$.
It follows from this that $\hcob$ has a $\z^2$ summand. However, much more is
true: Dai, Hom, Stoffregen, and Truong \cite{DHST1} proved that $\hcob$
has a $\z^\infty$ summand. Their proof uses involutive Heegaard Floer homology.

The monotonicity of the invariants $h,d,\qt$ leads to the following result.

\begin{thm}\label{thm:lin-indep}
Let $Y$ by an oriented homology $3$-sphere. If
\[\min(h(Y),d(Y))<0<\qt(Y)\] 
then $Y$ does not bound any definite $4$-manifold without elements of order~$4$
in its second cohomology.
\end{thm}

An explicit example to which the theorem applies is $2\Si(2,5,9)\#-3\Si(2,3,5)$.

A related result was obtained by Nozaki, Sato, and Taniguchi \cite{NST1}.
Using a filtered version of instanton homology they proved that certain linear
combinations of Brieskorn homology $3$--spheres do not bound any definite
$4$--manifold.

\begin{thm}\label{thm:h-qt-torsion}
If an oriented homology $3$-sphere $Y$ satisfies
\[h(Y)\le0<\qt(Y)\]
then $I^5(Y;\z)$ contains $2$--torsion, hence $Y$ is not homology cobordant
to any Brieskorn sphere $\Si(p,q,r)$.
\end{thm}

We conclude this introduction with two sample applications of the invariant
$\qt$.

\begin{thm}\label{thm:Poincare-sphere}
Let $X$ be a smooth compact oriented connected $4$-manifold whose boundary
is the Poincar\'e sphere $\Si(2,3,5)$.
Suppose the intersection form of $X$ is negative definite.
Let $\ti J_X$ be as in Corollary~\ref{cor:JtiJ}.
\begin{description}
\item[(i)] If $\ti J_X$ is even then $\ti J_X=0$ or $-E_8$.
\item[(ii)] If $\ti J_X$ is odd then $H^2(X;\z)$
  contains an element of order $4$.
\end{description}
\end{thm}

Earlier versions of this result were obtained using instanton homology
in \cite{Fr0} (assuming $X$
is simply-connected) and in \cite{Scaduto2} (assuming $X$ has no $2$--torsion
in its homology).

There are up to isomorphism two even, positive
definite, unimodular forms of rank~$16$, namely $2E_8$ and $\Ga_{16}$.
If $Z$ denotes the negative definite $E_8$--manifold then the boundary
connected sum $Z\#_{\prtl}Z$ has intersection form $-2E_8$.
It is then natural to ask whether $\Si(2,3,5)\#\Si(2,3,5)$ also bounds
$-\Ga_{16}$.
There appears to be no obstruction to this coming from
the correction term. 

\begin{thm}\label{thm:sisi}
Let $X$ be a smooth compact oriented $4$-manifold whose boundary
is $\Si(2,3,5)\#\Si(2,3,5)$. Suppose the intersection form of $X$ is negative
definite and $H^2(X;\z)$ contains no $2$--torsion. If $\ti J_X$ is even then
\[\text{$\ti J_X=0$, $-E_8$, or $-2E_8$.}\]
\end{thm}

Further results on the definite forms bounded by a given homology $3$--sphere
were obtained by Scaduto \cite{Scaduto2}.

Some of the results of this paper were announced in various talks several years
ago. The author apologizes for the long delay in publishing the results.

\subsection{Outline}

  To learn the definition of $\qt$ the reader may proceed directly
  to Sections~\ref{sec:inst-cohom} and \ref{section:defqt} and refer back
  to earlier sections for notation and set-up.
  The relationship of the invariant $\qt$ to definite $4$--manifold is
  discussed in Section~\ref{sec:definite-4mflds}. Note that in this section
  we work with $\SO3$ connections modulo all automorphisms of the bundle,
  not just those that lift to $\SU2$.

  The proof of additivity of $\qt$ is given in
  Sections~\ref{sec:operations}, \ref{sec:two-points-I}, and
  \ref{sec:two-points-II}. It involves various operations on instanton
  cohomology defined by cobordisms. The first two such operations, denoted
  $\phi$ and $\psi$, appear
  in Subsection~\ref{subsec:comm-cup-products}. In these cases the cobordism
  is just a cylinder. Proofs of the main properties of $\phi$ and $\psi$
  are deferred to the last two sections~\ref{sec:two-points-I}
  and \ref{sec:two-points-II}, which form the technically most
  difficult part of the paper. Operations on instanton cohomology defined
  by cobordisms with three boundary components are discussed in
  Section~\ref{sec:operations}, leading to a proof of
  additivity of $\qt$ assuming the results of the last two sections.
  
  The remaining results stated in this introduction are proved in
  Section~\ref{sec:further-props}.

{\bf Acknowledgement.} The author would like to thank Tom Mrowka for
helpful conversations.

\section{The base-point fibration}\label{sec:base-point-fibration}

Let $X$ be a connected smooth $n$--manifold, possibly with boundary, and
$P\to X$ a principal $\SO3$ bundle. Fix $p>n$ and let $A$ be a $\lloc p1$
connection in $P$. This means that $A$ differs from a smooth connection by a
$1$--form which lies locally in $L^p_1$. Let $\Ga_A$ be the group of
$\lloc p2$ automorphisms (or gauge transformations) of $P$ that preserve $A$.
The connection $A$ is called
\begin{itemize}
\item {\em irreducible} if $\Ga_A=\{1\}$, otherwise {\em reducible};
\item {\em Abelian} if $\Ga_A\approx\U1$;
\item {\em twisted reducible} if $\Ga_A\approx\z/2$.
\end{itemize}
Note that a non-flat reducible connection in $P$ is either
Abelian or twisted reducible.

Recall that automorphisms of $P$ can be regarded as sections of the bundle
  $P\uset{\SO3}\times\SO3$ of Lie groups, where $\SO3$ acts on itself by
  conjugation.  An automorphism is called {\em even} if it lifts to a
  section of $P\uset{\SO3}\times\SU2$. A connection $A$ in $P$ is called
  {\em even-irreducible} if its stabilizer $\Ga_A$ contains no non-trivial
  even automorpisms, otherwise $A$ is called {\em even-reducible}.
A non-flat connection is even-reducible if and only if it is Abelian.

  Now suppose $X$ is compact and let $\ca$ be the space of 
  all $L^p_1$ connections in $P$. The affine Banach space $\ca$ is acted upon
  by the Banach Lie group $\cg$ consisting of all $L^p_2$ automorphisms of $P$.
  Let $\ca^*\subset\ca$ be subset of irreducible connections and define
  $\cb=\ca/\cg$. The irreducible part $\cb^*\subset\cb$ is a Banach manifold,
  and it admits smooth partitions of unity provided $p>n$ is an even integer,
  which we assume from now on. Instead of $\cb^*$ we often write $\cb^*(P)$,
  or $\cb^*(X)$ if the bundle $P$ is trivial. Similarly for $\ca,\cg$ etc. 

  Let $\ca^*\ev$ be the space of all even-irreducible
   $L^p_1$ connections in $P$.
  Let $\cg\ev$ be the group of even $\lloc p2$ 
  automorphisms of $P$. As explained in \cite[p.\ 235]{BD1}, there is an exact
  sequence
  \[1\to\cg\ev\to\cg\to H^1(X;\z/2)\to0.\]
  The quotient $\cb^*\ev=\ca^*\ev/\cg\ev$ is a Banach manifold.

\begin{defn}\label{defn:admissible}
Let $X$ be a topological space.
\begin{description}
\item[(i)] A class $v\in H^2(X;\z/2)$ is called {\em admissible} if
  $v$ has a non-trivial pairing with a class in $H_2(X;\z)$, or equivalently,
  if there exist a closed oriented $2$--manifold $\Si$ and a continuous map
  $f:\Si\to X$ such that $f^*v\neq0$. If $\Si$ and $f$ can be chosen such that,
  in addition,
\begin{equation}\label{eqn:fa0}
  f^*a=0\quad\text{for every $a\in H^1(X;\z/2)$,}
  \end{equation}
  then $v$ is called
  {\em strongly admissible}.
\item[(ii)] An $\SO3$ bundle $E\to X$ is called
  {\em (strongly) admissible} if
  the Stiefel-Whitney class $w_2(E)$ is (strongly) admissible.
  \end{description}
\end{defn}

For example, a finite sum $v=\sum_ia_i\cup b_i$ with
$a_i,b_i\in H^1(X;\z/2)$ is never strongly admissible.

\begin{prop}\label{prop:lift-to-U2}
  Let $X$ be a compact, oriented, connected smooth $4$--manifold with base-point
  $x\in X$. Let $P\to X$ be an $\SO3$ bundle.
\begin{description}
\item[(i)] If $P$ is admissible then the $\SO3$ base-point fibration over
  $\cb^*\ev(P)$ lifts to a $\U2$ bundle.
\item[(ii)] If $P$ is strongly admissible then the $\SO3$ base-point
  fibration over $\cb^*(P)$ lifts to a $\U2$ bundle.
\end{description}
  \end{prop}

\proof We spell out the proof of (ii), the proof of (i) being similar
(or easier). Let $\Si$ be a closed oriented surface and $f:\Si\to X$
a continuous map such that $f^*P$ is non-trivial and \refp{eqn:fa0} holds.
We can clearly arrange that $\Si$ is
connected. Because $\dim X\ge2\dim\Si$ it follows
from \cite[Theorems~2.6 and 2.12]{Hirsch} that $f$ can be uniformly approximated
by (smooth) immersions $f_0$. Moreover, if the approximation is sufficiently
good then $f_0$ will be homotopic to $f$. Therefore, we may assume $f$ is an
immersion. 
Since base-point fibrations associated to different base-points in
$X$ are isomorphic we may also assume that $x$ lies in the image of $f$,
say $x=f(z)$.

We adapt the proof of \cite[Proposition~2.6]{Kotsch1}, see also
\cite[Proposition~2.3]{FS2}. Let $\bbe\to\cb^*:=\cb^*(P)$ be the
oriented Euclidean $3$--plane bundle associated to the
base-point fibration. We must find an Hermitian 2-plane bundle
$\ti\bbe$ such that
$\bbe$ is isomorphic to the bundle $\frg^0_{\ti\bbe}$ 
of trace-free
skew-Hermitian endomorphisms of $\ti\bbe$.

Let $E\to X$ be the standard $3$--plane bundle associated to $P$.
Choose an Hermitian $2$--plane bundle $W\to\Si$ together with an isomorphism
$\phi:\frg^0_W\oset\approx\to f^*E$, and fix a connection $A_{\Si,\text{det}}$
in $\det(W)$. Any (orthogonal)
connection $A$ in $E$ induces a connection in $f^*E$
which in turn induces a connection $A_\Si$ in $W$ with central part
$A_{\Si,\text{det}}$. Choose a spin structure on $\Si$ and let $S^*\pm$ be
the corresponding spin bundles over $\Si$. For any
connection $A$ in $E$ let
\[\npart_{\Si,A}:S^+\otimes W\to S^-\otimes W\]
be the Dirac operator
coupled to $A_\Si$. If $A$ is an $L^p_1$ connection, $p>4$, and $A_0$ is a
smooth connection in $E$ then $A-A_0$ is continuous, hence
$\npart_{\Si,A}-\npart_{\Si,A_0}$ defines a bounded operator $L^2\to L^2$ and
therefore a compact operator $L^2_1\to L^2$. Let
\[\cl:=\det\ind(\npart_{\Si,W})\]
be the determinant line bundle over $\ca(E)$
associated to the family of Fredholm operators
\[\npart_{\Si,A}:L^2_1\to L^2.\]
Then automorphism $(-1)$ of $W$ acts on $\cl$ with weight equal to the
numerical index of $\npart_{\Si,A}$. According to Atiyah-Singer's theorem
\cite{AS1} this index is
\[\ind(\npart_{\Si,A})=\{\text{ch}(W)\hat A(\Si)\}\cdot[\Si]=c_1(W)\cdot[\Si].\]
But the mod~$2$ reduction of $c_1(W)$ equals $f^*(w_2(E))$,
which is non-zero by assumption, so the index is odd.

The assumption \refp{eqn:fa0} means that every automorphism of $E$
pulls back to an
{\em even} automorphism of $f^*E$. Moreover, every even automorphism of
$f^*E\approx\frg^0_W$
lifts to an automorphism of $W$ of determinant~$1$, the lift being well-defined
up to an overall sign since $\Si$ is connected. Because the automorphism
$(-1)$ of $W$ acts trivially on $\cl\otimes W_z$ this yields an action of
$\cg(E)$ on $\cl\otimes W_z$. The quotient
\[\ti\bbe:=(\cl\otimes W_z)/\cg(E)\]
is a complex 2-plane bundle over $\cb^*(E)$.

We claim that there is an Hermitian
metric on $\ti\bbe$ such that on every fibre $\cl_A$ there is an Hermitian
metric for which the projection $\cl_A\otimes W_z\to\ti\bbe_{[A]}$ is an
isometry. To see this, let $S\subset\ca(E)$ be any local slice for the action
of $\cg(E)$, so that $S$ projects diffeomorphically onto an open subset
$U\subset\cb^*(E)$. Choose any Hermitian metric on $\cl|_S$ and let
$g_U$ be the induced Hermitian metric on $\ti\bbe_U\approx(\cl\otimes W_z)|_S$.
Now cover $\cb^*(E)$ by such open sets $U$ and patch together the corresponding
metrics $g_U$ to obtain the desired metric on $\ti\bbe$.

Given any Hermitian metric on a fibre $\cl_A$ there are linear isometries
\[\frg^0_{\cl_A\otimes W_z}\oset\approx\to\frg^0_{W_z}\oset\approx\to E_x,\]
where the first isometry is canonical and independent of the chosen metric
on $\cl_A$ and the second one is given by $\phi$. This yields an isomorphism
$\frg^0_{\ti\bbe}\oset\approx\to\bbe$.\square

\section{Moduli spaces}\label{section:mod-sp}

Let $P\to Y$ be a principal $\SO3$ bundle, where $Y$ is a closed oriented
$3$--manifold. The Chern-Simons functional
\[\cs:\ca(P)\to\R/\z\]
is determined up to an additive constant by the property that if $A$ is any
connection in the pull-back of $P$ to the band $[0,1]\times Y$ then
\begin{equation}\label{eqn:cs-int-band}
  \cs(A_1)-\cs(A_0)=\fct\int_{[t_0,t_1]\times Y}\la F_A\wedge F_A\ra,
  \end{equation}
where $A_t$ denotes the restriction of $A$ to the slice
$\{t\}\times Y$, and $\la\cdot\wedge\cdot\ra$ is formed by combining the wedge
product on forms with minus the Killing form on the Lie algebra of $\SO3$.
If $P=Y\times\SO3$ then we normalize $\cs$ so that its value on
the product connection $\theta$ is zero.  If $v$ is any automorphism of $P$
then for any connection $B$ in $P$ one has
\begin{equation}\label{eqn:csdeg}
  \cs(v(B))-\cs(B)=-\frac12\deg(v),
\end{equation}
where the degree $\deg(v)$ is defined to be the intersection number of $v$
with the image of the constant section $1$.

Equation~\refp{eqn:csdeg},
up to an overall sign, was stated without proof in \cite[Proposition~1.13]{BD1}.
A proof of \refp{eqn:csdeg}
can be obtained by first observing that the left-hand side of the
equation is independent of $B$, and both sides define homomorphisms from the
automorphism group of $P$ into $\R$. Replacing $v$ by $v^2$ it then
only remains to verify the equation
for even gauge transformations, which is easy.

If $v$ lifts to a section $\ti v$ of $P\uset{\SO3}\times\SU2$ then
\[\deg(v)=2\deg(\ti v),\]
where $\deg(\ti v)$ is the intersection number of $\ti v$ with the image
of the constant section $1$. In particular, every even automorphism of
$P$ has even degree.

The critical points of the Chern-Simons functional $\cs$
are the flat connections in $P$. In practice, we will add a small holonomy
perturbation to $\cs$ as in \cite{F1,D5}, but this will usually not be
reflected in our notation.
Let $\fl(P)$ denote the space of all critical points of $\cs$ modulo
even automorphisms of $P$. The even-reducible part of $\fl(P)$ is denoted
by $\fl^*(P)$. If $Y$ is an (integral) homology sphere then $P$ is
necessarily trivial and we write $\fl(Y)=\fl(P)$.

Now let $X$ be an oriented Riemannian $4$--manifold with tubular ends
$[0,\infty)\times Y_i$, $i=0,\dots,r$, such that the complement of
  \[\xend:=\bigcup_i\,[0,\infty)\times Y_i\]
is precompact. We review the standard set-up of moduli spaces of anti-self-dual
connections in a principal $\SO3$ bundle $Q\to X$, see \cite{D5}. Given a
flat connection $\rho$ in $Q|_{\xend}$, we define the moduli space
$M(X,Q;\rho)$ as follows. Choose a smooth connection $A_0$ in $Q$ which agrees
with $\rho$ outside a compact subset of $X$.
We use the connection $A_0$ to define Sobolev norms on forms with values
in the adoint bundle $\frg_Q$ of Lie algebras associated to $Q$.
Fix an even integer $p>4$.
Let $\ca=\ca(Q)$ be the space of connections in $Q$ of the form $A_0+a$ with
$a\in\llw pw1$, where $w$ is a small, positive exponential weight as in
\cite[Section~2.1]{Fr13}. There is a smooth action on $\ca$ by the Banach
Lie group $\cg$ consisting of all $\lloc p2$ gauge transformation $u$
of $Q$ such that $\nabla_{A_0}u\cdot u\inv\in\llw pw1$.
Let $\cb:=\ca/\cg$ and let $M(X,Q;\rho)$ be the subset of $\cb$ consisting
of gauge equivalence classes of connections $A$ satisfying $F^+_A=0$.
In practice, we will often add a small holonomy perturbation to the
ASD equation, but this will usually be suppressed from notation.

We observe that the value of the Chern-Simons integral
\begin{equation}\label{eqn:ka-int}
  \ka(Q,\rho):=-\frac1{8\pi^2}\int_X\la F_A\wedge F_A\ra
  \end{equation}
is the same for all $A\in\ca$. (If $X$ is closed then the right hand side
of Equation~\refp{eqn:ka-int} equals the value of $-p_1(Q)$ on the
fundamental class of $X$. This normalization will be convenient in
Section~\ref{sec:definite-4mflds}.)
If $u$ is an automorphism of $Q|_{\xend}$ then
from Equations~\refp{eqn:cs-int-band} and \refp{eqn:csdeg} we deduce that
\[\ka(Q,u(\rho))-\ka(Q,\rho)=2\sum_i\deg(u_i),\]
where $u_i$ is the restriction of $u$to the slice $\{0\}\times Y_i$.
Similarly, for the expected dimensions we have
\[\dim\,M(X,Q;u(\rho))-M(X,Q;\rho)=4\sum_i\deg(u_i).\]
On the other hand, if $u$ extends to a smooth automorphism of all of $Q$
then $\sum\deg(u_i)=0$,
and the converse holds at least if $u$ is even.

Given the reference connection $A_0$, we can identify the restriction of the
bundle $Q$ to an end $[0,\infty)\times Y_i$ with the pull-back of a bundle
$P_i\to Y_i$. 
Let $\al_i\in\fl(P_i)$ be the element obtained by restricting
$\rho$ to any slice $\{t\}\times Y_i$ where $t>0$. We will usually
assume that each $\al_i$ is non-degenerate.
The above remarks show that
the moduli space $M(X,Q;\rho)$ can be specified by the
$r$--tuple $\va=(\al_1,\dots,\al_r)$ together with one extra piece of data:
Either the Chern-Simons value $\ka=\ka(Q,\rho)$ or the expected dimension $d$
of $M(X,Q;\rho)$. We denote such a moduli space by 
\[M_{\ka}(X,Q;\va)\quad\text{or}\quad M_{(d)}(X,Q;\va).\]
Note that for given $\va$ there is exactly one moduli space $M_{(d)}(X,Q;\va)$
with $0\le d\le7$; this moduli space will just be denoted by $M(X,Q;\va)$.

For any anti-self-dual connection $A$ over $X$, the {\em energy} $\ce_A(Z)$
of $A$ over a measurable subset $Z\subset X$ is defined by
\begin{equation}\label{eqn:def-energy}
\ce_A(Z):=-\fct\int_Z\la F_A\wedge F_A\ra
=\fct\int_Z|F_A|^2.
\end{equation}
If $X=\ry$ and $Z=I\times Y$ for some interval $I$ then we write
$\ce_A(I)$ instead of $\ce_A(I\times Y)$.

\section{Spaces of linearly dependent vectors}\label{sec:spaces-lin-dep}

This section provides background for the definition of the cup product $u_2$
as well as results which will be used in the proof of
Proposition~\ref{prop:psiv3v2}.

The main result of this section is Proposition~\ref{prop:E-lla}.
To put this result into context, we first consider the Stiefel-Whitney classes
$w_j(E)$ of a real $n$--plane bundle $E\to B$ over a space $B$. For
$1\le k\le n$ let
$\ell\to\RP k$ be the tautological real line bundle.
Let $\uline E:=E\times\RP k$ and $\uline\ell:=B\times\ell$
be the pull-backs of
the bundles $E$ and $\ell$, respectively, to $B\times\RP k$.
If $B$ is paracompact, then by applying the Splitting Principle one finds that
$w_{n-k}(E)$ can be described as a slant product
\[w_{n-k}(E)=w_n(\uEue)/[\RP k],\]
where $[\RP k]$ denotes the fundamental class of $\RP k$ with
coefficients in $\z/2$. If $B$ is a closed manifold
then one can rephrase this formula in terms of Poincar\'e duals as follows:
The map on homology induced by the projection
\[p:B\times\RP k\to B\]
takes the Poincar\'e dual of $w_n(\uEue)$ to the Poincar\'e dual of $w_{n-k}(E)$,
i.e.
\begin{equation}\label{eqn:pdwn}
  p_*(\text{P.D.}(w_n(\uEue)))=\text{P.D.}(w_{n-k}(E)).
\end{equation}
Our main interest lies in the case when $n=3$, $k=1$, and $B$ has dimension
at most $5$. In this
case, the above formula can be deduced from Proposition~\ref{prop:E-lla} below,
as we will explain at the end of this section.

In order to state that proposition, we need some notation.
For any finite-dimensional real vector space $V$ set
\begin{equation}\label{eqn:LV}
L(V):=\{(v,w)\in V\oplus V\st\text{$v,w$ are linearly dependent in $V$}\}.
\end{equation}
Then $L(V)$ is closed in $V\oplus V$, and
\[L^*(V):=L(V)\setminus\{(0,0)\}\]
is a smooth submanifold of $V\oplus V$ of codimension~$n-1$, where $n$ is the
dimension of $V$.

As a short-hand notation we will often write $v\wedge w=0$ to express that
$v,w$ are linearly dependent.

For the remainder of this section we assume $B$ is a smooth Banach manifold and
$\pi:E\to B$ a smooth real vector
bundle of finite rank. Let $L^*(E)\to B$ be the associated smooth fibre bundle
whose fibre over a point $x\in B$ is $L^*(E_x)$, where $E_x=\pi\inv(x)$.
Similarly, let $L(E)\to B$ be the topological fibre bundle with fibre
$L(E_x)$ over $x$.

We now take $k=1$ in the discussion above and replace $\RP1$ with $S^1$.
Then $\ell\to S^1$ is the non-trivial real line bundle such that for
$z\in S^1$ the fibre of $\ell$ over $z^2$ is the line $\R z$ in $\bbc$.
We identify $\R^2=\bbc$, so that $(a,b)=a+bi$ for real numbers $a,b$.

\begin{prop}\label{prop:E-lla}
Suppose $s=(s_1,s_2)$ is a nowhere vanishing smooth section of $E\oplus E$.
Let $\si$ be the section of $\uline E\otimes\uline\ell$ such that for any
$p\in B$ and $z=(x_1,x_2)\in S^1$ one has
\[\si(p,z^2)=(x_1s_1(p)+x_2s_2(p))\otimes z.\]
\begin{description}
\item[(i)] The projection $B\times S^1\to B$ maps the zero-set of $\si$
  bijectively onto the locus in $B$ where $s_1$ and $s_2$
  are linearly dependent.
\item[(ii)] A zero $(p,w)$ of $\si$ is regular if and
  only if $s$ is transverse to $L^*(E)$ at $p$.
  \end{description}
\end{prop}

\proof
The proof of (i) is left as an exercise. To prove (ii) we may assume
$E$ is trivial, so that $s_j$ is represented by a smooth map $f_j:B\to V$
for some finite-dimensional real vector space $V$. We observe that
for any $u_1,u_2\in V$ and $z=(x_1,x_2)\in S^1$ one has
\begin{equation}\label{eqn:u1u2}
  (u_1,u_2)=(x_1u_1+x_2u_2)\otimes z+(x_1u_2-x_2u_1)\otimes iz
  \end{equation}
as elements of $V\oplus V=V\otimes_{\R}\bbc$.
It follows that the tangent space of $L^*(V)$ at a point $(v_1,v_2)$
which satisfies $x_1v_1+x_2v_2=0$ is given by
\begin{equation}\label{eqn:tlv}
  T_{(v_1,v_2)}L^*(V)=V\otimes iz+\R(x_1v_2-x_2v_1)\otimes z.
  \end{equation}
Now suppose $(p,w)$ is a zero of $\si$ and $s(p)=(v_1,v_2)$, $z^2=w$.
Then \refp{eqn:tlv} holds. Let $L_j:T_pB\to V$ be the derivative of $f_j$ at $p$.
Then $(p,w)$ is
a regular zero of $\si$ precisely when $V$ is spanned by the vector
$x_1v_2-x_2v_1$ together with the image of the map $x_1L_2+x_2L_2$.
From \refp{eqn:u1u2} we see that
the latter condition is also equivalent to $s$ being transverse to
$L^*(V)$ at $p$.\square

We record here a description of the sections of
$\uline E\otimes\uline\ell$ which
will be used in the proof of Proposition~\ref{prop:psiv3v2} below.
Let $\Ga_a(\uline E)$ denote
the space of all sections $s\in\Ga(\uline E)$ such that
\[s(p,-z)=-s(p,z)\]
for all $(p,z)\in B\times S^1$.

\begin{lemma}\label{lemma:Ga-a}
Then there is a canonical real linear isomorphism
\[\Ga(\uline E\otimes\uline\ell)\to\Ga_a(\uline E),\quad\si\mapsto\hat\si\]
characterized by the fact that
\[\si(p,z^2)=\hat\si(p,z)\otimes z\]
for all $(p,z)\in B\times S^1$.\square
\end{lemma}

We will now relate Proposition~\ref{prop:E-lla} to \refp{eqn:pdwn}, assuming
$B$ is finite-dimensional and closed.
First recall that the Poincar\'e dual of the top Stiefel-Whitney class
$w_n(E)$ is represented by the zero-set of a generic section of $E$.
Now suppose $E$ has rank $n=3$, and that $B$ has dimension at most $5$.
If $s$ is a generic smooth
section of $E\oplus E$ then $s$ does not vanish anywhere, and $s\inv(L(E))$
represents the Poincar\'e dual of $w_2(E)$.
Proposition~\ref{prop:E-lla} now yields \refp{eqn:pdwn} for $k=1$.

\section{``Generic'' sections}

Let $B$ be a smooth Banach manifold and $\pi:E\to B$ a smooth 
real vector bundle of finite rank. If $B$ is infinite-dimensional then we
do not define a topology on the space $\Ga(E)$ of (smooth) sections of
$E$, so it makes no sense to speak about residual subsets of $\Ga(E)$.
Instead, we will say
a subset $Z\subset\Ga(E)$ is ``residual'' (in quotation marks) if
there is a finite-dimensional subspace $\frP\subset\Ga(E)$ such that
for every finite-dimensional subspace $\frP'\subset\Ga(E)$ containing
$\frP$ and every section $s$ of $E$ there is a residual subset
$\frR\subset\frP'$ such that $s+\frR\subset Z$. Note that ``residual'' subsets
are non-empty, and any finite intersection of ``residual'' subsets is again
``residual''. We will say a given property
holds for a ``generic'' section of $E$ if it holds for every section belonging
to a ``residual'' subset of $\Ga(E)$.

We indicate one way of constructing such subspaces $\frP$.
Suppose $B$ supports smooth bump functions, i.e.\ for any point
$x\in B$ and any neighbourhood $U$ of $x$ there exists a smooth function
$c:B\to\R$ such that $c(x)\neq0$ and $c=0$ outside $U$. Given a compact subset
$K$ of $B$, one can easily construct a finite-dimensional subspace
$\frP\subset\Ga(E)$ such that, for every $x\in K$, the evaluation map
\[\frP\to E_x,\quad s\mapsto s(x)\]
is surjective. Therefore, if we are given a collection of smoooth maps
$f_k:M_k\to B$, $k=1,2,\dots$, where each $M_k$ is a
finite-dimensional manifold and the image of each $f_k$ is
contained in $K$ then, for a ``generic'' section $s$ of $E$, the map
\[s\circ f_k:M_k\to E\]
is transverse to the zero-section in $E$ for each $k$.

\section{Instanton cohomology and cup products}\label{sec:inst-cohom}

In this section we will work with $\SO3$ connections modulo
{\em even} gauge transformation (see Section~\ref{sec:base-point-fibration}),
although this {\em will not be
reflected in our notation}.
In particular, we write $\cb^*$ instead of $\cb^*\ev$. This notational
convention applies only to this section.
(In Subsection~\ref{subsec:comm-cup-products}, which only deals with
homology spheres, the convention is irrelevant.)

\subsection{Instanton cohomology}\label{subsec:inst-cohom}

Let $Y$ be a closed oriented connected Riemannian $3$-manifold and $P\to Y$ an 
$\SO3$ bundle. If $Y$ is not an homology sphere then we assume $P$ is
admissible. Let $\ry$ have the product Riemannian metric and
for any $\al,\beta\in\fl(P)$ let $M(\al,\beta)$ denote the
moduli space of instantons in the bundle $\R\times P\to\ry$ with flat
limits $\al$ at $-\infty$ and $\beta$ at $\infty$ and with expected dimension
in the interval $[0,7]$. Let
\[\chm(\al,\beta)=M(\al,\beta)/\R,\]
where $\R$ acts by translation. If $\al,\beta$ are irreducible then
the {\em relative index}
$\ind(\al,\beta)\in\z/8$ is defined by
\[\ind(\al,\beta)=\dim\,M(\al,\beta)\mod8.\]
For any commutative ring $R$ with unit we denote by $I(P;R)$ the relatively
$\z/8$ graded 
instanton cohomology with coefficients in $R$ as defined in
\cite{D5}. Recall that this is the cohomology of a cochain complex
$(C(P;R),d)$ where $C(P;R)$ is the free $R$--module generated by $\fl^*(P)$
and the differential $d$ is defined by
\[d\al=\sum_\beta\#\chm(\al,\beta)\cdot\beta.\]
Here, $\#$ means the number of points counted with sign,
and the sum is taken over all $\beta\in\fl^*(P)$ satisfying
$\ind(\al,\beta)=1$.
If $P$ is admissible then $\fl^*(P)=\fl(P)$. If instead $Y$ is an homology
sphere then $\fl(P)=\fl(Y)$ contains exactly one reducible point $\theta$,
represented by the trivial connection.
The presence of the trivial connection provides $C(P;R)=C(Y;R)$ with an absolute
$\z/8$ grading defined by
\[\ind(\al)=\dim\,M(\theta,\al)\mod8.\]
The trivial connection also gives rise to homomorphisms
\[C^4(Y;R)\oset{\del}\to R\oset{\del'}\to C^1(Y;R)\]
defined on generators by
\[\del\al=\#\chm(\al,\theta),\quad
\del1=\sum_\beta\#\chm(\theta,\beta)\cdot\beta,\]
where we sum over all $\beta\in\fl^*(Y)$ of index $1$.
These homomorphisms satisfy $\del d=0$ and $d\del'=0$ and therefore define
\[I^4(Y;R)\oset{\del_0}\to R\oset{\del_0'}\to I^1(Y;R).\]

We conclude this subsection with some notation for energy. If $A$ is any
ASD connection in the bundle $Q:=\R\times P$ and $I$ is any interval then
we write $\ce_A(I)$ instead of $\ce_A(I\times Y)$. Moreover, if
$\al,\beta\in\fl(Y)$ and the moduli space $M(\al,\beta)$ is expressed as
$M(\ry,Q;\rho)$ in the notation of Section~\ref{section:mod-sp} then
we define
\begin{equation}\label{eqn:cs-al-beta}
  \cs(\al,\beta):=\frac14\ka(Q,\rho),
  \end{equation}
which equals the total energy of any element of $M(\al,\beta)$. (Note,
however, that $M(\al,\beta)$ may be empty.)

\subsection{Cup products}\label{subsec:cup-prod}

We continue the discussion of the previous subsection, assuming $P$ is
admissible unless $Y$ is an homology sphere.
In most of this paper
the coefficient ring $R$ will be $\z/2$, and we write
\[I(P):=I(P;\z/2).\]
For $j=2,3$ we will define a degree~$j$ endomorphism
$u_j:I^*(P)\to I^{*+j}(P)$.
Insofar as the Floer cohomology is some kind of Morse
cohomology of $\cb^*(P)$, one may think of $u_j$ as cup product with
the $j$th Stiefel-Whitney class of the base-point fibration over $\cb^*(P)$.

The map $u_j$ will be
induced by an endomorphism
\[v_j:C^*(P)\to C^{*+j}(P)\]
which we now define. For any $t\in\R$ set
\begin{equation}\label{eqn:ybt-def}
  \yb t:=[t-1,t+1]\times Y.
\end{equation}
Let $P_0=[-1,1]\times P$ denote the pull-back of the bundle $P$ to $\yb0$.
For any $\al,\beta\in\fl(P)$ and any irreducible point $\om\in M(\al,\beta)$
let
\[\om[t]:=\om|_{Y[t]}\in\cb^*(P_0)\]
denote the restriction of $\om$ to the band $Y[t]$. (The fact that $\om[t]$
is irreducible follows from
Proposition~\refp{prop:unique-continuation-cylinder}.)
Choose a base-point $y_0\in Y$,
and let
\[\bbe\to\cb^*(P_0)\]
be the natural
real vector bundle of rank~$3$ associated to the base-point $(0,y_0)\in \yb0$.
To define $v_3$, choose a ``generic'' smooth section $s_1$ of $\bbe$.
For any $\al,\beta\in\fl^*(P)$
with $\ind(\beta)-\ind(\al)\equiv3\mod8$ the matrix coefficient
$\la v_3\al,\beta\ra$ is defined to be
\be{equation}\label{eqn:v3def}
\la v_3\al,\beta\ra:=\#\{\om\in M(\al,\beta)\st s_1(\om[0])=0\},
\end{equation}
where $\#$ means the number of points counted modulo $2$.
To define $v_2$, let $s_2,s_3$ be a pair of smooth sections of $\bbe$ which
defines a ``generic'' section of $\bbe\oplus\bbe$.
For any $\al,\beta\in\fl^*(P)$
with $\ind(\beta)-\ind(\al)\equiv2\mod8$ the matrix coefficient
$\la v_2\al,\beta\ra$ is defined to be
\begin{equation*}
\la v_2\al,\beta\ra:=
\#\{\om\in M(\al,\beta)\st\text{$s_2,s_3$ are linearly dependent at $\om[0]$}\}.
\end{equation*}
Note that, for dimensional reasons, neither $s_2$ nor $s_3$ will
vanish at $\om[0]$ for any $\om\in M(\al,\beta)$.

The following lemma gives an alternative description of $v_2$ in the spirit of
Section~\ref{sec:spaces-lin-dep}.

\begin{lemma}\label{lemma:v_2-defn-alt}
The matrix coefficient
$\la v_2\al,\beta\ra$ equals the number of points $(\om,z^2)$ in
$\cm(\al,\beta)\times S^1$, where $z=(a,b)\in S^1$, such that
\[as_2(\om[0])+bs_3(\om[0])=0.\square\]
\end{lemma}

\begin{prop}\label{prop:dvvd}
For $j=2,3$ one has
\[dv_j=v_jd\]
as homomorphisms $C^*(P)\to C^{*+j+1}(P)$.
\end{prop}

\proof To prove this for $j=2$, let $\al,\beta\in\fl^*(P)$ with
$\ind(\beta)-\ind(\al)\equiv3\mod8$.
The number of ends of the $1$-manifold
\[\{\om\in M(\al,\beta)\st
\text{$s_2,s_3$ are linearly dependent at $\om[0]$}\},\]
counted modulo $2$, is $\la(dv_2+v_2d)\al,\beta\ra$. Since the number of ends
must be even, this proves the assertion for $j=2$. The case $j=3$ is similar.
\square

The homomorphism $u_j:I^*(P)\to I^{*+j}(P)$ induced by $v_j$ is independent of
the sections $s_i$. For $u_3$ this will follow from
Lemma~\ref{lemma:v3-chain-hom} below,
and a similar argument works for $u_2$.
We consider again the bundle $P_0=[-1,1]\times P$ over $Y[0]=[-1,1]\times Y$.

\begin{defn}\label{defn:prop-T}
  Let $U$ be an open subset of $\cb^*(P_0)$ such that for all
  $\al,\beta\in\fl^*(P)$ with $\ind(\al,\beta)\le3$ and every
  $\om\in M(\al,\beta)$ one has that $\om[0]\in U$. A section $s$ of
  $\bbe|_U$ is said to {\em satisfy Property~$\pt3$} if for all $\al,\beta$
  as above the map
  \[M(\al,\beta)\to\bbe,\quad\om\mapsto s(\om[0])\]
  is transverse to the zero-section in $\bbe$.
  \end{defn}

\begin{lemma}\label{lemma:v3-chain-hom}
  Let $U\subset\cb^*(P_0)$ be as in Definition~\ref{defn:prop-T}
  and suppose $s,s'$ are sections of $\bbe|_U$ satisfying Property~$\pt3$.
  Let $v_3,v'_3$ be the corresponding cup products defined as in
  \refp{eqn:v3def}. Then there is an endomorphism
  \[H:C(P)\to C(P)\]
  such that
  \[v_3+v'_3=dH+Hd.\]
\end{lemma}

\proof For a ``generic'' section $\si$ of $\bbe$ the map
\begin{gather*}
  f_{\al\beta}:M(\al,\beta)\times[0,1]\to\bbe,\\
  \om\mapsto(1-t)s(\om[0])+ts'(\om[0])+t(1-t)\si(\om[0])
\end{gather*}
is transverse to the zero-section whenever $\ind(\al,\beta)\le3$.
Fix such a $\si$ and let $Z_{\al\beta}$ denote the zero-set of $f_{\al\beta}$.
If $\ind(\al,\beta)=2$ then $Z_{\al\beta}$ is a finite set. Let $H$ be the
homomorphism with matrix coefficients
\[\la H\al,\beta\ra=\#Z_{\al\beta}.\]
If $\ind(\al,\beta)=3$ then $Z_{\al\beta}$ is a compact
$1$--manifold-with-boundary. Counted modulo~$2$, the number of boundary points
of $Z_{\al\beta}$ is $\la(v_3+v'_3)\al,\beta\ra$, whereas the number of
ends is $\la(dH+Hd)\al,\beta\ra$. These two numbers must agree, proving
the lemma.\square

\begin{prop}\label{prop:u-del--funct}
  Let $W$ be a smooth, compact, oriented, connected $4$--manifold with two
  boundary components, say $\prtl W=-Y_0\cup Y_1$. Let $Q\to W$ be an
  $\SO3$ bundle, and let $P_i$ be the restriction of $Q$ to $Y_i$. Suppose
  one of the following two conditions holds.
  \begin{itemize}
  \item[(i)] At least one of the bundles $P_0,P_1$ is admissible.
  \item[(ii)] Both $Y_0$ and $Y_1$ are homology spheres, the bundle $Q$ is
    trivial, and $H_1(W;\z)=0$ and $b_+^2(W)=0$.
  \end{itemize}
  Then the homomorphism $T:I(P_0)\to I(P_1)$ induced by $(W,Q)$ satisfies
  \[Tu_j=u_jT\quad\text{for $j=2,3$.}\]
   Moreover, if (ii) holds then
  \[\del T=\del:I^4(Y_0)\to\z/2.\square\]
  \end{prop}

\begin{prop}\label{prop:u3zero}
If $P\to Y$ is an admissible $\SO3$ bundle then $u_3=0$ on $I(P)$.
\end{prop}

\proof By Proposition~\ref{prop:lift-to-U2} there is an Hermitian
$2$--plane bundle $\ti\bbe\to\bev^*$ such that $\bbe\approx\frg^0_{\ti\bbe}$.
For a ``generic'' section $\ti s$ of $\ti\bbe$, we have
$\ti s(\om[0])\neq0$ whenever $\om$ lies in a moduli space $M(\al,\beta)$
of dimension at most $3$. Given such a section $\ti s$, let $U$ be the
open subset of $\bev^*$ where $\ti s\neq0$. Then $\ti\bbe|_U$ splits
as an orthogonal sum
\[\ti\bbe|_U=\uline{\bbc}\oplus L\]
of two complex line bundles. Hence $\ti\bbe|_U$ has a nowhere vanishing
trace-free skew-Hermitian endomorphism
$\left(
\begin{array}{cc}
  i & 0 \\ 0 & -i
\end{array}
\right)$. This yields a non-vanishing section $s'$ of $\bbe|_U$.
Let $s$ be the restricion to $U$ of a ``generic'' section of $\bbe$,
and let $v_3,v'_3$ be the cup products defined by $s,s'$, respectively.
Then $v'_3=0$, so by Lemma~\ref{lemma:v3-chain-hom} we have
\[v_3=dH+Hd.\]
By definition, $v_3$ induces the cup product $u_3$ in cohomology,
so $u_3=0$.\square

\begin{prop}\label{prop:u3u2nzero}
  Let $Y$ be an oriented homology $3$--sphere and $Y'$ the result of $(\pm1)$
  surgery on a knot $\ga$ in $Y$. Let $n$ be a non-negative integer.
  \begin{description}
  \item[(i)] If $(u_3)^n=0$ on $I(Y)$ then $(u_3)^{n+1}=0$ on $I(Y')$.
  \item[(ii)] If $(u_2)^n=0$ on $I(Y)$ and $\ga$ has genus $1$
    then $(u_2)^{n+1}=0$ on $I(Y')$.
  \end{description}
\end{prop}

  \proof If $R$ is a commutative ring and
  \[A\longrightarrow B\longrightarrow C\]
  an exact sequence of modules over the polynomial ring $R[u]$ such
  that $u^m=0$ on $A$ and $u^n=0$ on $C$ for non-negative integers $m,n$ then
  $u^{m+n}=0$ on $B$. (Here, $u^0$ acts as the identity map.)

  Now suppose $Y'$ is $(-1)$ surgery on $\ga$. (If instead $Y'$ is $(+1)$
  surgery on $\ga$ then the proof is similar with the roles of $Y,Y'$ reversed.)
  Let $Y''$ be $0$ surgery on $\ga$ and $I(Y'')$ the instanton cohomology of
  the non-trivial $\SO3$ bundle over $Y''$. 
  We apply the above observation to the long exact surgery sequence
  (see \cite{BD1,Scaduto1})
  \[\cdots\to I(Y'')\to I(Y)\to I(Y')\to I(Y'')\to\cdots\]
  Statement~(i) now follows from Proposition~\ref{prop:u3zero}. To prove
  (ii), recall that if $P_{T^3}$ is a non-trivial $\SO3$ bundle over the
  $3$--torus then $I(P_{T^3})$ is non-zero in two degrees differing by
  $4$ modulo $8$ and zero in all other degrees. Therefore, $u_2=0$ on
  $I(P_{T^3})$.  If $\ga$ has genus~$1$ then by arguing as in the proof of
  \cite[Theorem~9]{Fr3} we find that $u_2=0$ on $I(Y'')$, from which
  (ii) follows.\square

  As a special case of Proposition~\ref{prop:u3u2nzero} we have the following
  corollary.
  
\begin{cor}\label{cor-u3-vanishes-on-surgery}
If $Y$ is $(\pm1)$ surgery on a knot in $S^3$ then $u_3=0$ on $I(Y)$.
  \end{cor}

\begin{prop} Let $P\to Y$ be an $\SO3$ bundle. We assume $P$ is admissible
  if $Y$ is not a homology sphere. Then the endomorphisms $u_2$ and $u_3$
  on $I(P)$
  are nilpotent. In other words, there is a positive integer $n$ such that
  \[u_2^n=0,\quad u_3^n=0\quad\text{on $I(P)$.}\]
  \end{prop}

\proof We use the same link reduction schemes as
in the proofs of \cite[Theorems~9 and 10]{Fr3}. 
In the present case there is no need to consider any reduced groups, as
the cup products $u_j$ are defined on all of $I(Y)$.\square

We include here a result for oriented homology $3$--spheres $Y$ obtained by
adapting the proof of Proposition~\ref{prop:dvvd} for $j=2$ to $2$--dimensional
moduli spaces $M(\al,\theta)$. This result will be used in
Proposition~\ref{prop:dPPQ} below.
For any $\ga\in\fl^*(Y)$ we introduce the
temporary notation
\[M_\ga:=\{\om\in M(\ga,\theta)\st
\text{$s_2\wedge s_3=0$ at $\om[0]$, and $\ce_\om([0,\infty))\ge\eps$}\},\]
where $\eps$ is a small positive constant.
If $\dim\,M(\ga,\theta)<6$ then $M_\ga$ is a manifold-with-boundary, and
$\prtl M_\ga$ has a description analogous to that of $M_\ga$, just replacing
the inequality $\ce_\om([0,\infty))\ge\eps$
  by an equality. We define homomorphisms
\[\dot\del:C^2(Y)\to\z/2,\quad\del^-:C^3(Y)\to\z/2\]
on generators by
\[\dot\del\al:=\#(\prtl M_\al),\quad \del^-\beta:=\# M_\beta.\]

\begin{prop}\label{prop:delv2deld}
$\del v_2+\del^-d=\dot\del$.
  \end{prop}
\proof Let $\al\in\fl^*(Y)$, $\ind(\al)=2$. Then $M_\al$ is a
  $1$--manifold-with-boundary. The number of boundary points,
  counted modulo~$2$, is $\dot\del\al$
  by definition, and this must agree with the number of ends of $M_\al$, which
  is $(\del v_2+\del^-d)\al$.\square

\subsection{Commutators of cup products}\label{subsec:comm-cup-products}

Let $Y$ be an oriented homology $3$--sphere.
We introduce a degree~$4$ endomorphism
\[\phi:C^*(Y)\to C^{*+4}(Y)\]
which will be used to describe the commutator of $v_2$ and $v_3$.

\be{defn}\label{defn:M23}
For any $\al,\beta\in\fl^*(Y)$
let $\mx23(\al,\beta)$ be the subspace of $\mab\times\R$ consisting of those
points $(\om,t)$ satisfying the following conditions:
\be{itemize}
\item $s_1(\om[-t])=0$,
\item $s_2(\om[t])$ and $s_3(\om[t])$ are linearly dependent.
\end{itemize}
\end{defn}

If $\ind(\beta)-\ind(\al)\equiv4\mod8$
then $\mx23(\al,\beta)$ consists of a finite number of points (see part~(I)
of the proof of Proposition~\ref{prop:v2v3phi} below), and we set
\[\la \phi\al,\beta\ra:=\#\mx23(\al,\beta).\]

\be{prop}\label{prop:v2v3phi}
If $Y$ is an oriented integral homology $3$-sphere then for ``generic''
sections $s_1,s_2,s_3$ one has
\be{equation}\label{eqn:v2v3chhom}
v_2v_3+v_3v_2+\del'\del=d\phi+\phi d.
\end{equation}
Hence, on $I(Y)$ one has
\be{equation}\label{eqn:u2u3}
u_2u_3+u_3u_2=\del_0'\del_0.
\end{equation}
\end{prop}

The proof will be given in Subsection~\ref{subsec:proof-v2v3phi}.

Let $v_3,v_3':C^*(Y)\to C^{*+3}(Y)$ be the cup products defined by ``generic''
sections $s,s'$ of $\bbe$. At least in degrees different from $3$ and $4$, the
commutator of $v_3$ and $v_3'$ is given by a formula analogous to
\refp{eqn:v2v3chhom}. This formula involves the homomorphism
\[\psi:C^p(Y)\to C^{p+5}(Y),\quad p\neq4\]
with matrix coefficients
\[\la\psi\al,\beta\ra=\#\{(\om,t)\in\mab\times\R\st s(\om[-t])=0=s'(\om[t])\}.\]
The condition $p\neq4$ is imposed to make sure that factorizations through
the trivial connection do not occur in the moduli spaces $M(\al,\beta)$.

\begin{prop}\label{prop:v3v3}
For $q\not\equiv3,4\mod8$ one has
\begin{equation}\label{eqn:dpsipsid}
  d\psi+\psi d=v_3v'_3+v'_3v_3
\end{equation}
as maps $C^q(Y)\to C^{q+6}(Y)$.
\end{prop}

The proof is given in Subsection~\ref{subsec:cochain-map-psi},
where the proposition
is restated as Proposition~\ref{prop:v3v3-restated}.

As Tom Mrowka pointed out to the author, Equation~\refp{eqn:dpsipsid}
is reminiscent of the cup-$i$ construction of Steenrod squares, see for
instance \cite[p.\,271]{Spanier}.

If the sections $s,s'$ are sufficiently close (in a certain
sense) then $v_3=v_3'$ (see Lemma~\ref{lemma:prop-T} below)
and the following hold.

\begin{prop}\label{prop:psixi}
  If the sections $s,s'$ are sufficiently close then there exist
  \begin{itemize}
  \item an extension of $\psi$ to a cochain map $C^*(Y)\to C^{*+5}(Y)$
    defined in all degrees, and
  \item a homomorphism $\Xi:C^*(Y)\to C^{*+4}(Y)$ such that
    \[\psi=v_2v_3+d\Xi+\Xi d,\]
    where the cup products $v_2,v_3$ are defined by three ``generic'' sections
    of $\bbe$.
  \end{itemize}
  \end{prop}

The proof will be given in Subsection~\ref{subsec:calc-psi}.

\section{Definition of the invariant $\qt$}\label{section:defqt}

Let $Y$ be any oriented homology $3$-sphere.

\be{defn}\label{defn:zetai}
We define a non-negative integer $\zeta_2(Y)$ as follows. If $\del_0=0$
on $\ker(u_3)\subset I(Y)$ set $\zeta_2(Y):=0$. Otherwise, let $\zeta_2(Y)$
be the largest positive integer $n$ for which there exists an
$x\in\ker(u_3)$ such that
\[\del_0u_2^kx=\begin{cases}
0 & \text{for $0\le k<n-1$,}\\
1 & \text{for $k=n-1$.}
\end{cases}\]
\end{defn}
Here, $u_2^k$ denotes the $k$'th power of the endomorphism $u_2$. Note that
if $x$ is as in Definition~\ref{defn:zetai} then using
the relation~\refp{eqn:u2u3} one finds that $u_3u_2^kx=0$ for
$0\le k\le n-1$.

\be{defn}Set $\qt(Y):=\zeta_2(Y)-\zeta_2(-Y)$.
\end{defn}

An alternative description of $\qt$ will be given in
Proposition~\ref{prop:alt-def-qt} below. 

\begin{lemma}\label{lemma:q2ymin}
  If $\im(\del'_0)\subset\im(u_3)$ in $I^1(Y)$ then $\zeta_2(-Y)=0$. Otherwise,
  $\zeta_2(-Y)$ is the largest positive integer $n$ for which the inclusion
  \begin{equation}\label{eqn:imu2delincl}
  \im(u_2^k\del'_0)\subset\im(u_3)+\sum_{j=0}^{k-1}\im(u_2^j\del'_0)
  \quad\text{in $I(-Y)$}
  \end{equation}
  holds for $0\le k<n-1$ but not for $k=n-1$.
\end{lemma}
Of course, in \refp{eqn:imu2delincl} it suffices to sum over those $j$ that are
congruent to $k$ mod~$4$, since $I(-Y)$ is mod~$8$ periodic.

\proof Recall that $I^q(Y)$ and $I^{5-q}(-Y)$ are dual vector spaces for any
$q\in\z/8$. Furthermore, the maps
\[\del_0:I^4(Y)\to\z/2,\quad u_3:I^q(Y)\to I^{q+j}(Y)\]
are dual to
\[\del'_0:\z/2\to I^1(-Y),\quad u_3:I^{5-q-j}(-Y)\to I^{5-q}(-Y),\]
respectively. In general, the kernel of a linear map between finite-dimensional
vector spaces is equal to the annihilator of the image of the dual map.
Applying this to $\del_0u_2^j:I^{4-2j}(Y)\to\z/2$ we see that the inclusion
\refp{eqn:imu2delincl} holds if and only if
\[\ker(\del_0u_2^k)\supset\ker(u_3)\cap\bigcap_{j=0}^{k-1}\ker(\del_0u_2^j)
\quad\text{in $I(Y)$.}\]
This proves the lemma.\square

\be{prop}
Either $\zeta_2(Y)=0$ or $\zeta_2(-Y)=0$.
\end{prop}

\proof Suppose $\zeta_2(Y)>0$, so there is an $x\in I^4(Y)$ such that
$u_3x=0$ and $\del_0x=1$. Then Proposition~\ref{prop:v2v3phi} yields
$\del'_0(1)=u_3u_2x$, hence $\zeta(-Y)=0$ by Lemma~\ref{lemma:q2ymin}.\square

We now reformulate the definition of $\zeta_2$ in terms of the mapping cone of
$v_3$. This alternative definition will display a clear analogy with the 
instanton $h$-invariant and will be essential for handling the algebra involved
in the proof of additivity of $\qt$.
For $q\in\z/8$ set
\[MC^q(Y):=C^{q-2}(Y)\oplus C^q(Y),\]
and define
\[D:MC^q(Y)\to MC^{q+1}(Y),\quad(x,y)\mapsto(dx,v_3x+dy).\]
Then $D\circ D=0$, and we define $MI(Y)$ to be the cohomology of the
cochain complex $(MC(Y),D)$. The short exact sequence of cochain
complexes
\[0\to C^*(Y)\oset\si\to MC^*(Y)\oset\tau\to C^{*-2}(Y)\to0,\]
where $\si(y)=(0,y)$ and $\tau(x,y)=x$,
gives rise to a long exact sequence
\be{equation}\label{eqn:ujexact}
\cdots\to I^{q-3}(Y)\oset{u_3}\to I^q(Y)\oset{\si_*}
\to MI^q(Y)\oset{\tau_*}\to I^{q-2}(Y)\to\cdots.
\end{equation}
We introduce some extra structure on $\ix*j(Y)$. Firstly,
the homomorphisms
\be{gather*}
\Del:=\del\circ\tau:MC^6(Y)\to\z/2,\\
\Del':=\si\circ\del':\z/2\to MC^1(Y)
\end{gather*}
induce homomorphisms
\[MI^6(Y)\oset{\Del_0}\longrightarrow\z/2\oset{\Del'_0}\longrightarrow
MI^1(Y).\]
We extend $\Del$ trivially to all of $MC(Y)$, and similarly for $\Del_0$.
Furthermore, we define a homomorphism
\[V:MC^*(Y)\to MC^{*+2}(Y),\quad(x,y)\mapsto(v_2x,\phi x+v_2y).\]
A simple calculation yields
\be{equation}\label{eqn:dvvdx}
DV+VD=\Del'\Del,
\end{equation}
which is analogous to the relation \cite[Theorem~4~(ii)]{Fr3} in rational
instanton homology. It follows that $V$ induces homomorphisms
\be{gather*}
MI^q(Y)\to MI^{q+2}(Y),\quad q\not\equiv6,7\mod8,\\
MI^6(Y)\cap\ker(\Del_0)\to MI^0(Y),
\end{gather*}
each of which will be denoted by $U$.

\begin{prop}\label{prop:alt-def-qt}
If $\Del_0=0$ on $MI^6(Y)$ then $\zeta_2(Y)=0$. Otherwise,
$\zeta_2(Y)$ is the
largest positive integer $n$ for which there exists a $z\in MI(Y)$
such that 
\[\Del_0 U^kz=\be{cases}
0 & \text{for $0\le k<n-1$,}\\
1 & \text{for $k=n-1$.}
\end{cases}\]
\end{prop}

\proof This follows immediately from the definitions.\square


\section{Definite $4$-manifolds}\label{sec:definite-4mflds}

The goal of this section is to prove Theorem~\ref{thm:defub}.
Let $X$ be an oriented, connected
Riemannian $4$--manifold with a cylindrical end $[0,\infty)\times Y$,
where $Y$ is an integral homology sphere. 
Suppose
\[b_1(X)=0=b^+(X).\]
Let $E\to X$ be an oriented Euclidean $3$--plane bundle and $w_2(E)$
its second Stiefel-Whitney class. We will count reducibles in
ASD moduli spaces for $E$
with trivial asymptotic limit.

Let $\ti w\in H^2(X,\xend;\z/2)$ be the unique
lift of $w_2(E)$. Abusing notation, we denote by $w_2(E)^2\in\z/4$
the value of the Pontryagin square
\[\ti w\vphantom{w}^2\in H^4(X,\xend;\z/4)\]
on the fundamental class in
$H_4(X;\xend;\z/4)$.  Then for $\al\in\fl^*(Y)$ the expected dimension of
a moduli space for $E$ with asymptotic limit $\al$ satisfies
\[\dim\,M_\ka(X,E;\al)\equiv\ind(\al)-2w_2(E)^2\mod8.\]

If $\rho$ is a trivial connection in $E|_{\xend}$ then $\ka(E,\rho)$ is an
integer reducing to $-w_2(E)^2$ modulo~$4$. Hence,
\[M_k:=M_k(X,E;\theta)\]
is defined for integers $k$ satisfying $k\equiv-w_2(E)^2\mod4$. Moreover,
$M_k$ is empty for $k<0$, and $M_0$ (when defined) consists of flat connections.
The expected dimension is
\[\dim\,M_k=2k-3.\]

\subsection{Reducibles}

In this subsection we restrict to $k>0$.
After perturbing the Riemannian metric on $X$ in a small ball we can arrange
that $M_k$ contains no twisted reducibles (see \cite{Fr14}).

The set $M\red_k$ of reducible (i.e.\ Abelian)
points in $M_k$ has a well known description
in terms of the cohomology of $X$, which we now recall. Let
\[\ti P:=\{c\in H^2(X;\z)\st [c]_2=w_2(E),\;c^2=-k\},\]
where $[c]_2$ denotes the image of $c$ in $H^2(X;\z/2)$.
Let $P:=\ti P/\pm1$ be the quotient of $\ti P$ by the involution
$c\mapsto-c$.
 
\begin{prop}
There is a canonical bijection $M\red_k\to P$.
  \end{prop}

\proof If $[A]\in M\red_k$ then $A$ respects a unique splitting
\[E=\lla\oplus L,\]
where $\lla$ is a trivial rank~$1$ subbundle of $E$. A choice of orientation
of $\lla$ defines a complex structure on $L$. Mapping $[A]$ to the point in
$P$ represented by $c_1(L)$ yields the desired bijection. For further
details see \cite[Lemma~2]{Fr3} and \cite[Proposition~4.1]{Fr14}.\square

Assuming $P$ is non-empty we now express the number $|P|$ of elements of $P$
in terms of the intersection form of $X$ and the torsion subgroup
$\ct$ of $H^2(X;\z)$. For any $v\in H^2(X;\z)$ let $\bar v$ denote
the image of $v$ in $H^2(X;\z)/\ct$. Choose $a\in\ti P$ and let
\[\ti Q_a:=\{r\in H^2(X;\z)/\ct\st r\equiv\bar a\;\text{mod}\;2,\;r^2=-k\}.\]
Define $Q_a:=\ti Q_a/\pm1$.

\begin{prop}\label{prop:PtauQa}
$|P|=|2\ct|\cdot|Q_a|$.
  \end{prop}
Note that $2\ct$ has even order precisely when $H^2(X;\z)$ contains an element
of order $4$.

\proof Because $k>0$ we have that $(-1)$ acts without fixed-points on both
$\ti P$ and $\ti Q_a$. Therefore,
\begin{equation}\label{eqn:2PQ}
  |\ti P|=2|P|,\quad|\ti Q_a|=2|Q_a|.
  \end{equation}
The short exact sequence $0\to\z\oset2\to\z\to\z/2\to0$ gives rise to a long
exact sequence
\begin{equation}\label{eqn:2long-exact-seq}
  \cdots\to H^2(X;\z)\oset2\to H^2(X;\z)\to H^2(X;\z/2)\to H^3(X;\z)\to\cdots.
\end{equation}
From this sequence we see that there is a well defined map
\[\ti P\to\ti Q_a,\quad c\mapsto\bar c\]
which descends to an injective map
\[f:\ti P/2\ct\to\ti Q_a.\]
In fact, $f$ is bijective. To see that $f$ is surjective, let $r\in\ti Q_a$.
Then
\[r=\bar a+2\bar x=\oline{a+2x}\]
for some $x\in H^2(X;\z)$, and $a+2x\in\ti P$. This shows that
\[|\ti P|=|2\ct|\cdot|\ti Q_a|.\]
Combining this with \refp{eqn:2PQ} we obtain the proposition.\square

\subsection{$2$--torsion invariants of $4$--manifolds}
\label{subsec:2-torsion-inv}

The proof of Theorem~\ref{thm:defub} will involve certain $2$--torsion
Donaldson invariants which we now define. Let $d_0$ be the smallest expected
dimension of any moduli space $M_k=M_k(X,E;\theta)$ that contains a reducible,
where $k$ is a non-negative integer.
For any pair
$(r,s)$ of non-negative integers satisfying
\[2r+3s\le d_0+2\]
we will define an element
\[\twod rs=\twod rs(X,E)\in I(Y)\]
which will be independent of the Riemannian
metric on $X$ and also independent of the choice of small holonomy
perturbations.

To define $\twod rs$, choose disjoint compact codimension~$0$ submanifolds
$Z_1,\dots,Z_{r+s}$ of $X$ and base-points $z_j\in Z_j$.
It is convenient to assume that each of these submanifolds contains a band
$[t_j,t_j+1]\times Y$ for some $t_j\ge1$. (We assume that the perturbed
ASD equation is of gradient flow type in the region $[1,\infty)\times Y$.)
Then
Proposition~\ref{prop:unique-continuation-cylinder} guarantees that
every perturbed ASD connection in $E$ with irreducible limit will
restrict to an irreducible connection over each $Z_j$.

Choose ``generic''
sections $\{\si_{ij}\}_{i=1,2,3}$ of the canonical $3$--plane bundle
$\bbe_j\to\cb^*(Z_j,E_j)$, where $E_j:=E|_{Z_j}$. For any $\al\in\fl^*(Y)$
let $d=d(\al)$ be the integer such that
\begin{gather*}
  0\le d-2r-3s\le7,\\
  d\equiv\ind(\al)-2w_2(E)^2\mod8.
\end{gather*}
Let $M_{r,s}(X,E;\al)$ be the set of all $\om\in M_{(d)}(X,E;\al)$ such that
\begin{itemize}
\item $\si_{2,j},\si_{3,j}$ are linearly dependent at $\om|_{Z_j}$ for
  $j=1,\dots,r$, and
\item $\si_{1,j}(\om|_{Z_j})=0$ for $j=r+1,\dots,r+s$.
\end{itemize}
Let
\begin{equation}\label{eqn:qdef}
  q_{r,s}:=\sum_\al\#M_{r,s}(X,E;\al)\cdot\al\in C(Y),
  \end{equation}
where the sum is taken over all generators in $C(Y)$ of index
$2w_2(E)^2+2r+3s$. Then $q_{r,s}$ is a cocycle, and we define
\[\twod rs(X,E):=[q_{r,s}]\in I(Y).\]
Standard arguments show that $\twod rs$ is
independent of the choice of submanifolds
$Z_j$ and sections $\si_{ij}$.

\begin{prop}\label{prop:delDred}
  Let $k$ be an integer greater than one.
  If $M\red_\ell$ is empty for $\ell<k$ then
  \[\del\twod{k-2}0=\#M\red_k.\]
  \end{prop}

\proof Deleting from $M_k$ a small neighbourhood of each reducible point
we obtain a manifold-with-boundary $W$ with one boundary component $P_\eta$
for each reducible $\eta$, each such component being diffeomorphic to
$\CP{k-2}$. Let
\[\hat W:=W\cap M_{k-2,0}(X,E;\theta)\]
be the set of all $\om\in W$ such that $\si_{2,j}$ and $\si_{3,j}$ are linearly
dependent at $\om|_{Z_j}$ for $j=1,\dots,k-2$.
Then $\hat W$ is a $1$--manifold-with-boundary. For dimensional reasons and
because of the condition that $M\red_\ell$ be empty for $\ell<k$, bubbling
cannot occur in sequences in $\hat W$. Therefore, the only source of
non-compactness in $\hat W$ is factorization over the end of $X$, so
the number of ends of $\hat W$ equals $\del\twod{k-2}0$ modulo $2$.
As for the boundary points of $\hat W$,
observe that for every $x\in X$ the restriction of the $3$--plane bundle
$\bbe_{\theta,x}\to M^*_k$ to $P_\eta$ is isomorphic to the direct sum
$\uline{\R}\oplus L$ of a trivial real line bundle and the tautological
complex line bundle. It follows easily from this that $P_\eta\cap\hat W$
has an odd number of points for every reducible $\eta$, hence
\[|\prtl\hat W|\equiv|M\red_k|\mod2.\]
Since the number of boundary points of $\hat W$ must agree with the number of
ends when counted modulo $2$, this proves the proposition.\square

In the proof of the following proposition and at many places later we will
make use of a certain kind of cut-off function. This should be a smooth
function $b:\R\to\R$ such that
\begin{equation}\label{eqn:b-prop1}
  b(t)=
  \begin{cases}
0&\text{for $t\le-1$,}\\
1&\text{for $t\ge1$.}
\end{cases}
\end{equation}

\begin{prop}\label{prop:slide-points}
  Suppose $2r+3s\le d_0+2$, so that $\twod rs$ is defined.
  \begin{description}
  \item[(i)] $\twod rs=u_2\twod{r-1}s$ if $r\ge1$.
  \item[(ii)] $\twod rs=u_3\twod r{s-1}$ if $s\ge1$.
    \end{description}
  \end{prop}

\proof We only spell out the proof of (ii), the proof of (i) being similar.
Let $M_{r,s-1}(X,E;\al)$ be defined as above, but using only the submanifolds
$Z_1,\dots,Z_{r+s-1}$ and the corresponding sections $\si_{ij}$.
Choose a path $\ga:[-1,\infty)\to X$ such that $\ga(-1)=z_{r+1}$ and
$\ga(t)=(t,y_0)$ for $t\ge0$, where $y_0\in Y$ is a base-point.
  For any $\al\in\fl^*(Y)$ and $x\in X$ let
  \[\bbe_{\al,x}\to M_{r,s-1}(X,E;\al)\]
  be the canonical $3$--plane bundle associated to the base-point $x$.
  For any $\om=[A]\in M_{r,s-1}(X,E;\al)$ and $t\ge-1$ let
  \[\hol_{\om,t}:(\bbe_{\al,\ga(t)})_\om\to(\bbe_{\al,\ga(-1)})_\om\]
  be the isomorphism defined by the holonomy of $A$ along $\ga$.
  Here, $(\bbe_{\al,x})_\om$ denotes the fibre of the bundle $\bbe_{\al,x}$ at
  the point $\om$.
  Given a ``generic'' section $s$ of $\bbe\to\cb^*(Y[0])$ we define
  a section $s_\al$ of the bundle
  \[\bbe_{\al,\ga(-1)}\times[-1,\infty)\to M_{r,s-1}(X,E;\al)\times[-1,\infty)\]
by
\[s_\al(\om,t):=(1-b(t-2))\cdot\si_{1,r+s}(\om|_{Z_{r+s}})
+b(t-2)\cdot\hol_{\om,t}(s(\om[t])),\]
where $b$ is as in \refp{eqn:b-prop1}.
Let $j:=2w_2(E)^2+2r+3s\in\z/8$. If $\ind(\al)=j-1$ then the zero set
$s_\al\inv(0)$ is a finite set. Summing over such $\al$ we define
\[h_{r,s}:=\sum_\al(\#s_\al\inv(0))\cdot\al\in I^j(Y).\]
Counting ends and boundary points of the $1$--manifolds $s_\beta\inv(0)$
for $\ind(\beta)=j$ we see that
\[dh_{r,s}+v_3q_{r,s-1}=q_{r,s}.\]
Passing to cohomology, we obtain (ii).\square

\begin{prop}\label{prop:ESnon-trivial}
  If $E$ is strongly admissible then $D_{r,s}(X,E)=0$ for $s>0$.
  \end{prop}

\proof Let $f:\Si\to X$ be as in Definition~\ref{defn:admissible}
with $v=w_2(E)$. For $t\ge0$ let $X\sco t$ be the
result of deleting from $X$ the open subset $(t,\infty)\times Y$.
Choose $t>0$ so large that $X\sco t$ contains $f(\Si)$. Then
$E|_{X\sco t}$ is strongly admissible.
Choose the submanifolds $Z_1,\dots,Z_{r+s}$ such that $Z_{r+s}=X\sco t$.
By Proposition~\ref{prop:lift-to-U2}
the (frame bundle of) $\bbe_j\to\cb^*(E_{r+s})$ 
lifts to a $\U2$ bundle.
For $j=1,\dots,r+s-1$ choose ``generic'' sections $\{\si_{ij}\}_{i=1,2,3}$
of $\bbe_j$. Arguing as in the proof of Proposition~\ref{prop:u3zero}
we see that there is an open subset $U\subset\cb^*(Z_{r+s},E_{r+s})$
and a section $\si$ of $\bbe_{r+s}$ such that if $\om$ is any element of
a $3$--dimensional moduli space $M_{r,s-1}(X,E;\al)$ then $\om|_{Z_{r+s}}\in U$
and $\si(\om|_{Z_{r+s}})\neq0$. Taking $\si_{1,r+s}:=\si$ we have that all
$0$--dimensional moduli spaces $M_{r,s}(X,E;\al)$ are empty. 
Reasoning as in the proof of Lemma~\ref{lemma:v3-chain-hom} we conclude
that $D_{r,s}=0$.\square

\subsection{Lower bound on $\qt$}

Recall Definition~\ref{defn:admissible} above.

\begin{defn} Given a space, $X$, a non-zero class $w\in H^2(X;\z)/torsion$
  is called {\em strongly admissible}
  if some (hence every) lift of $w$ to $H^2(X;\z)$ maps to a strongly
  admissible class in $H^2(X;\z/2)$.
\end{defn}

\begin{thm}\label{thm:defub-gen}
Let $V$ be a smooth compact oriented connected $4$-manifold whose boundary
is a homology sphere $Y$ and whose intersection form is negative
definite. Let $w$ be an element of
\[J_V:=H^2(V;\z)/torsion\]
which is not divisible by $2$ and suppose
at least one of the following two conditions holds:
\begin{itemize}
\item[(i)] $H^2(V;\z)$ contains no $2$--torsion.
\item[(ii)] $H^2(V;\z)$ contains no element of order $4$, and
  $w^2\not\equiv0\mod4$. Furthermore, either $w$ is strongly admissible or
  $u_3=0$ on $I(Y)$ (or both).
\end{itemize}
Let $k$ be the minimal square norm (with
respect to the intersection form) of any element
of $w+2J_V$. Let $n$ be the number of elements of $w+2J_V$ of square norm $k$.
If $k\ge2$ and $n/2$ is odd then
\be{equation}\label{eqn:q2ineq}
\qt(Y)\ge k-1.
\end{equation}
\end{thm}

Note that if we leave out case~(ii) then the theorem says the same as
Theorem~\ref{thm:defub}.

\proof After performing surgery on a collection of loops in $V$ representing
a basis for $H_1(V;\z)/\torsion$ we may assume that $b_1(V)=0$.
From the exact sequence \refp{eqn:2long-exact-seq} we see that the
$2$--torsion subgroup of $H^2(V;\z)$ is isomorphic to $H^1(V;\z/2)$.
Let
\[X:=V\cup(0,\infty)\times Y\]
be the result of adding a half-infinite cylinder to $V$, and choose a
Riemannian metric on $X$ which is of cylindrical form over the end.
We identify the (co)homology of $X$ with that of $V$. Choose a
complex line bundle $L\to X$ whose Chern class represents $w$. Choose a
Euclidean metric on the $3$--plane bundle
\[E:=\uline{\R}\oplus L.\]
Since we assume that $H^2(X;\z)$
contains no element of order $4$, it follows from Proposition~\ref{prop:PtauQa}
that $M_\ell$ contains an odd number of reducibles for $\ell=k$
but no reducibles for $0<\ell<k$.

We now show that if $w^2\equiv0\,(4)$, so that $M_0$ is defined, then $M_0$
is free of reducibles. Suppose $A$ is a connection in $E$
representing a reducible point in $M_0$. Then $A$ preserves some orthogonal
splitting $E=\lla\oplus L'$, where $\lla\to X$ is a real line bundle.
Because Condition~(i) of the proposition must hold, the bundle $\lla$ is
trivial. Choose a complex structure on $L'$. Since $L'$ admits a flat
connection, its Chern class $c_1(L')$ is a torsion class in $H^2(X;\z)$.
But $c_1(L)$ and $c_1(L')$ map to the same element of $H^2(X;\z/2)$, namely
$w_2(E)$, hence
\[c_1(L)=c_1(L')+2a\]
for some $a\in H^2(X;\z)$. This contradicts our assumption that $w\in J_V$
is not divible by $2$. Thus, $M_0$ is free of reducibles as claimed.

By Proposition~\ref{prop:delDred} we have
\[\del D_{k-2,0}\neq0,\]
and Proposition~\ref{prop:slide-points} says that
\[D_{k-2,0}=u_2^{k-2}D_{0,0}.\]
Now suppose $w$ is strongly admissible (which is trivially the case
if Condition~(i) holds). Then
the bundle $E$ is strongly admissible, so by
Propositions~\ref{prop:slide-points} and \ref{prop:ESnon-trivial} we have
\[u_3D_{0,0}=D_{0,1}=0.\]
This proves \refp{eqn:q2ineq}.\square

\section{Operations defined by cobordisms}\label{sec:operations}

\subsection{Cutting down moduli spaces}\label{subsec:cutting-down-mod-sp}

Let $Y_0,Y_1,Y_2$ be oriented (integral) homology $3$--spheres and $W$ a
smooth compact connected oriented $4$--manifold such that
$H_i(W;\z)=0$ for $i=1,2$ and $\prtl W=(-Y_0)\cup(-Y_1)\cup Y_2$. Then we call
$W$ a ($4$--dimensional) {\em pair-of-pants cobordism} from $Y_0\cup Y_1$ to
$Y_2$, or a pair-of-pants cobordism from $Y_1$ to $(-Y_0)\cup Y_2$.

We will consider various operations on Floer cochain complexes induced by
pair-of-pants cobordism. To define these we first introduce some notation.

Let $X$ be an oriented connected Riemannian $4$--manifold with incoming
tubular ends $(-\infty,0]\times Y_j$, $j=0,\dots,r$ and outgoing tubular ends
$[0,\infty)\times Y_j$, $j=r+1,\dots,r'$, where each $Y_j$ is an
homology sphere. For $t\ge0$ let $X\sco t$ be the result of deleting
from $X$ the open pieces $(-\infty,-t)\times Y_j$, $j=0,\dots,r$ and
$(t,\infty)\times Y_j$, $j=r+1,\dots,r'$. We assume
$X\sco0$ is compact. For $i=0,\dots,r'$ let $y_i\in Y_i$ be a base-point
and set
\[e_i:=\begin{cases}
-1, & i=0,\dots,r,\\
1, & i=r+1,\dots,r'.
\end{cases}\]
For any integers $j,k$ in the interval $[0,r']$ such that $j<k$
let $\ga_{jk}:\R\to X$ be a smooth path satisfying
$\ga_{jk}(t)\in X\sco1$ for $|t|\le1$ and
\[\ga_{jk}(t)=\begin{cases}
(-e_jt,y_j), & t\le-1,\\
(e_kt,y_k), & t\ge1.
\end{cases}\]

\noindent Loosely speaking, the path $\ga_{jk}$ enters along the $j$th end
and leaves along the $k$th end of $X$.

Let $\va=(\al_1,\dots,\al_{r'})$, where $\al_j\in\fl(Y_j)$ and at least one
$\al_j$ is irreducible. For the remainder of this subsection we write
\[M:=M(X,E;\va),\]
where $E\to X$ is the product $\SO3$ bundle.
The unique continuation result of
Proposition~\refp{prop:unique-continuation-cylinder} ensures that if
$\al_j$ is irreducible then the restriction of any
element of $M$ to a band on the $j$th end of $X$ will be irreducible.

Let $\bbu\to M\times X$ be the universal (real) $3$--plane bundle (see
\cite[Subsection~5.1]{DK}).
For any $t\ge0$ let $\bbu\sco t$ denote the
restriction of $\bbu$ to $M\times X\sco t$. Given a base-point $x_0\in X$ let
$\bbe_{X,x_0;\va}\to M$ be the canonical $3$--plane bundle,
which can be identified
with the restriction of $\bbu$ to $M\times\{x_0\}$.

If $\ga:J\to X$ is a smooth path in $X$ defined on some interval $J$ then
a section $\si$ of the pull-back bundle $(\id\times\ga)^*\bbu$ over
$M\times J$ is called {\it holonomy invariant} if
for all $\om=[A]\in M$ and real numbers $s<t$ one has that $\si(\om,s)$
is mapped to $\si(\om,t)$ by the isomorphism
\begin{equation*}
\bbu_{(\om,\ga(s))}\to\bbu_{(\om,\ga(t))}
\end{equation*}
defined by holonomy of $A$ along the path $\ga|_{[s,t]}$.

Suppose $Z\subset X$ is a compact codimension~$0$ submanifold-with-boundary
such that $A|_Z$ is irreducible for every $[A]\in M$. Given a base-point
$z_0\in Z$, let $\bbe_{Z,z_0}\to\cb^*(E|_Z)$
be the base-point fibration, and let
\[R_Z:M\to\cb^*(E|_Z),\quad\om\mapsto\om|_Z.\]
Then the pull-back bundle $R_Z^*\bbe_{Z,z_0}$ is canonically isomorphic to
$\bbe_{X,z_0;\va}$, and we will usually identify the two bundles without further
comment.

Choose (smooth) sections $z_1,z_2,z_3$ of $\bbu\sco2$ and
for any $x\in X\sco2$ let
\begin{align*}
  &  M\cap w_3(x):=
  \{\om\in M\st z_1(\om,x)=0\},\\
  & M\cap w_2(x):=
  \{\om\in M\st\\
  & \hsp\text{$z_2,z_3$ are linearly
    dependent at $(\om,x)$}\}.
\end{align*}

For $j=0,\dots,r'$ let $\bbe_j\to\cb^*(Y_j[0])$ be the canonical
$3$--plane bundle associated to a base-point $(0,y_j)$.
For $j<k$, any $j'$, and $i=1,2,3$ choose
\begin{itemize}
\item a section $\uz ijk$ of $\bbe_j$
  and a section $\oz ijk$ of $\bbe_k$,
\item  a section $\zz ijk$ of $\bbu\sco2$,
  \item a section $s_{ij'}$ of $\bbe_{j'}$.
\end{itemize}

Let $b_{-1},b_0,b_1$ be a partion
of unity of $\R$ subordinate to the open cover
$\{(-\infty,-1),(-2,2),(1,\infty)\}$. 
If $j<k$ and both $\al_j,\al_k$ are
irreducible we introduce, for $i=1,2,3$, a section
of the bundle $(\id\times\ga_{jk})^*\bbu$ associated, loosely speaking,
to a base-point moving along the path $\ga_{jk}$. Precisely, we define
\[s_{ijk}(\om,t):=b_{-1}(t)\uz ijk(\om|_{Y_j[-e_jt]})
+b_0(t)\zz ijk(\om|_{X\sco2},\ga_{jk}(t))
+b_1(t)\oz ijk(\om|_{Y_k[e_kt]}).\]
Using these sections, we define cut-down moduli spaces
\begin{align*}
  & M\cap w_3(\ga_{jk}):=
  \{(\om,t)\in M\times\R\st s_{1jk}(\om,t)=0\},\\
  & M\cap w_2(\ga_{jk}):=
  \{(\om,t)\in M\times\R\st\\
  & \hsp\text{$s_{2jk},\,s_{3jk}$ are linearly
    dependent at $(\om,t)$}\}.
\end{align*}
We now consider the case of a base-point moving along the $j$th end.
For $t\ge0$ let $\ga_j(t):=(e_jt,y_j)$. If $\al_j$ is irreducible let
\begin{align*}
  & M\cap w_2(\ga_j):=\{(\om,t)\in M\times[0,\infty)\st\\
    & \hsp\text{$s_{2j},s_{3j}$ are linearly dependent at $\om|_{Y_j[e_jt]}$}\}.
\end{align*}

  We omit the definition of $M\cap w_3(\ga_j)$ since it will not be
  needed in the remainder of this paper
  (although something close to it was used in the proof of
  Proposition~\ref{prop:slide-points}).

We can also combine the ways moduli spaces are cut down in 
the above definitions. Namely, for $\ell,\ell'\in\{2,3\}$ let
\begin{align*}
  & M\cap w_\ell(x)\cap w_{\ell'}(\ga_{jk}):=
  \{(\om,t)\in M\cap w_{\ell'}(\ga_{jk})\st\\
  & \hsp\om\in M\cap w_\ell(x)\},\\
  & M\cap w_\ell(\ga_{jk})\cap w_{\ell'}(\ga_{j'k'}):=
  \{(\om,t,t')\in M\times\R\times\R\st\\
  & \hsp(\om,t)\in M\cap w_\ell(\ga_{jk}),\;
  (\om,t')\in M\cap w_{\ell'}(\ga_{j'k'})\},\\
  & M\cap w_\ell(\ga_{jk})\cap w_2(\ga_{j'}):=
  \{(\om,t,t')\in M\times\R\times[0,\infty)\st\\
  & \hsp(\om,t)\in M\cap w_\ell(\ga_{jk}),\;
  (\om,t')\in M\cap w_2(\ga_{j'})\}.
\end{align*}

If one of the $\al_j$s is trivial, say $\al_h=\theta$, and $\dim\,M<8$
(to prevent bubbling) then one can also cut
down $M$ by, loosely speaking, evaluating $w_2$ or $w_3$ over the
``link of $\theta$ at infinity'' over the $h$th end of $X$. We now make this
precise in the case of $w_2$ and an outgoing end $[0,\infty)\times Y_h$. The
  definitions for $w_3$ or incoming ends are similar. To simplify notation
  write $Y:=Y_h$.

We introduce a function $\tau^+=\tau^+_h$ on $M$ related to the energy
distribution of elements over the $h$th end.
  Choose $\eps>0$ so small that for any $\beta\in\fl(Y)$ the Chern-Simons
  value $\cs(\beta)\in\R/\z$ has no real lift in the interval $(0,\eps]$.
(Recall that we assume $\cs(\theta)=0$.)
Given $\om\in M$, if there exists a $t>0$ such that
$\ce_\om([t-2,\infty)\times Y)=\eps$ then $t$ is unique, and we write
  $t^+(\om):=t$. This defines $t^+$ implicitly as a smooth function on an open
  subset of $M$. We modify $t^+$ to get a smooth function
  $\tau^+:M\to[1,\infty)$ by
  \[\tau^+(\om):=
    \begin{cases}
      1+b(t^+(\om)-2)\cdot(t^+(\om)-1) & \text{if $t^+(\om)$ is defined,}\\
        1 & \text{else,}
    \end{cases}\]
    where the cut-off function $b$ is as in \refp{eqn:b-prop1}.
  Note that $\tau^+(\om)<3$ if $t^+(\om)<3$ and
  $\tau^+(\om)=t^+(\om)$ if $t^+(\om)\ge3$.
  The restriction of $\om$ to the band $Y[\tau^+(\om)]$ will be denoted by
  $R^+(\om)\in\cb(Y[0])$.

  \begin{lemma}
    In the above situation there is a real number $T_0$
    such that if $\om$ is any element of $M$ satisfying
    $\tau^+(\om)>T_0-1$ then $R^+(\om)$ is irreducible.
  \end{lemma}
  \proof Suppose the lemma is false. Then we can find a sequence $\om_n$
  in $M$ such that $\tau^+(\om_n)\to\infty$ and
  $R^+(\om_n)$ is reducible for every $n$. Let $A_n$ be a smooth connection
  representing $\om_n$, and let $t_n=\tau^+(\om_n)$.
  By assumption, there is no bubbling in $M$, so
  we can find gauge transformations $u_n$ defined over $[0,\infty)\times Y$
  and a smooth connection $A'$ over $\ry$ such that, for every constant
  $c>0$, the sequence $u_n(A_n)|_{[t_n-c,t_n+c]}$ converges in $C^\infty$
  to $A'|_{[-c,c]}$. The assumption on $\eps$ means that no energy can be
  lost over the end $[0,\infty)\times Y$ in the limit, hence
  \[\ce_{A'}([-2,\infty)\times Y)=\eps.\]
  In particular, $A'$ is not trivial. But there are no non-trivial reducible
  finite-energy instantons over $\ry$ (as long as the perturbation of the
  Chern-Simons functional is so small that there are no
  non-trivial reducible critical points).
  Therefore, $A'$ must be irreducible.  From the unique continuation result of
  Proposition~\ref{prop:unique-continuation-cylinder} it follows that
  $A'|_{\{0\}\times Y}$ is
  also irreducible, so $A_n$ is irreducible for large $n$.
  This contradiction proves the lemma.
  \square
  
  Let $T_0$ be as in the lemma. For any element of $M$ for which
  $R^+(\om)$ is irreducible, let
  $s'_{ih}(\om)$ denote the holonomy invariant section of
  $(\id\times\ga_h)^*\bbu$ such that $s'_{ih}(\om,\tau^+(\om))=s_{ih}(R^+(\om))$.
Let $x_h:=(0,y_h)$ and define a section of $\bbe_{X,x_h;\va}$ by  
  \[s_{ih}(\om):=(1-b(\tau^+(\om)-T_0))\cdot z_i(\om|_{X\sco2},x_h)
  +b(\tau^+(\om)-T_0)\cdot s'_{ih}(R^+(\om)),\]
  where again $b$ is as in \refp{eqn:b-prop1}.
Let
\[M\cap w_2(\tau^+):=\{\om\in M\st\text{$s_{2h},s_{3h}$ linearly dependent
    at $\om$}\}.\]
If $j<k$ and both $\al_j,\al_k$ are irreducible let
\[M\cap w_{\ell}(\ga_{jk})\cap w_2(\tau^+):=
  \{(\om,t)\in M\cap w_\ell(\ga_{jk})\st\om\in M\cap w_2(\tau^+)\}.\]
  
If $M$ is regular, then the various cut down moduli
spaces defined above will be transversely cut out when the sections involved
are ``generic''.

\subsection{Operations, I}\label{subsec:opI}

We now specialize to the case when $X$ has two incoming ends
$(-\infty,0]\times Y_j$, $j=0,1$ and one outgoing end $[0,\infty)\times Y_2$,
    and
    \[H_i(X;\z)=0,\quad i=1,2.\]
Such a cobordism gives rise to a homomorphism
\begin{equation}\label{eqn:cv}
  A:C^p(Y_0)\otimes C^q(Y_1)\to C^{p+q}(Y_2)
  \end{equation}
for any $p,q\in\z/8$, with matrix coefficients
\[\la A(\al_0\otimes\al_1),\al_2\ra:=\#M(X;\va)\]
for generators $\al_0\in C^p(Y_0)$, $\al_1\in C^q(Y_1)$, and
$\al_2\in C^{p+q}(Y_2)$, where $\va=(\al_0,\al_1,\al_2)$. We can construct
more homomorphisms using the sections $s_{ijk}$ chosen above.
For any path $\ga_{jk}$ as above and $k=2,3$ let
\[T_{i,j,k}:C^p(Y_0)\otimes C^q(Y_1)\to C^{p+q+i-1}(Y_2)\]
be defined on generators by
\[\la T_{i,j,k}(\al_0\otimes\al_1),\al_2\ra:=
\#[M(X;\va)\cap w_i(\ga_{jk})].\]
For the cases used in this paper we introduce the simpler notation
\[B:=T_{3,0,1},\quad E:=T_{3,0,2},\quad A':=T_{2,1,2}.\]
We will also consider homomorphisms defined using two base-points, each
moving along a path in $X$. At this point we only define
\[B':C^p(Y_0)\otimes C^q(Y_1)\to C^{p+q+3}(Y_2)\]
by
\[\la B'(\al_0\otimes\al_1),\al_2\ra:=
\#[M(X;\va)\cap w_3(\ga_{01})\cap w_2(\ga_{12})].\]

In the next proposition, the differential in the cochain complex
$C(Y_i)$ will be denoted by $d$ (for $i=0,1,2$), and
\[\ti d=d\otimes1+1\otimes d\]
will denote the differential in $C(Y_0)\otimes C(Y_1)$.
Let
\[\ti v_3:=v_3\otimes1+1\otimes v_3,\]
regarded as a degree~$3$ cochain map from $C(Y_0)\otimes C(Y_1)$ to itself.

\begin{prop}\label{prop:dabe}
\begin{description}
\item[(i)] $dA+A\ti d=0$.
\item[(ii)] $dB+B\ti d=A\ti v_3$.
\item[(iii)] $dE+E\ti d=A(v_3\otimes1)+v_3A$.
\item[(iv)] $dA'+A'\ti d=A(1\otimes v_2)+v_2A$.
\item[(v)] $dB'+B'\ti d=B(1\otimes v_2)+v_2B
  +A'\ti v_3+A(1\otimes\phi)+A_\theta(1\otimes\del)$.
\end{description}
\end{prop}

\proof The only non-trivial part here is (v), where one encounters
factorization through the trivial connection over the end
$(-\infty,0]\times Y_1$. This can be handled as in the proof of
  Proposition~\ref{prop:v2v3phi} given in Subsection~\ref{subsec:proof-v2v3phi},
  to which we refer for details.\square

\begin{prop}\label{prop:clmc}
  The homomorphism
  \begin{align*}
    \cl:MC^*(Y_0)\otimes MC^*(Y_1)&\to C^*(Y_2),\\
    (x_0,y_0)\otimes(x_1,y_1)&\mapsto B(x_0,x_1)+A(x_0\otimes y_1+y_0\otimes x_1)
  \end{align*}
  is a cochain map of degree $-2$.
\end{prop}

\proof Let $\ti D=D\otimes1+1\otimes D$ be the differential in the complex
$MC(Y_1)\otimes MC(Y_2)$. Then
\begin{align*}
  \cl\ti D[(x_0,y_0)\otimes(x_1,y_1)]&=
  \cl[(dx_0,v_3x_0+dy_0)\otimes(x_1,y_1)+(x_0,y_0)\otimes(dx_1,v_3x_1+dy_1)]\\
  &=B(dx_0\otimes x_1+x_0\otimes dx_1)\\
  &+A[dx_0\otimes y_1+(v_3x_0+dy_0)\otimes x_1+
  x_0\otimes(v_3x_1+dy_1)+y_0\otimes dx_1]\\
  &=B\ti d(x_0\otimes x_1)
  +A\left[\ti v_3(x_0\otimes x_1)+\ti d(x_0\otimes y_1+y_0\otimes x_1)\right]\\
  &=d\cl[(x_0,y_0)\otimes(x_1,y_1)],  
\end{align*}
where the last equality follows from Proposition~\ref{prop:dabe}.\square

The homomorphism
  \[MI^*(Y_0)\otimes MI^*(Y_1)\to I^*(Y_2)\]
obtained from Proposition~\ref{prop:clmc} will also be denoted by $\cl$.

In order to simplify notation we will often write $\del,\Del$
instead of $\del_0,\Del_0$ if no confusion can arise.

  \begin{prop}\label{prop:clv}
    For all $a\in MI(Y_0)$, $b\in MI(Y_1)$, the following hold.
\begin{description}
\item[(i)] If $\Del a=0$ then $\cl(Ua,b)=u_2\cl(a,b)$.
\item[(ii)] If $\Del b=0$ then $\cl(a,Ub)=u_2\cl(a,b)$.  
\end{description}
\end{prop}

  \proof We spell out the proof of (ii). Reversing the roles of $Y_0,Y_1$
  yields a proof of (i). Let
  \[\cl',\ce:MC^*(Y_0)\otimes MC^*(Y_1)\to C^*(Y_2)\]
  be given by
  \begin{align*}
   \cl'[(x_0,y_0)\otimes(x_1,y_1)]
    &:= B'(x_0,x_1)+A'(x_0\otimes y_1+y_0\otimes x_1),\\
    \ce[(x_0,y_0)\otimes(x_1,y_1)]
    &:=(\del x_1)A_\theta(x_0).
\end{align*}
  Let $\ti D$ be as in the proof of Proposition~\ref{prop:clmc}. We show that
  \[d\cl'+\cl'\ti D=v_2\cl+\cl(1\times V)+\ce,\]
  from which (ii) follows. Observe that the first four lines
  in the calculation of $\cl\ti D$ in 
  Proposition~\ref{prop:clmc} carry over to $\cl'\ti D$. 
  That proposition then gives
  \begin{align*}
    \cl'\ti D&[(x_0,y_0)\otimes(x_1,y_1)]\\
    &=(B'\ti d+A'\ti v_3)(x_0\otimes x_1)
    +A'\ti d(x_0\otimes y_1+y_0\otimes x_1)\\
    &=dB'(x_0\otimes x_1)+B(x_0\otimes v_2x_1)+v_2B(x_0\otimes x_1)
    +A(x_0\otimes\phi x_1)+(\del x_1)A_\theta(x_0)\\
    &\hphantom{=}+[dA'+A(1\otimes v_2)+v_2A](x_0\otimes y_1+y_0\otimes x_1)\\
    &=[d\cl'+v_2\cl+\cl(1\times V)+\ce][(x_0,y_0)\otimes(x_1,y_1)].\square
    \end{align*}

Our next goal is to compute $\del u_2\cl$. To this end we introduce some
variants $\dot A,\dot B,A^+,B^+$ of the operators $A,B$. Each of these
variants is a homomorphism
\[C^p(Y_0)\otimes C^q(Y_1)\to C^{p+q+d}(Y_2)\]
for $d=2,4,1,3$, respectively, defined for all $p,q$, and the matrix
coefficients are
\begin{align*}
  \la\dot A(\al_0\otimes\al_1),\al_2\ra&:=
  \#[M(X;\va)\cap w_2(x_2)],\\
  \la\dot B(\al_0\otimes\al_1),\al_2\ra&:=
  \#[M(X;\va)\cap w_2(x_2)\cap w_3(\ga_{01})],\\
 \la A^+(\al_0\otimes\al_1),\al_2\ra&:=
 \#[M(X;\va)\cap w_2(\ga_2)],\\
 \la B^+(\al_0\otimes\al_1),\al_2\ra&:=
 \#[M(X;\va)\cap w_3(\ga_{01})\cap w_2(\ga_2)],
\end{align*}
where $\va=(\al_0,\al_1,\al_2)$ as before, $x_2=\ga_2(0)\in X$, and
$\ga_i,\ga_{ij}$ are as in Subsection~\ref{subsec:cutting-down-mod-sp}.

\begin{prop}\label{prop:ddotA}
  \begin{description}
  \item[(i)] $d\dot A+\dot A\ti d=0$.
  \item[(ii)] $d\dot B+\dot B\ti d=\dot A\ti v_3$.
  \item[(iii)] $dA^++A^+\ti d=v_2A+\dot A$.
  \item[(iv)] $dB^++B^+\ti d=A^+\ti v_3+v_2B+\dot B$.
  \end{description}
 \end{prop}

\proof Standard.\square

\begin{prop}\label{prop:dotcl}
  The homomorphism
  \begin{align*}
   \dot\cl:MC^*(Y_0)\otimes MC^*(Y_1)&\to C^*(Y_2),\\
   (x_0,y_0)\otimes(x_1,y_1)&\mapsto\dot B(x_0,x_1)
   +\dot A(x_0\otimes y_1+y_0\otimes x_1)
  \end{align*}
  is a (degree preserving) cochain map.
  \end{prop}
\proof The same as for Proposition~\ref{prop:clmc}, using
Proposition~\ref{prop:ddotA} (i), (ii).\square

The homomorphism
  \[MI^*(Y_0)\otimes MI^*(Y_1)\to I^*(Y_2)\]
  obtained from Proposition~\ref{prop:dotcl} will also be denoted by $\dot\cl$.

 \begin{prop}\label{prop:u2cl}
  As maps $MI^*(Y_0)\otimes MI^*(Y_1)\to I^*(Y_2)$ one has
  \[\dot\cl=u_2\cl.\]
\end{prop}

\proof This is analogous to the proof of Proposition~\ref{prop:clv}. Let
\[\cl^+:MC^*(Y_0)\otimes MC^*(Y_1)\to C^*(Y_2)\]
be given by
\[\cl^+[(x_0,y_0)\otimes(x_1,y_1)]
    := B^+(x_0,x_1)+A^+(x_0\otimes y_1+y_0\otimes x_1).\]
We show that
\[d\cl^++\cl^+\ti d=v_2\cl+\dot\cl.\]
From Proposition~\ref{prop:ddotA} we get
\begin{align*}
  \cl^+\ti D(x_0,y_0)\otimes(x_1,y_1)
  &=(B^+\ti d+A^+\ti v_3)(x_0\otimes x_1)+A^+\ti d(x_0\otimes y_1+
  y_0\otimes x_1)\\
  &=(dB^++v_2B+\dot B)(x_0\otimes x_1)
  +(dA^++v_2A+\dot A)(x_0\otimes x_1)\\
  &=(d\cl^++v_2\cl+\dot\cl)(x_0,y_0)\otimes(x_1,y_1).\square
  \end{align*}

We also need to bring in moduli spaces over $X$ with trivial limit over the
end $\R_+\times Y_2$. These give rise to homomorphisms
\[A^\theta,B^\theta,\dot A^\theta,\dot B^\theta:C^p(Y_0)\otimes C^{d-p}(Y_1)\to\z/2\]
where $d=5,3,3,1$, respectively. They are defined on generators by
\begin{align*}
  A^\theta(\al_0\otimes\al_1)&:=\#M(\al_0,\al_1,\theta),\\
  B^\theta(\al_0\otimes\al_1)&:=\#[M(\al_0,\al_1,\theta)\cap w_3(\ga_{01})],\\
  \dot A^\theta(\al_0\otimes\al_1)&:=\#[M(\al_0,\al_1,\theta)\cap w_2(x_0),\\
    \dot B^\theta(\al_0\otimes\al_1)
    &:=\#[M(\al_0,\al_1,\theta)\cap w_2(x_0)\cap w_3(\ga_{01}).
  \end{align*}

\begin{prop}\label{prop:deldotA}
  \begin{description}
  \item[(i)] $\del A+ A^\theta\ti d=0$.
  \item[(ii)] $\del B+B^\theta\ti d=A^\theta\ti v_3$.
  \item[(iii)] $\del\dot A+\dot A^\theta\ti d=0$.
  \item[(iv)] $\del\dot B+\dot B^\theta\ti d=
    \dot A^\theta\ti v_3+\del\otimes\del$.
 \end{description}
\end{prop}
Here, $(\del\otimes\del)(x_0\otimes x_1)=(\del x_0)(\del x_1)$.

\proof The term $\del\otimes\del$ in (iv) accounts for factorization through
the trivial connection over $X$, see Subsection~\ref{subsec:proof-v2v3phi}
below. The remaining parts of the proof are standard.\square

\begin{prop}\label{prop:delu2cl}
   \begin{description}
  \item[(i)] $\del\cl=0$.
  \item[(ii)] $\del u_2\cl=\Del\otimes\Del$.
 \end{description}
  \end{prop}

\proof Statement (i) is proved just as Proposition~\ref{prop:clmc}, replacing
Proposition~\ref{prop:dabe} by Proposition~\ref{prop:deldotA}.
We now prove (ii). For $g_i=(x_i,y_i)\in MC(C_i)$, $i=0,1$ let
 \[\dot\cl^\theta(g_0\otimes g_1):=\dot B^\theta(x_0\otimes x_1)
 +\dot A^\theta(x_0\otimes y_1+y_0\otimes x_1).\]
 Arguing as in the proof of Proposition~\ref{prop:clmc} and using
 Proposition~\ref{prop:deldotA} we obtain
 \begin{align*}
   \dot\cl^\theta\ti D(g_0\otimes g_1)
   &=(\dot B^\theta\ti d+\dot A\ti v_3)(x_0\otimes x_1)
   +\dot A^\theta\ti d(x_0\otimes y_1+y_0\otimes x_1)\\
   &=\del\dot B(x_0\otimes x_1)+\del x_0\cdot\del x_1
   +\del\dot A(x_0\otimes y_1+y_0\otimes x_1)\\
   &=(\del\dot\cl+\Del\otimes\Del)(g_0\otimes g_1).
   \end{align*}
If $g_0,g_1$ are cocycles then by Proposition~\ref{prop:u2cl} we have
\[\del v_2\cl(g_0\otimes g_1)=\del\dot\cl(g_0\otimes g_1)
=\Del g_0\cdot\Del g_1.\square\]

  For $p\neq4$ let
\begin{equation}\label{eqn:Fdef}
  F:C^p(Y_0)\otimes C^q(Y_1)\to C^{p+q+4}(Y_2)
\end{equation}
be defined by
  \[\la F(\al_0\otimes\al_1),\al_2\ra:=
  \#[M(X;\va)\cap w_3(\ga_{01})\cap w_3(\ga_{02})].\]
  For $p=4$ the map $F$ may not be well-defined due to possible factorizations
  through the trivial connection over the end $\R_-\times Y_0$.
  
The definition of $F$ involves two different sections of the bundle
$\bbe_0\to\cb^*(Y_0[0])$, namely
\[s_k:=\uz 10k,\quad k=1,2.\]
From now on we assume $s_1,s_2$ are so close that
they define the same cup product $v_3:C^*(Y_0)\to C^{*+3}(Y_0)$.

\begin{prop}\label{prop:Fthm}
  If the sections $s_1,s_2$ are sufficiently close then the map
  $F$ in \refp{eqn:Fdef} can be extended to all bidegrees $(p,q)$ such that
  \begin{equation}\label{eqn:Fthm}
    dF+F\ti d=B(v_3\otimes1)+v_3B+E\ti v_3+A(\psi\otimes1),
    \end{equation}
  where $\psi$ is as in Proposition~\ref{prop:psixi}.
\end{prop}

The main difficulty in extending the map $F$ to degree~$p=4$,
related to factorization through the trivial connection over the end
$(-\infty,0]\times Y_0$, is the same as in
extending the map $\psi$ to degree $4$, and the main difficulty in
proving \refp{eqn:Fthm} is the same as in proving that $\psi$ is a cochain map
(Proposition~\ref{prop:psichainmap}). As we prefer to explain the ideas involved
in the simplest possible setting, we will not spell out the proof
of Proposition~\ref{prop:Fthm} but instead refer to
Subsection~\ref{subsec:cochain-map-psi} for details.

Sometimes we will fix the variable $\al_1$ in the expressions defining
$A,B,E,F$. Thus, for any $y\in C^r(Y)$ we define a homomorphism
\[A_y:C^*(Y_0)\to C^{*-r}(Y_2),\quad x\mapsto A(x\otimes y),\]
and we define $B_y,E_y,F_y$ similarly. Looking at moduli spaces over $X$
with trivial limit over the end $\R_-\times Y_1$ we obtain homomorphisms
\begin{align*}
A_\theta&:C^*(Y_0)\to C^*(Y_2),\\
E_\theta&:C^*(Y_0)\to C^{*+2}(Y_2).
\end{align*}
with matrix coefficients
\begin{align*}
  \la A_\theta(\al_0),\al_2\ra&:=\#M(X;\al_0,\theta,\al_2),\\
  \la E_\theta(\al_0),\al_2\ra&:=\#[M(X;\al_0,\theta,\al_2)\cap w_3(\ga_{02})].
  \end{align*}

We consider a variant of Floer's complex introduced by Donaldson
\cite[p.\ 169]{D5}.
For any oriented homology $3$--sphere $Y$ let $\oc*(Y)$ be the complex
with cochain groups
\begin{align*}
  \oc p(Y)&=C^p(Y),\quad p\neq0,\\
\oc0(Y)&=C^0(Y)\oplus\z/2  
  \end{align*}
and differential $\bar d=d+\del'$.
Now take $Y:=Y_1$. For $y=(z,t)\in\oc0(Y_1)$ let
\[A_y:=A_z+tA_\theta,\quad E_y:=E_z+tE_\theta.\]

\begin{lemma}\label{lemma:dby}
  For any $x\in C(Y_1)$ and $y\in\oc*(Y_1)$ we have
  \begin{align*}
    [d,A_y]+A_{\bar dy}&=0,\\
    [d,E_y]+E_{\bar dy}&=[A_y,v_3],\\
    [d,B_x]+B_{dx}&=A_xv_3+A_{v_3x},\\
    [d,F_x]+F_{dx}&=[B_x,v_3]+E_xv_3+E_{v_3x}+A_x\psi.    
    \end{align*}
\end{lemma}
Here, $[d,A_y]=dA_y+A_yd$, and similarly for the other commutators.

\proof For $y\in C(Y_1)$ this follows from Propositions~\ref{prop:dabe} and
\ref{prop:Fthm}, whereas the case $y=(0,1)\in\oc0(Y_1)$ is easy.\square

\begin{lemma}\label{lemma:ckcm}
  Suppose $x\in C^{-2}(Y_1)$ and $y=(z,t)\in\oc0(Y_1)$ satisfy
  \[dx=0,\quad v_3x=\bar dy.\]
  Then the homomorphism $\ck:MC^*(Y_0)\to MC^*(Y_2)$ given by the matrix
  \[\left(
  \begin{array}{cc}
    A_y+B_x & A_x\\
    E_y+F_x+A_x\Xi & A_y+B_x+E_x+A_xv_2
    \end{array}
  \right)\]
  is a cochain map. Here, $\Xi$ is as in Proposition~\ref{prop:psixi}.
  \end{lemma}

\proof Writing $\ck=\left(\begin{array}{cc} P & Q \\ R & S\end{array}\right)$
  we have
\[\ti d\ck+\ck\ti d=\left(
\begin{array}{cc}
  dP+Pd+Qv_3 & dQ+Qd \\
  dR+Rd+v_3P+Sv_3 & dS+Sd+v_3Q
    \end{array}
\right).\]
The fact that this matrix vanishes is easily deduced from
Propositions~\ref{prop:dvvd} and \ref{prop:psixi} and Lemma~\ref{lemma:dby}.
We write out the calculation only for the bottom
left entry.
\begin{align*}
  [d,E_y&+F_x+A_x\Xi]\\
  &=E_{v_3x}+[v_3,A_y]+[v_3,B_x]+E_{v_3x}+E_xv_3+A_x\psi+A_x[d,\Xi]\\
  &=v_3(A_y+B_x)+(A_y+B_x+E_x+A_xv_2)v_3,
\end{align*}
hence $[d,R]=v_3P+Sv_3$ as claimed.\square

\begin{prop}\label{prop:u3clz}
  As maps $MI^*(Y_0)\otimes MI^*(Y_1)\to I^*(Y_2)$ one has
  \[u_3\cl=0.\]
\end{prop}

\proof For $j=0,1$ let $(x_j,y_j)$ be a cocycle in $MC(Y_j)$, i.e.
\[dx_j=0,\quad v_3x_j=dy_j.\]
Let the map $\ck$ of Lemma~\ref{lemma:ckcm} be defined with
$x=x_1$, $y=y_1$, and let $(x_2,y_2):=\ck(x_0,y_0)$. Then
\[\cl((x_0,y_0)\otimes(x_1,y_1))=B_{x_1}(x_0)+A_{y_1}(x_0)+A_{x_1}(y_0)=x_2.\]
Since $(x_2,y_2)$ is a cocycle, we have $v_3x_2=dy_2$, proving the proposition.
\square

 \begin{prop}\label{prop:qtsubaddI}
If $\qt(Y_j)\ge1$ for $j=0,1$ then
\[\qt(Y_2)\ge\qt(Y_0)+\qt(Y_1).\]
  \end{prop}

 \proof For $j=0,1$ let $n_j:=\qt(Y_j)$ and choose $z_j\in MI(Y_j)$ such that
 \[\Del U^kz_j=\be{cases}
0 & \text{for $0\le k<n_j-1$,}\\
1 & \text{for $k=n_j-1$.}
\end{cases}\]
Let $x:=\cl(z_0\otimes z_1)\in I(Y_2)$. Then $u_3x=0$ by
Proposition~\ref{prop:u3clz}. For $0\le k_j\le n_j-1$, repeated
application of Proposition~\ref{prop:clv} yields
\[u_2^{k_0+k_1}x=\cl(U^{k_0}z_0\otimes U^{k_1}z_1),\]
hence $\del u_2^{k_0+k_1}x=0$ by Proposition~\ref{prop:delu2cl}. Therefore,
\[\del u_2^mx=0,\quad 0\le m\le n_1+n_2-2.\]
On the other hand,
\begin{align*}
  \del u_2^{n_1+n_2-1}x&=\del u_2u_2^{n_0-1}u_2^{n_1-1}x\\
  &=\del u_2\cl(U^{n_0-1}z_0\otimes U^{n_1-1}z_1)\\
  &=(\Del U^{n_0-1}z_0)(\Del U^{n_1-1}z_1)\\
  &=1.
\end{align*}
Therefore, $\qt(Y_2)\ge n_0+n_1$ as claimed.\square

We will give a second application of Lemma~\ref{lemma:ckcm}, but first we need
some preparation. Let $A^\theta_\theta:C^5(Y_0)\to\z/2$
be defined on generators by
\[A^\theta_\theta(\al):=\#M(\al,\theta,\theta).\]
For $y=(z,t)\in\oc q(Y_1)$ define
$A^\theta_y:C^{5-q}(Y_0)\to\z/2$ and $B^\theta_z:C^{3-q}(Y_0)\to\z/2$
by
\[A^\theta_y(x):=A(x\otimes z)+tA^\theta_\theta(x),\quad
B^\theta_z(x):=B^\theta(x\otimes z).\]

\begin{lemma}\label{lemma:delAB}
  \begin{description}
  \item[(i)] $\del A_\theta+A^\theta_\theta d+A^\theta_{\del'(1)}=\del$.
  \item[(ii)] $\del A_y+A^\theta_y d+A^\theta_{\bar dy}=t\del$.
  \item[(iii)] $\del B_z+B^\theta_z d+B^\theta_{dz}=A^\theta_zv_3+A^\theta_{v_3z}$.
  \end{description}
  \end{lemma}

\proof Standard.\square

 \begin{prop}\label{prop:onezero1}
   If $\qt(Y_0)\ge1$ and $\qt(Y_1)=0$ then $\qt(Y_2)\ge1$.
   \end{prop}

 \proof Since $\qt(Y_0)\ge1$ we can find $(x_0,y_0)\in MC^6(Y_0)$ such
 that
 \[dx_0=0,\quad v_3x_0=dy_0,\quad\del x_0=1.\]
 Since $\qt(Y_1)=0$, Lemma~\ref{lemma:q2ymin} says that there exist
 $x_1\in C^{-2}(Y_1)$ and $y_1=(z_1,1)\in\oc 0(Y_1)$ such that
 \[dx_1=0,\quad v_3x_1=\bar dy_1.\]
 Let $\ck$ be as in Lemma~\ref{lemma:ckcm}. Then $\ck(x_0,y_0)$ is a cocycle
 in $MC(Y_2)$, and by Lemma~\ref{lemma:delAB} we have
 \begin{align*}
   \Del\ck(x_0,y_0)&=\del(A_{y_1}+B_{x_1})x_0+\del A_{x_1}y_0\\
   &=(\ath_{\bar dy_1}+\del+\ath_{x_1}v_3+\ath_{v_3x_1})x_0+\ath_{x_1}dy_0\\
   &=1.
   \end{align*}
Therefore, $\qt(Y_2)\ge1$.\square

\subsection{Operations, II}\label{subsec:OpII}

We now consider the case when $X$ has one incoming end $(-\infty,0]\times Y_0$
and two outgoing ends $[0,\infty)\times Y_1$ and $[0,\infty)\times Y_2$,
    where $Y_2=\Si=\Si(2,3,5)$ is the Poincar\'e homology sphere oriented as the
boundary of the negative definite $E_8$--manifold. We again assume that
    \[H_i(X;\z)=0,\quad i=1,2.\]
We will define homomorphisms
\[P,P',Q:C^*(Y_0)\to C^{*+d}(Y_1)\]
where $d=2,3,4$, respectively, making use of cut-down moduli spaces
introduced at the end of Subsection~\ref{subsec:cutting-down-mod-sp} with
$h=2$, so that $\tau^+=\tau^+_2$. 
We define $P,P',Q$ on generators by
\begin{align*}
  \la P\al_0,\al_1\ra&:=\#[M(X;\al_0,\al_1,\theta)\cap w_2(\tau^+)],\\
  \la P'\al_0,\al_1\ra&:=
  \#[M(X;\al_0,\al_1,\theta)\cap w_2(\ga_{01})\cap w_2(\tau^+)],\\
   \la Q\al_0,\al_1\ra&:=
  \#[M(X;\al_0,\al_1,\theta)\cap w_3(\ga_{01})\cap w_2(\tau^+)].  
\end{align*}

\begin{prop}\label{prop:dPPQ}
  As maps $C(Y_0)\to C(Y_1)$ the following hold.
  \begin{description}
  \item[(i)] $[d,P]=0$.
  \item[(ii)] $[d,P']=[v_2,P]$.
  \item[(iii)] $[d,Q]=[v_3,P]+\del'\del$.
    \item[(iv)] $\del P+Pd=\dot\del$.
    \end{description}
  \end{prop}
Here, $\dot\del$ is as defined at the end of Subsection~\ref{subsec:cup-prod}.

\proof In (iii), argue as in the proof of Proposition~\ref{prop:v2v3phi}
to handle
factorization through the trivial connection over $X$.\square

Note that statements (i), (iii) are equivalent to the fact that the
homomorphism
\[\Psi=
\left(\begin{array}{cc}
  P & 0\\
  Q & P
    \end{array}\right)
:MC^*(Y_0)\to MC^{*+2}(Y_1)\]
satisfies
\[[D,\Psi]=\Del'\Del.\]


The homomorphism $I^*(Y_0)\to I^{*+2}(Y_1)$ induced by $P$ will also be denoted
by $P$.

\begin{prop}\label{prop:u2u3P}
    As maps $I(Y_0)\to I(Y_1)$ the following hold.
  \begin{description}
  \item[(i)] $[u_2,P]=0$.
  \item[(ii)] $[u_3,P]=\del'\del$.
  \item[(iii)] $\del P=\del u_2$.
    \end{description}
  \end{prop}

\proof Combine Propositions~\ref{prop:delv2deld} and \ref{prop:dPPQ}.\square

\begin{prop}\label{prop:q2ge2}
  If $\qt(Y_0)\ge2$ then
  \[\qt(Y_1)\ge\qt(Y_0)-1.\]
\end{prop}

\proof Let $n:=\qt(Y_0)$ and choose $x\in I(Y_0)$ such that $u_3x=0$ and
\[\del u_2^kx=\begin{cases}
0 & \text{for $0\le k<n-1$,}\\
1 & \text{for $k=n-1$.}
\end{cases}\]
By Proposition~\ref{prop:u2u3P} we have $u_3Px=0$ and
\[\del u_2^kPx=\del Pu_2^kx=\del u_2^{k+1}x=\begin{cases}
0 & \text{for $0\le k<n-2$,}\\
1 & \text{for $k=n-2$.}
\end{cases}\]
This shows that $\qt(Y_1)\ge n-1$.\square

\subsection{Additivity of $\qt$}

Throughout this subsection, $Y,Y_0,Y_1$ will denote oriented homology
$3$--spheres. As before, $\Si$ will denote the Poincar\'e homology sphere.

\begin{prop}\label{prop:qtsubaddII}
If $\qt(Y_j)\ge1$ for $j=1,2$ then
\[\qt(Y_0\# Y_1)\ge\qt(Y_0)+\qt(Y_1).\]
\end{prop}

\proof Recall that there is a standard cobordism $W$ from $(-Y_0)\cup(-Y_1)$
to $Y_0\# Y_1$. By attaching half-infinite tubular ends to $W$ we obtain
a manifold $X$ to which we can apply the results of
Subsection~\ref{subsec:opI}. The proposition now follows from
Proposition~\ref{prop:qtsubaddI}.\square

\begin{prop}\label{prop:onezeroII}
If $\qt(Y_0)\ge1$ and $\qt(Y_1\#(-Y_0))=0$ then $\qt(Y_1)\ge1$.
\end{prop}

\proof This follows from Proposition~\ref{prop:onezero1}.\square

\begin{prop}\label{prop:ysi1}
If $\qt(Y\#\Si)\ge2$ then
\[\qt(Y)\ge\qt(Y\#\Si)-1.\]
\end{prop}

\proof This follows from Proposition~\ref{prop:q2ge2}
with $Y_0=Y\#\Si$ and $Y_1=Y$.
\square

In the following, we write $Y_0\sim Y_1$ to indicate that $Y_0$ and
$Y_1$ are homology cobordant.

\begin{lemma}\label{lemma:poinc-add}
If $Y_0\# Y_1\sim\Si$ then $\qt(Y_0)+\qt(Y_1)=1$.
\end{lemma}

\proof Let $k_j:=\qt(Y_j)$.

{\it Case 1:} $n_0n_1=0$. Without loss of generality we may assume that
$n_1=0$. By Proposition~\ref{prop:onezeroII} we have $n_0\ge1$.
If $n_0\ge2$ then,
since $Y_0\sim\Si\#(-Y_1)$, Proposition~\ref{prop:ysi1} would give
\[-n_1=\qt(-Y_1)\ge\qt(\Si\#(-Y_1)-1\ge1,\]
a contradiction. Hence, $n_0=1$, so the lemma holds in this case.

{\it Case 2:} $n_0n_1>0$. We show that this cannot occur.
If $k_j>0$ then Proposition~\ref{prop:qtsubaddII} yields
\[1=\qt(\Si)\ge n_0+n_1\ge2,\]
a contradiction. Similarly, if $k_j<0$ then the same proposition yields
$-1=\qt(-\Si)\ge2$. 

{\it Case 3:} $n_0n_1<0$. Then we may assume that $n_0>0$.
Applying Proposition~\ref{prop:qtsubaddII} we obtain
\[n_0=\qt(\Si\#(-Y_1))\ge1-n_1\ge2.\]
Proposition~\ref{prop:ysi1} now gives $-n_1\ge n_0-1$.
Altogether, this shows that
$n_0+n_1=1$.\square

\begin{cor}\label{cor:Siadd}
$\qt(Y\#\Si)=\qt(Y)+1$.
\end{cor}

\proof Apply the lemma with $Y_0=Y\#\Si$ and $Y_1=-Y$.\square

\begin{thm}
For any oriented integral homology $3$--spheres $Y_0,Y_1$ one has
\[\qt(Y_0\# Y_1)=\qt(Y_0)+\qt(Y_1).\]
\end{thm}

\proof Let $k_j:=\qt(Y_j)$ and $Z_j:=Y_j\#(-k_j\Si)$.
By Corollary~\ref{cor:Siadd} we have $\qt(Z_j)$=0, so by
Proposition~\ref{prop:onezeroII},
\[0=\qt(Z_0\# Z_1)=\qt(Y_0\# Y_1\#(-n_0-n_1)\Si)=\qt(Y_0\# Y_1)-n_0-n_1.
\square\]

\section{Further properties of $\qt$. Examples}\label{sec:further-props}

\subsection{Proof of Theorem~\ref{thm:monotone}}

Let $W'$ be the result of connecting the two boundary components of $W$
by a $1$--handle. Then $W$ and $W'$ have the same second cohomology group
and the same intersection form.

Let $Z$ be the negative definite $E_8$--manifold (i.e. the result of
plumbing on the $E_8$ graph), so that the boundary of $Z$
is the Poincar\'e sphere $\Si$. We will apply Theorem~\ref{thm:defub-gen}
to the boundary-connected sum
\[V:=W'\#_\partial Z.\]
Let $S,S'\subset Z$ be embedded oriented $2$--spheres corresponding to adjacent
nodes on the $E_8$ graph. These spheres both have self-intersection number
$-2$, and $S\cdot S'=1$. Let
\[v=\text{P.D.}([S])\in H^2(V,\prtl V)\approx H^2(V)\]
be the Poincar\'e dual of the homology class in $V$ represented by $S$. Then
$v\cdot[S']=1$, hence $v$ is strongly admissible. The class
$w\in J_V$ represented by $v$ satisfies $w^2=-2$, and
$\pm w$ are the only classes in $w+2J_V$ with square norm $2$.
Theorem~\ref{thm:defub-gen}
and Proposition~\ref{prop:Brieskorn-calc} now yield
\[\qt(Y)+1=\qt(Y\#\Si)\ge1,\]
hence $\qt(Y)\ge0$ as claimed.\square

\subsection{Proof of Theorem~\ref{thm:surgery-bound}}

Theorem~\ref{thm:surgery-bound} is an immediate consequence of the following
two propositions.

\begin{prop}\label{prop:crossing-change-bound} 
  Let $K,K'$ be knots in $S^3$ such that $K'$ is obtained from $K$ by changing
  a positive crossing. Let $Y,Y'$ be $(-1)$ surgeries on $K,K'$, respectively.
  Then
  \[0\le\qt(Y')-\qt(Y)\le1.\]
  \end{prop}

\proof We observe that $Y'$ is obtained from $Y$ by $(-1)$ surgery on a linking
circle $\ga$ of the crossing such that
$\ga$ bounds a surface in $Y$ of genus $1$.

The surgery cobordism $W$ from $Y$ to $Y'$ satisfies $H_1(W;\z)=0$ and
$b^+_2(W)=0$, hence $\qt(Y')\ge\qt(Y)$ by Theorem~\ref{thm:monotone}. Since $Y$
bounds a simply-connected negative definite $4$--manifold (the trace of the
surgery on $K$) we have $\qt(Y)\ge0$ by the same theorem.

Let $Y''$ be $0$--surgery on $\ga$.
By Floer's surgery theorem \cite{BD1,Scaduto1} there is a long exact sequence
\[\cdots\to I(Y'')\to I(Y)\oset\phi\to I(Y')\oset\psi\to I(Y'')\to\cdots\]
where $\phi$ is induced by the cobordism $W$.
Let $n:=\qt(Y')$ and suppose $n\ge2$, the proposition already being proved
for $n=0,1$. Then there is a $b\in I(Y')$ such that
\[\del u_2^jb=\begin{cases}
0,\quad 0\le j<n-1,\\
1,\quad j=n-1.
\end{cases}\]
By Proposition~\ref{prop:u-del--funct} we have
\[\psi u_2b=u_2\psi b=0,\]
hence $u_2b=\phi a$ for some $a\in I(Y)$. For $j\ge0$ we have
\[\del u_2^j a=\del u_2^j\phi a=\del u_2^{j+1} b.\]
Combining this with Corollary~\ref{cor-u3-vanishes-on-surgery} we
obtain $\qt(Y)\ge n-1=\qt(Y')-1$ and the proposition is proved.\square

\begin{prop}\label{prop:negpos} 
If $Y$ is $(-1)$ surgery on a positive knot $K$ in $S^3$ then $\qt(Y)=0$.
\end{prop}

\proof This follows from Theorem~\ref{thm:monotone} because $Y$ bounds 
simply-connected $4$--manifolds $V_\pm$ where $V_+$ is 
positive definite and $V_-$ is negative definite. As $V_-$ one can take the
trace of the $(-1)$ surgery on $K$. On the other hand, since $K$ can be
unknotted by changing a collection of positive crossings, the observation in
the beginning of the proof of Proposition~\ref{prop:crossing-change-bound}
yields $V_+$.\square

\subsection{Proof of Proposition~\ref{prop:Brieskorn-calc}}

Let $Y_k:=\Si(2,2k-1,4k-3)$. Then $Y_k$ bounds the simply-connected
$4$--manifold $V_k$ obtained by plumbing according the weighted graph
in Figure~$1$,

\begin{figure}[h]
   \caption{The plumbing graph for $\Si(2,2k-1,4k-3)$}
\centerline{\includegraphics*[width=0.6\textwidth]{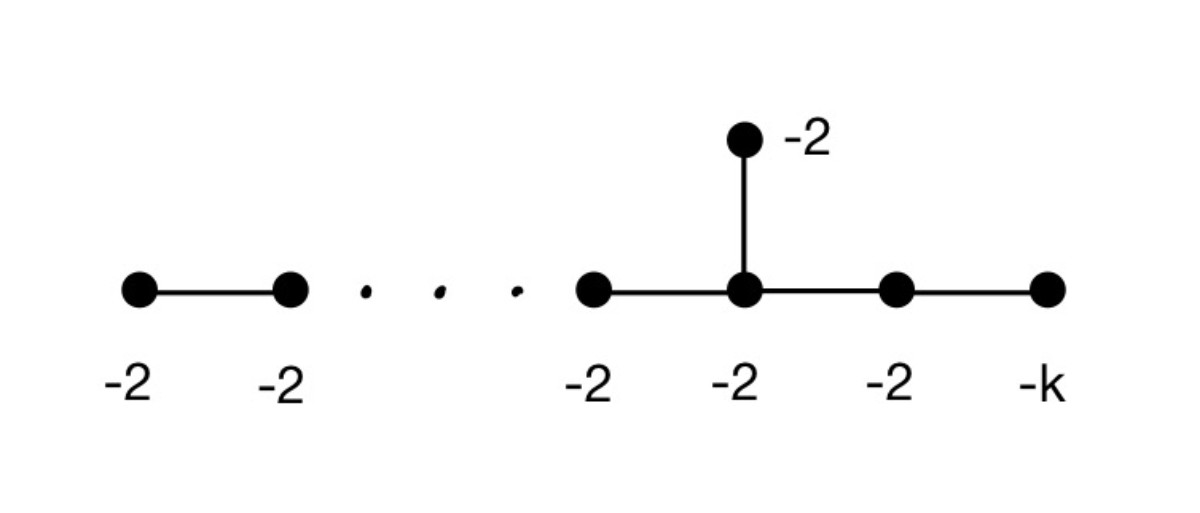}}
\end{figure}
\noindent where the total number of nodes is $4k$.
Let $e_1,\dots,e_{4k}$ be an orthonormal basis for $\R^{4k}$. The
intersection form of $V_k$ is isomorphic to the lattice
\[\Ga_{4k}:=
\left\{\sum_i x_ie_i\st
2x_i\in\z,\;x_i-x_j\in\z,\;\sum_i x_i\in2\z\right\},\]
with the nodes of the plumbing graph corresponding to the following
elements of $\Ga_{4k}$:
\[\frac12\sum_{i=1}^{4k}e_i,\quad e_2+e_3,\quad
(-1)^j(e_{j-1}-e_j),\;j=3,\dots,4k.\]
Let $w\in J_k=H^2(V_k;\z)$ be the element corresponding to
$\frac12\sum_{i=1}^{4k}e_i$. Since $\pm w$ are the only elements
of minimal square norm
in $w+2J_k$ it follows from Theorem~\ref{thm:defub-gen} that
\[\qt(Y_k)\ge k-1.\]
On the other hand, $Y_k$ is also the result of $(-1)$ surgery on the
torus knot $T_{2,2k-1}$. Since $T_{2,2k-1}$ can be unknotted by changing $k-1$
crossings we deduce from Theorem~\ref{thm:surgery-bound} that
\[\qt(Y_k)\le k-1.\]
This proves the proposition.\square

\subsection{Proof of Theorem~\ref{thm:h-qt-torsion}}

Since we will use different coefficient rings $R$, the homomorphism
\[\del:C^4(Y;R)\to R\]
defined in Subsection~\ref{subsec:inst-cohom} will now be denoted by
$\del_R$.

By definition, the condition $h(Y)>0$ means that there exists a cocycle
$w\in C^4(Y;\q)$ such that $\del_\q w\neq0$. Note that replacing the
coefficient group $\q$ by $\z$ yields an equivalent condition.

On the other hand, the condition $\qt(Y)>0$ means that there exists a
cocycle $z\in C^4(Y;\z/2)$ such that $\del_{\z/2}z\neq0$ and such that the
cohomology class of $z$ is annihilated by $u_3$. If in addition $z$ lifts
to an integral cocycle $\ti z\in C^4(Y;\z)$ then $\del_\z\ti z$ must be odd,
in particular non-zero, hence $h(Y)>0$.

Now suppose $\qt(Y)>0$ and $h(Y)\le0$. 
The above discussion shows that the homomorphism
$I^4(Y;\z)\to I^4(Y;\z/2)$ is not surjective, hence the Bockstein homomorphism
$I^4(Y;\z/2)\to I^5(Y;\z)$ is non-zero. This proves the theorem.\square

\subsection{Proofs of Theorems~\ref{thm:Poincare-sphere} and
  \ref{thm:sisi}}\label{subsec:proof-sisi}

{\em Proof of Theorem~\ref{thm:Poincare-sphere}:}
Part~(i) was proved in \cite{Fr1,Fr4} using Seiberg-Witten
theory. To prove (ii), let $\Si=\Si(2,3,5)$. Then $\qt(\Si)=1$ by
Proposition~\ref{prop:Brieskorn-calc}. If $H^2(X;\z)$ contains no $2$--torsion
then (ii) follows from Corollary~\ref{cor:JtiJ}. Under the weaker assumption
that $H^2(X;\z)$ contains no element of order $4$, we can appeal to
Theorem~\ref{thm:defub-gen} since $u_3=0$ on $I(\Si)$.\square

{\em Proof of Theorem~\ref{thm:sisi}:}
Let $\bfh$ be the 
monopole $h$--invariant defined in \cite{Fr4}. (One could
equally well use the correction term $d$.) Then $\bfh(\Si)=-1$, and
additivity of $\bfh$ yields $\bfh(\Si\#\Si)=-2$. If $\xi$ is any
characteristic vector for $J_X$ then by \cite[Theorem~4]{Fr4} one has
\[-\bfh(Y)\ge\frac18(b_2(X)+\xi\cdot\xi).\]
Let $J_X=m\la-1\ra\oplus\ti J_X$ as in Corollary~\ref{cor:JtiJ}.  By
assumption, $\ti J_X$ is even, so
$J_X$ has characteristic vectors $\xi$ with $\xi\cdot\xi=-m$. Therefore,
\[\rank\,\ti J_X=b_2(X)-m\le16.\]
By the classification
of even unimodular definite forms of rank~$\le16$ (see \cite{HM}) one has
\[\text{$\ti J_X=0$, $-E_8$, $-2E_8$, or $-\Ga_{16}$.}\]
It only remains to rule out $\ti J_X=-\Ga_{16}$.
Recalling that $\Si$ is the result of
$(-1)$ surgery on the negative trefoil knot and applying
Proposition~\ref{prop:u3u2nzero}
twice we find that $u_2^2=0$ on $I^*(\Si\#\Si)$, hence
$\qt(\Si\#\Si)\le2$. On the other hand,
if $\ti J_X=-\Ga_{16}$ then applying Theorem~\ref{thm:defub-gen} as in the proof
of Proposition~\ref{prop:Brieskorn-calc} we would obtain
$\qt(\Si\#\Si)\ge3$, a contradiction. This proves the theorem.\square

\section{Two points moving on a cylinder, I}\label{sec:two-points-I}

The main goal of this section is to prove
Proposition~\ref{prop:v2v3phi}. The first two subsections
will introduce some concepts used in the proof, which appears in the final
subsection.


\subsection{Energy and holonomy}
\label{subsec:enhol}

Let $Y$ be an oriented (integral) homology $3$--sphere with base-point $y_0$.
Let
\[\bbe\to\cb^*(Y[0])\]
be the canonical oriented Euclidean $3$--plane bundle, where
$Y[0]=[-1,1]\times Y$ as in \refp{eqn:ybt-def}.

Let $\al,\beta\in\fl(Y)$, not both reducible. Over $M(\al,\beta)\times\R$ there
is a canonical $3$--plane bundle $\bbei\al\beta$
obtained by pulling back the universal bundle over
$M(\al,\beta)\times\R\times Y$ by the map $(\om,t)\mapsto(\om,t,y_0)$.
There is a canonical isomorphism $\bbei\al\beta\to R^*\bbe$ where
\begin{equation}\label{eqn:Rproj}
R:M(\al,\beta)\times\R\to\cb^*(\yb0),\quad(\om,t)\mapsto\om[t],
\end{equation}
so we can identify the fibre of $\bbei\al\beta$ at $(\om,t)$ with
the fibre $\bbe_{\om[t]}$ of $\bbe$ at $\om[t]$.

Recall from Subsection~\ref{subsec:cutting-down-mod-sp}
that a section $\si$ of $\bbei\al\beta$ is called holonomy invariant if
for all $\om=[A]\in\mab$ and real numbers $s<t$ one has that $\si(\om,s)$
is mapped to $\si(\om,t)$ by the isomorphism
\be{equation*}
\bbe_{\om[s]}\to\bbe_{\om[t]}.
\end{equation*}
defined by holonomy of $A$ along the path $[s,t]\times\{y_0\}$.

Let $\ufl$ be the set of elements of $\cb^*(\yb0)$ that can be
represented by flat connections. 
Choose three sections $\rho_1,\rho_2,\rho_3$ of $\bbe$ which form a positive
orthonormal basis at every point in some neighbourhood of $\ufl$.
Choose $\eps>0$ so small that the following three conditions hold:
\be{description}
\item[(i)]If $A$ is any instanton over $(-\infty,2]\times Y$ satisfying
$\energy A{(-\infty,2]}<\eps$ such that the flat limit $\al$ of $A$ is
irreducible then $\rho_1,\rho_2,\rho_3$ are orthonormal at $A[0]$.
\item[(ii)]If $A$ is any instanton over $[-2,\infty)\times Y$ satisfying
$\energy A{[-2,\infty)}<\eps$ such that the flat limit $\beta$ of $A$ is
irreducible then $\rho_1,\rho_2,\rho_3$ are orthonormal at $A[0]$.
\item[(iii)]For each pair $\al,\beta\in\fl(Y)$ the difference
$\cs(\al)-\cs(\beta)\in\R/\z$ has no real lift in the half-open interval
$(0,2\eps]$.
  \end{description}
Here, $\ce_A$ refers to the energy of $A$ as defined in
\refp{eqn:def-energy}.

Let $\al,\beta$ be distinct elements of $\fl(Y)$. If $[A]\in M(\al,\beta)$ then
\[\ce_A(\R)>2\eps,\]
since the left hand side is a positive real lift of
$\cs(\al)-\cs(\beta)$. We can therefore define smooth functions
\[\tau^-,\tau^+:M(\al,\beta)\to\R\]
implicitly by
\[\ce_A((-\infty,\tau^-(A)+2])=\eps
=\ce_A([\tau^+(A)-2,\infty)).\]
We will consider the average and difference
\[\tau_a:=\frac12(\tau^++\tau^-),\quad\tau_d:=\tau^+-\tau^-.\]
Clearly, $\tau_d>0$.
There are translationary invariant smooth restriction maps
\[R^\pm:M(\al,\beta)\to\cb^*(\yb0),\quad\om\mapsto\om[\tau^\pm(\om)]\]
which, by the unique continuation result of
Proposition~\refp{prop:unique-continuation-cylinder}, descend to injective maps
$\check R^\pm:\chm(\al,\beta)\to\cb^*(\yb0)$. 

If $\al$ is irreducible then for any $\om=[A]\in M(\al,\beta)$ the vectors
\be{equation}
\rho_i(R^-(\om)),\quad i=1,2,3
\end{equation}
form an orthonormal basis for $\bbe_{R^-(\om)}$, by choice of $\eps$.
Let $\rho^-_i$ be the holonomy invariant section of $\bbei\al\beta$ whose
value at $(\om,\tau^-(\om))$ is $\rho_i(R^-(\om))$.

Similarly, if $\beta$ is irreducible, then the vectors
$\rho_i(R^+(\om))$ form an orthonormal basis for $\bbe_{R^+(\om)}$.
Let $\rho^+_i$ be the holonomy invariant section of $\bbei\al\beta$ whose
value at $(\om,\tau^+(\om))$ is $\rho_i(R^+(\om))$.

If $\al,\beta$ are both irreducible let
\[h=(h_{ij}):M(\al,\beta)\to\SO3\]
be the map whose value at $[A]$ is the holonomy of $A$ along
$[\tau^-(A),\tau^+(A)]\times\{y_0\}$ with respect to the bases described above,
so that
\[\rho^-_j(\om,t)=\sum_ih_{ij}(\om)\rho^+_i(\om,t).\]

\subsection{Factorization through the trivial connection}
\label{subsec:facttrivconn}

Now assume $\ind(\al)=4$, $\ind(\beta)=1$. We will introduce real valued
functions $\lla^\pm$ on $M(\al,\beta)$ which measure the extent
to which a given element factors through the trivial connection over $Y$.
Set
\[M_{\al,\theta}:=R^-(M(\al,\theta)),\]
which is a finite subset of $\cb^*(\yb0)$.
Let $M_\al$ be the union of all subsets
$R^-(M(\al,\beta'))\subset\cb^*(\yb0)$ where $\beta'\in\fl^*(Y)$ and
$\dim\,M(\al,\beta')\le4$. Note that $M_\al$ is compact.
Choose an open neighbourhood $U_\al$ of $M_{\al,\theta}$ in $\cb^*(\yb0)$ such that
\be{itemize}
\item the closure of $U_\al$ is disjoint from $M_\al$,
\item $U_\al$ is the disjoint union of open sets
$U_{\al,i}$, $i=1,\dots,r$,  each of which
contains exactly one point from $M_{\al,\theta}$.
\end{itemize}

Choose a closed neighbourhood $U'_\al$ of $M_{\al,\theta}$ contained in $U_\al$
and a smooth function
\be{equation}\label{eqn:eal}
e_\al:\cbs\to[0,\infty)
\end{equation}
such that $e_\al=1$ on $U'_\al$ and $e_\al=0$ outside $U_\al$. Define the
translationary invariant function
\[\lambda^-:M(\al,\beta)\to[0,\infty),\quad
\om\mapsto e_\al(R^-(\om))\cdot\tau_d(\om).\]

The function $\lla^+$ is defined in a symmetrical fashion (corresponding to
reversing the orientation of $Y$). 
Let $M_\beta$ be the union of all subsets
$R^+(M(\al',\beta))\subset\cb^*(\yb0)$ where $\al'\in\fl^*(Y)$ and
$\dim\,M(\al',\beta)\le4$.
Choose an open neighbourhood $V_\beta$ of
$M_{\theta,\beta}:=R^+(M(\theta,\beta)$ in $\cb^*(\yb0)$ such that
the closure of $V_\beta$ is disjoint from $M_\beta$, and such that
$V_\beta$ is the disjoint union of open sets
$V_{\beta,j}$, $j=1,\dots,s$, each of which
contains exactly one point from $M_{\theta,\beta}$.
Choose a
closed neighbourhood $V'_\beta$ of $M_{\theta,\beta}$ contained in $V_\beta$
and a smooth function
\[e_\beta:\cbs\to[0,\infty)\]
such that $e_\beta=1$ on $V'_\beta$ and $e_\beta=0$ outside $V_\beta$. Set
\[\lambda^+:M(\al,\beta)\to[0,\infty),\quad
\om\mapsto e_\beta(R^+(\om))\cdot\tau_d(\om).\]

\be{lemma}\label{lemma:llamin-llaplus}
There is a constant $C<\infty$ such that for any $\om\in M(\al,\beta)$ 
satisfying $\lla^-(\om)+\lla^+(\om)>C$ one has $\lla^-(\om)=\lla^+(\om)$.
\end{lemma}

\proof Suppose the lemma does not hold. Then one can find a sequence $\om_n$
in $M(\al,\beta)$ such that $\lla^-(\om_n)+\lla^+(\om_n)\to\infty$ and
$\lla^-(\om_n)\neq\lla^+(\om_n)$. After passing to a subsequence we may assume
that the sequence $\om_n$ chain-converges. If the chain-limit lay in
$\chm(\al,\beta)$, or if the chain-limit involved factorization through
an irreducible critical point, then $\lla^\pm(\om_n)$ would be bounded.
Therefore, the chain-limit must lie in
$\chm(\al,\theta)\times\chm(\theta,\beta)$ and, consequently,
$\lla^-(\om_n)=\tau_d(\om_n)=\lla^+(\om_n)$ for $n\gg0$, a contradiction.\square

In the course of the proof we also obtained the following:

\be{lemma}\label{lemma:delinfty}
For a chain-convergent sequence $\om_n$ in $M(\al,\beta)$ the following are
equivalent:
\be{description}
\item[(i)] $\lambda^-(\om_n)\to\infty$.
\item[(ii)]  $\lambda^+(\om_n)\to\infty$.
\item[(iii)] The chain-limit of $\om_n$ lies in
$\chm(\al,\theta)\times\chm(\theta,\beta)$.\square
\end{description}
\end{lemma}

Since $\lla^+$ will not appear again in the text, we set
\[\lla:=\lla^-\]
to simplify notation. For any real number $T$ set
\[\mab_{\lla=T}:=\{\om\in\mab\st\lla(\om)=T\}.\]

Given $\om\in M(\al,\beta)$, one has $R^-(\om)\in U_\al$ if $\lla(\om)>0$
(by definition of $\lla$), and $R^+(\om)\in V_\beta$ if $\lla(\om)\gg0$
(by Lemma~\ref{lemma:delinfty}).
Therefore, if $\lla(\om)\gg0$ then there is a map
\[d:M(\al,\beta)_{\lla=T}\to\chm(\al,\theta)\times\chm(\theta,\beta)\]
characterized by the fact that if $d(\om)=(\om_1,\om_2)$ then
$R^-(\om)$ and $\check R^-(\om_1)$ lie in the same set $U_{\al,i}$, and
$R^+(\om)$ and $\check R^+(\om_2)$ lie in the same set $V_{\beta,j}$.

Gluing theory (see \cite{D5,Fr13}) provides the following result:
\be{lemma}\label{lemma:glthm}
There is a $T_0>0$ such that for any $T\ge T_0$ the map
\[d\times h\times\tau_a:
\mab_{\lla=T}\to(\chm(\al,\theta)\times\chm(\theta,\beta))\times\SO3\times\R\]
is a diffeomorphism.\square
\end{lemma}

\subsection{Proof of Proposition~\ref{prop:v2v3phi}}
\label{subsec:proof-v2v3phi}

Let $\al,\beta\in\fl^*(Y)$ with
$\ind(\beta)-\ind(\al)\equiv5\mod8$. To compute the
matrix coefficient $\la(v_2v_3+v_3v_2)\al,\beta\ra$ we distinguish between two
cases. If $\ind(\al)\not\equiv4\mod8$ the calculation will consist in counting
modulo~$2$ the number of ends of the $1$-manifold $\mx23(\al,\beta)$.
If $\ind(\al)\equiv4\mod8$ then $M(\al,\beta)$ may contain
sequences factoring through the trivial connection over $Y$. To deal
with this we consider the subspace of
$M(\al,\beta)\times\R$ consisting of points $(\om,t)$ with
$\lla(\om)\le T$ for some large 
$T$. By carefully cutting down this subspace to a $1$-manifold and then
counting the number of ends and boundary points modulo~$2$ we obtain
\refp{eqn:v2v3chhom}.

For $s\in\R$ we define the translation map
\[\ct_s:\ry\to\ry,\quad(t,y)\mapsto(t+s,y).\]

{\bf Part (I)} Suppose $\ind(\al)\not\equiv4\mod8$. Then
no sequence in $M(\al,\beta)$ can have a chain-limit involving
factorization through the trivial connection.
We will determine the ends of the smooth $1$-manifold $\mx23(\al,\beta)$
introduced in Definition~\ref{defn:M23}.
Let $(\om_n,t_n)$ be a sequence in
$\mx23(\al,\beta)$. After passing to a subsequence we may assume that
the following hold:
\be{description}
\item[(i)] The sequence $\ct^*_{-t_n}(\om_n)$
  converges over compact subsets of $\ry$ to some
  $\om^-\in M(\al^-,\beta^-)$. (By this we mean
  that there are connections
$A_n,\bar A$ representing $\om_n,\om^-$ respectively, such that
$A_n\to \bar A$ in $C^\infty$ over compact subsets of $\ry$.)
\item[(ii)] The sequence $\ct^*_{t_n}(\om_n)$ converges over compact subsets of
$\ry$ to some $\om^+\in M(\al^+,\beta^+)$.
\item[(iii)] The sequence $t_n$ converges in $[-\infty,\infty]$ to some point
$t_\infty$.
\end{description}
Here, $[-\infty,\infty]$ denotes the compactification of the real line obtained
by adding two points $\pm\infty$.

Suppose $(\om_n,t_n)$ does not converge in $\mx23(\al,\beta)$.

{\em Case 1:} $t_\infty$ is finite. Then $M(\al^-,\beta^-)$ has
dimension $4$ and either $\al^-=\al$ or $\beta^-=\beta$. The corresponding
number of ends of $\mx23(\al,\beta)$, counted modulo~$2$, is
\[\la(d\phi+\phi d)\al,\beta\ra.\]

{\em Case 2:} $t_\infty=\infty$. Let $n^\pm$ be the dimension of
$M(\al^\pm,\beta^\pm)$. Because
\[s_1(\om^-[0])=0,\quad s_2(\om^+[0])\wedge s_3(\om^+[0])=0\]
we must have $n^-\ge3$ and $n^+\ge2$. On the other hand,
\[n^-+n^+\le\dim M(\al,\beta)=5,\]
so $n^-=3$, $n^+=2$. It follows that
\[\al=\al^-,\;\beta^-=\al^+,\;\beta^+=\beta.\]
The corresponding number of ends of $\mx23(\al,\beta)$ is
$\la v_2v_3\al,\beta\ra$ modulo~$2$.

{\em Case 3:} $t_\infty=-\infty$. Arguing as in Case~2 one finds that the number
of such ends of $\mx23(\al,\beta)$ is
$\la v_3v_2\al,\beta\ra$ modulo~$2$.

Since the total number of ends of $\mx23(\al,\beta)$ must be zero modulo~$2$,
we obtain the equation~\refp{eqn:v2v3chhom} in the case
$\ind(\al)\not\equiv4\mod8$.

{\bf Part (II)} Now suppose $\ind(\al)\equiv4\mod8$. In this case, the
$1$--manifold $\mx23(\al,\beta)$ may have additional ends corresponding to
factorization through the trivial connection. Instead of attempting to
count these ends directly, we will replace $\mx23(\al,\beta)$ with
another $1$--manifold $\tmx23(\al,\beta)$ defined in the same way as
$\mx23(\al,\beta)$ except that the equations for cutting down $M(\al,\beta)$
are deformed in part of the region of $M(\al,\beta)$ containing 
instantons that ``tend to'' factor through the trival connection.
We then cut off part of $\tmx23(\al,\beta)$ to obtain a
$1$--manifold-with-boundary $\mxl23(\al,\beta)$ in which factorization
through the trivial connection does not occur. Counting the ends and
boundary points of the latter manifold yields a proof of the proposition
in the case $\ind(\al)\equiv4\mod8$.

We will again make use of a cut-off function $b$ as in \refp{eqn:b-prop1} in
Subsection~\ref{subsec:2-torsion-inv},
but we now impose two further conditions, namely
\begin{equation}\label{eqn:b-prop2}
b(0)=\frac12,\quad\text{$b'(t)>0$ for $-1<t<1$.}
\end{equation}
Set
\begin{equation}\label{eqn:c23def}
c:\mab\times\R\to\R,\quad(\om,t)\mapsto b(t-\tau_a(\om)).
\end{equation}
Choose generic $3\times3$ matrices $A^+=(a^+_{ij})$ and $A^-=(a^-_{ij})$ and
for $j=1,2,3$ define a section $\ti\rho_j$ of the bundle $R^*\bbe$
over $M(\al,\beta)\times\R$ by
\begin{equation}\label{eqn:rhorho}
\ti\rho_j:=(1-c)\sum_ia^-_{ij}\rho^-_i+c\sum_ia^+_{ij}\rho^+_i.
\end{equation}
Define a function $g:M(\al,\beta)\times\R\to[0,1]$ by
\begin{equation}\label{eqn:gomt}
  g(\om,t):=b(\lla(\om)-1)\cdot b(\tau^+(\om)-t)\cdot b(t-\tau^-(\om)).
  \end{equation}

For $j=1,2,3$ we now define a section $\ti s_j$ of $R^*\bbe$ by
\[\ti s_j(\om,t):=(1-g(\om,t))\cdot s_j(\om[t])+g(\om,t)\cdot\ti\rho_j(\om,t).\]

\be{defn}
Let $\tmx23(\al,\beta)$ be the subspace of $\mab\times\R$ consisting of those
points $(\om,t)$ that satisfy the following conditions:
\be{itemize}
\item $\ti s_1(\om,-t)=0$,
\item $\ti s_2(\om,t)$ and $\ti s_3(\om,t)$ are linearly dependent.
\end{itemize}
\end{defn}

To understand the ends of $\tmx23(\al,\beta)$
we will need to know that certain subspaces of $M(\al,\theta)$ and
$M(\theta,\beta)$, respectively, are ``generically'' empty.
These subspaces are defined as
follows. For $\om\in M(\al,\theta)$ and $j=1,2,3$ let
\[\ti s_j(\om):=(1-b(-\tau^-(\om)))\cdot s_j(\om[0])+b(-\tau^-(\om))
\sum_ia^-_{ij}\rho^-_i(\om,0),\]
and for $\om\in M(\theta,\beta)$ let
\[\ti s_j(\om):=(1-b(\tau^+(\om)))\cdot s_j(\om[0])+b(\tau^+(\om))
\sum_ia^+_{ij}\rho^+_i(\om,0).\]
Set
\begin{align*}
\ti M_2(\al,\theta)&:=\{\om\in M(\al,\theta)\st\ti s_2(\om)\wedge\ti s_3(\om)=0\},\\
\ti M_3(\al,\theta)&:=\{\om\in M(\al,\theta)\st\ti s_1(\om)=0\}.
\end{align*}
Replacing $(\al,\theta)$ by $(\theta,\beta)$ in the last two definitions
we obtain subspaces $\ti M_k(\theta,\beta)$ of $M(\theta,\beta)$.
For $k=2,3$, each of the spaces $\ti M_k(\al,\theta)$ and $\ti M_k(\theta,\beta)$
has expected dimension
$1-k$ and is therefore empty for ``generic'' choices of sections $s_j$ and
matrices $A^\pm$. 

\begin{lemma}\label{lemma:ttauest1}
  There is a constant $C_0<\infty$ such that for all
  $(\om,t)\in\tmx23(\al,\beta)$ one has
  \[|t|\le\min(-\tau^-(\om),\tau^+(\om))+C_0.\] 
\end{lemma}

\proof We must prove that both quantities $|t|+\tau^-(\om)$ and
$|t|-\tau^+(\om)$ are uniformly bounded above for $(\om,t)\in\tmx23(\al,\beta)$.
The proof is essentially the same in both cases, so we will only spell it out
in the first case. Suppose, for contradiction, that $(\om_n,t_n)$
is a sequence in $\tmx23(\al,\beta)$
with $|t_n|+\tau^-(\om_n)\to\infty$.
After passing to a subsequence we may assume
that the sign of $t_n$ is constant, so $|t_n|=-et_n$ for some constant
$e=\pm1$. Then $\om[et_n]\to\ual$ by exponential decay
(see \cite[Subsection~4.1]{D5}), and
\[\ti s_j(\om,et_n)=s_j(\om_n[et_n])\quad\text{for $n\gg0$.}\]
If $e=1$ then this gives
\[0=s_2(\om_n[t_n])\wedge s_3(\om_n[t_n])\to s_2(\ual)\wedge s_3(\ual),\]
as $n\to\infty$, whereas if $e=-1$ we get
\[0=s_1(\om_n[-t_n])\to s_1(\ual).\]
However, for ``generic'' sections $s_j$, both $s_2(\ual)\wedge s_3(\ual)$
and $s_1(\ual)$ are non-zero. This contradiction proves the lemma.
\square

\begin{lemma}\label{lemma:ttauest2}
  For any constant $C_1<\infty$ there is constant $L>0$ such that for
  all $(\om,t)\in\tmx23(\al,\beta)$ satisfying $\lla(\om)\ge L$ one has
  \[|t|\le\min(-\tau^-(\om),\tau^+(\om))-C_1.\]
\end{lemma}

\proof Suppose to the contrary that there is a constant $C_1<\infty$ and a
sequence $(\om_n,t_n)$ in $\tmx23(\al,\beta)$ such that $\lla(\om_n)\to\infty$
and
\[|t_n|>\min(-\tau^-(\om_n),\tau^+(\om_n))-C_1.\]
After passing to a subsequence we may assume that at least one of the following
two conditions holds:
\begin{description}
\item[(i)] $|t_n|>-\tau^-(\om_n)-C_1$ for all $n$,
  \item[(ii)] $|t_n|>\tau^+(\om_n)-C_1$ for all $n$.
\end{description}
The argument is essentially the same in both cases, so suppose (i) holds. By
Lemma~\ref{lemma:ttauest1} we also have
\[|t_n|\le-\tau^-(\om_n)+C_0,\]
hence the sequence $\tau^-(\om_n)+|t_n|$ is bounded. Since
$\lla(\om_n)\to\infty$ we have $\tau_d(\om_n)\to\infty$, so
\[\tau^+(\om_n)+|t_n|=\tau_d(\om_n)+(\tau^-(\om_n)+|t_n|)\to\infty.\]
After passing to a subsequence we may assume that
\begin{itemize}
\item the sequence $\om_n$ chain-converges;
\item the sequence $\tau^-(\om_n)+|t_n|$ converges to a real number;
  \item $|t_n|=-et_n$ for some constant $e=\pm1$.
  \end{itemize}
From Lemma~\ref{lemma:delinfty} we deduce that $\om'_n:=\ct^*_{et_n}\om_n$
converges over compact subsets of $\ry$ to
some $\om\in M(\al,\theta)$. For large $n$ we have $c(\om_n,et_n)=0$
and
\[g(\om_n,et_n)=b(et_n-\tau^-(\om_n))=b(-\tau^-(\om'_n))\to b(-\tau^-(\om)).\]
For $j=1,2,3$ we now get
\[\ti s_j(\om_n,et_n)\to \ti s_j(\om).\]
But then $\om$ lies in
$\ti M_2(\al,\theta)$ (if $e=1$) or in $\ti M_3(\al,\theta)$ (if $e=-1$),
contradicting the fact that the latter two spaces are empty.\square

Choose $L\ge2$ such that for all $(\om,t)\in\tmx23(\al,\beta)$ with
$\lla(\om)\ge L$ one has
\[|t|\le\min(-\tau^-(\om),\tau^+(\om))-1,\]
which implies that $\ti s_j(\om,t)=\ti\rho_j(\om,t)$. Set
\[\mxl23(\al,\beta):=\{(\om,t)\in\tmx23(\al,\beta)\st\lla(\om)\le L\}.\]

We will show that $\mxl23(\al,\beta)$ is transversely cut and therefore
a one-manifold with boundary, and determine the number of boundary
points and ends modulo~$2$. We will see that the number of ends is given by
the same formula as in Part~(I), whereas the boundary points contribute the
new term $\del'\del$ of \refp{eqn:v2v3chhom}.

{\bf Ends of $\mxl23(\al,\beta)$:} 
Let $(\om_n,t_n)$ be a sequence in
$\mxl23(\al,\beta)$. After passing to a subsequence we may assume that
(i),(ii), (iii) of Part~(I) as well as the following hold:

\be{description}
\item[(iv)] The sequence $\om_n$ is chain-convergent.
\item[(v)] The sequence $\tau_a(\om_n)$ converges in $[-\infty,\infty]$.
\item[(vi)] Either $\lla(\om_n)>0$ for all $n$, or $\lla(\om_n)=0$ for all $n$.
\end{description}

Suppose $(\om_n,t_n)$ does not converge in $\mxl23(\al,\beta)$.

{\em Case 1:} $\lla(\om_n)=0$ for all $n$. Then $g(\om_n,t_n)=0$ and therefore
\[\ti s_j(\om_n,t_n)=s_j(\om_n[t_n]).\]
This case is similar to Part~(I) and the corresponding number of ends of
$\mxl23(\al,\beta)$, counted modulo~$2$, is
\[\la(v_2v_3+v_3v_2+d\phi+\phi d)\al,\beta\ra,\]
where $\phi$ is defined as before.

{\em Case 2:} $\lla(\om_n)>0$ for all $n$. We show this is impossible.
By definition of $\lla$ the
chain-limit of $\om_n$ must lie in $\chm(\al,\beta)$, so
$\tau_d(\om_n)$ is bounded. By Lemma~\ref{lemma:ttauest1}, the sequence
$\tau^-(\om_n)$ is bounded above whereas $\tau^+(\om_n)$ is bounded below,
hence both sequences must be bounded.
Applying Lemma~\ref{lemma:ttauest1}
again we see that $t_n$ is bounded. Therefore, both sequences
$\tau_a(\om_n)$ and $t_n$ converge in $\R$, so $(\om_n,t_n)$ converges
in $M(\al,\beta)\times\R$ and hence in $\mxl23(\al,\beta)$,
which we assumed was not the case.

{\bf Boundary points of $\mxl23(\al,\beta)$:} Let $M=M(3,\R)$ be the space of
all $3\times3$ real matrices, and let $U\subset M$ be the open subset
consisting of those matrices $B$ satisfying
\[B_1\neq0,\quad B_2\wedge B_3\neq0,\]
where $B_j$ denotes the $j$th column of $B$. Then $M\setminus U$ is the union
of three submanifolds of codimension at least two, hence $U$ is a connected
subspace and a dense subset of $M$. Let
\begin{align*}
  F:\SO3\times\R\times\R\times U\times U&\to\R^3\times\R^3\times\R^3,\\
(H,v,w,B^+,B^-)&\mapsto(F_1,F_2,F_3),
\end{align*}
where
\begin{align*}
  F_1&=(1-b(v))HB^-_1+b(v)B^+_1,\\
  F_j&=(1-b(w))HB^-_j+b(w)B^+_j,\quad j=2,3.
  \end{align*}
Then $F$ is a submersion, so $F\inv(0,0,0)$ is empty. Moreover, the set
\[Z:=F\inv(\{0\}\times L(\R^3),\]
consisting of those points in the domain of $F$ for which
\begin{equation}\label{eqn:FFF}
  F_1=0,\quad F_2\wedge F_3=0,
  \end{equation}
is a codimension $5$ submanifold and a closed subset of
$\SO3\times\R^2\times U^2$.

\begin{claim}
  The projection $\pi:Z\to U^2$ is a proper map whose mod~$2$ degree is
  \[\deg_2(\pi)=1.\]
  \end{claim}

\proof The equations \refp{eqn:FFF} imply $-1<v,w<1$, hence $\pi$ is
proper. To compute its degree,
let $e_1,e_2,e_3$ be the standard basis for $\R^3$ and let $B^\pm$ be given by
\begin{gather*}
  B^-_1=B^-_2=e_1,\quad B^-_3=e_2,\\
  B^+_1=-e_1,\quad B^+_2=e_1,\quad B^+_3=-e_2.
\end{gather*}
We show that the preimage
$Z':=\pi\inv(B^+,B^-)$ consists of precisely one point.
Suppose $(H,v,w)\in Z'$. Because $0\le b\le1$, the equation $F_1=0$ implies
$b(v)=1/2$ and hence $v=0$, $He_1=e_1$, $F_2=e_1$. Because
$He_2\perp e_1$, the vectors $F_2,F_3$ are linearly dependent if and only if
$F_3=0$, which yields $w=0$, $He_2=e_2$. Thus,
\[Z'=\{(I,0,0)\},\]
where $I$ is the identity matrix. 
Using the fact that $f(I,0,0)=(0,e_1,0)$
and that the tangent space to $L^*(\R^3)$ at $(e_1,0)$ is
$\R^3\times\{0\}+\R e_1$ it is easy to see that the map
\[F(\:\cdot\:,\:\cdot\:,\:\cdot\:,B^+,B^-):\SO3\times\R\times\R\to\R^9\]
is transverse to
$\{0\}\times L^*(\R^3)$ at $(I,0,0)$, or equivalently, that $(B^+,B^-)$
is a regular value of $\pi$. This proves the claim.\square

By Lemma~\ref{lemma:glthm} we can identify
\[\partial\mxl23(\al,\beta)=
\chm(\al,\theta)\times\chm(\theta,\beta)\times\pi\inv(A^+,A^-),\]
where $(H,v,w)$ corresponds to $(h(\om),-t-\tau_a(\om),t-\tau_a(\om))$ for
$(\om,t)\in\partial\mxl23(\al,\beta)$.
Hence, for generic matrices $A^\pm$ the number of boundary points of
$\mxl23(\al,\beta)$, counted modulo~$2$, is $\la\del'\del\al,\beta\ra$.

This completes the proof of Proposition~\ref{prop:v2v3phi}.
\square

\section{Two points moving on a cylinder, II}\label{sec:two-points-II}

Let $Y$ be an oriented homology $3$--sphere.
In this section we will prove Proposition~\ref{prop:psixi}, which concerns
a certain cochain map
\[\psi:C^*(Y)\to C^{*+5}(Y)\]
appearing in the proof of additivity of $\qt$.
We will continue using the notation
introduced in Section~\ref{sec:two-points-I}.

\subsection{The cochain map $\psi$}\label{subsec:cochain-map-psi}

We begin by recalling the definition of $\psi$ in degrees
different from $4$ mod~$8$ given in Subsection~\ref{subsec:comm-cup-products}.
Let $s_1,s_2$ be "generic" sections of the canonical $3$--plane bundle
\[\bbe\to\cb^*(Y[0]).\]
(Later we will impose further conditions on $s_1,s_2$.)
For any $\al,\beta\in\fl^*(Y)$ set
\[\mx33(\al,\beta):=\{(\om,t)\in\mab\times\R\st s_1(\om[-t])=0=s_2(\om[t])\}.\]
If $\ind(\al,\beta)=5$ and $\ind(\al)\not\equiv4\mod8$ then
arguing as in Part~(I) of the proof of Proposition~\ref{prop:v2v3phi}
one finds that
$\mx33(\al,\beta)$ is a finite set. We define the matrix coefficient
$\la\psi\al,\beta\ra$ by
\[\la\psi\al,\beta\ra:=\#\mx33(\al,\beta).\]
Recall that any "generic" section of $\bbe$ defines a cup product
$C^*(Y)\to C^{*+3}(Y)$ by the formula \refp{eqn:v3def}. Let $v_3$ and $v'_3$
be the cup products defined by $s_1$ and $s_2$, respectively.

\begin{prop}\label{prop:v3v3-restated}
For $q\not\equiv3,4\mod8$ one has
\[d\psi+\psi d=v_3v'_3+v'_3v_3\]
as maps $C^q(Y)\to C^{q+6}(Y)$.
\end{prop}

\proof Let $\al,\ga\in\fl^*(Y)$ with $\ind(\al,\ga)=6$ and
$\ind(\al)\not\equiv3,4\mod8$. Note that no sequence in $M(\al,\ga)$ can
have a chain-limit involving factorization through the trivial connection.
Now let $(\om_n,t_n)$ be a sequence
in $\mx33(\al,\ga)$. After passing to a subsequence we may assume that
\be{description}
\item[(i)] The sequence $\ct^*_{t_n}\om_n$ converges over compact subsets of
$\ry$ to some point $\om^+\in M(\al^+,\ga^+)$.
\item[(ii)] The sequence $\ct^*_{-t_n}\om_n$ converges over compact subsets of
$\ry$ to some point $\om^-\in M(\al^-,\ga^-)$.
\item[(iii)] The sequence $t_n$ converges in $[-\infty,\infty]$ to some point
$t_\infty$.
\end{description}

Clearly, $s_1(\om^+[0]=0=s_2(\om^-[0])$, hence $\ind(\al^\pm,\ga^\pm)\ge3$.

{\em Case 1:} $t_\infty$ finite. Then $\ind(\al^+,\ga^+)=5$ and either
$\al^+=\al$ or $\ga^+=\ga$. The corresponding number of ends of
$\mx33(\al,\ga)$, counted modulo $2$, is
\[\la(d\psi+\psi d)\al,\ga\ra.\]

{\em Case 2:} $t_\infty=\infty$. Then $\ind(\al^\pm,\ga^\pm)=3$, so
$\al^-=\al$, $\ga^-=\al^+$, and $\ga^+=\ga$. The corresponding number of ends of
$\mx33(\al,\ga)$ is $\la v_3v'_3\al,\ga\ra$ modulo $2$.

{\em Case 3:} $t_\infty=-\infty$. As in Case~2 one finds that the number of such
ends is $\la v'_3v_3\al,\ga\ra$ modulo $2$.

Since the total number of ends of $\mx33(\al,\ga)$ must be zero modulo $2$,
we obtain the proposition.\square

We now show that $v_3=v'_3$ if the sections $s_1,s_2$
are close enough in a certain
sense. To make this precise, we introduce the following
terminology: We will say a section $s$ of $\bbe$ has
{\em Property~$\pt4$} if for all
$\al,\beta\in\fl^*(Y)$ with $\ind(\al,\beta)\le4$ the map
\[s_{\al\beta}:M(\al,\beta)\to\bbe,\quad\om\mapsto s(\om[0])\]
is transverse to the zero-section in $\bbe$.

\begin{lemma}\label{lemma:prop-T}
Suppose $s\in\Ga(\bbe)$ has Property~$\pt4$, and let $\frP$ be any
finite-dimensional linear
subspace of $\Ga(\bbe)$. Then for any sufficiently small $\frp\in\frP$ 
the following hold:
\be{description}
\item[(i)]The section $s':=s+\frp$ has Property~$\pt4$.
\item[(ii)]The sections $s$ and $s'$ define the same cup product
$C^*(Y)\to C^{*+3}(Y)$.
\end{description}
\end{lemma}

\proof Let $\ind(\al,\beta)=3$.
Combining the transversality assumption with a compactness argument
one finds that the zero-set $Z$ of $s_{\al\beta}$ is a finite set.
Now observe that the map
\be{equation}\label{eqn:sfrpmap}
M(\al,\beta)\times\frP\to\bbe,\quad(\om,\frp)\mapsto(s+\frp)(\om[0])
\end{equation}
is smooth, since $\frP$ has finite dimension. Therefore, given any
neighbourhood $U$ of $Z$ in $M(\al,\beta)$ then the zero-set of
$(s+\frp)_{\al\beta}$ is contained in $U$ for all sufficiently small $\frp$.
The lemma now follows by applying the implicit function theorem to the map
\refp{eqn:sfrpmap}.\square

From now on we assume that $s_1,s_2$ are sufficiently close in the sense of the
lemma, so that in particular $v_3=v'_3$. Since we are taking coefficients 
in $\z/2$, we deduce from Proposition~\ref{prop:v3v3} that $d\psi=\psi d$
in degrees different from $3$ and $4$ modulo $8$.

We now extend the definition of $\psi$ to degree~$4$.
Let $\al,\beta\in\fl^*(Y)$ with
$\ind(\al)=4$ and $\ind(\beta)=1$. To define the matrix coefficient
$\la\psi\al,\beta\ra$ we use
the set-up of Subsections~\ref{subsec:enhol} and \ref{subsec:facttrivconn}
and define $\ti\rho_j,\ti s_j$ for $j=1,2$
as in Subsection~\ref{subsec:proof-v2v3phi}, where $A^\pm$ should now be
generic $3\times2$ real matrices. In particular, we require that
$A^\pm$ should have non-zero columns and that the angle between the columns
of $A^+$ should be different from the angle between the columns
of $A^-$. For any $3\times2$ real matrix $B$ with non-zero columns $B_j$
set
\begin{equation}\label{eqn:nuB}
  \nu(B):=\frac{\la B_1,B_2\ra}{\|B_1\|\|B_2\|},
  \end{equation}
using the standard scalar product and norm on $\R^3$. 
Then the above assumption on the angles means that $\nu(A^+)\neq\nu(A^-)$.
Now define
\[\tmx33(\al,\beta):=\{(\om,t)\in\mab\times\R\st
\ti s_1(\om,-t)=0,\;\ti s_2(\om,t)=0\}.\]

\be{prop}\label{prop:tmx-finite-set}
$\tmx33(\al,\beta)$ is a finite set.
\end{prop}

\proof It is easy to see that Lemmas~\ref{lemma:ttauest1}
and \ref{lemma:ttauest2} hold with $\tmx33(\al,\beta)$ in place of
$\tmx23(\al,\beta)$. Arguing as in Subsection~\ref{subsec:proof-v2v3phi}
one finds that for any $L>0$ there are only finitely many points
$(\om,t)\in\tmx33(\al,\beta)$ with $\lla(\om)\le L$. Choose $L\ge2$ such that
for all $(\om,t)\in\tmx33(\al,\beta)$ with $\lla(\om)\ge L$ one has
\[|t|\le\min(-\tau^-(\om),\tau^+(\om))-1,\]
which implies that $\ti s_j(\om,t)=\ti\rho_j(\om,t)$. We claim that
there are no such $(\om,t)$. For suppose $(\om,t)$
is such an element and set
\[(H,v_1,v_2):=(h(\om),-t-\tau_a(\om),t-\tau_a(\om))\in\SO3\times\R\times\R.\]
Then for $j=1,2$ one has
\begin{equation*}
  (1-b(v_j))HA^-_j+b(v_j)A^+_j=0.
\end{equation*}
However, there is no solution $(H,v_1,v_2)$ to these equations, since we
assume the columns $A^\pm_j$ are non-zero and $\nu(A^+)\neq\nu(A^-)$.\square

We define $\psi$ in degree $4$ by
\[\la\psi\al,\beta\ra:=\#\tmx33(\al,\beta).\]

\be{prop}\label{prop:psichainmap}
If the endomorphism $\psi$ is defined in terms of ``generic'' sections $s_1,s_2$
that are sufficiently close then
\[d\psi=\psi d\]
as maps $C^*(Y)\to C^{*+6}(Y)$.
\end{prop}

Although we could deduce this from Proposition~\ref{prop:psiv3v2} below,
we prefer to give a direct proof, partly because the techniques involved
are also needed in the proof of Proposition~\ref{prop:Fthm}.

It only remains to prove this in degrees $3$ and $4$ modulo $8$. There is a
complete symmetry between these two cases because of
Lemma~\ref{lemma:llamin-llaplus}, so we will spell out the proof only in
degree $4$. Let $\al,\ga\in\fl^*(Y)$ with $\ind(\al)=4$, $\ind(\ga)=2$.
We will show that $\la(d\psi+\psi d)\al,\ga\ra=0$ by counting the ends of
a certain $1$--dimensional submanifold
$\tmx33(\al,\ga)$ of $M(\al,\ga)\times\R$.

For any $\al'\in\fl(Y)$ we define a smooth function
\[\ttp:M(\al',\ga)\to\R\]
as follows. 
For each $\beta\in\fl^1_Y$ let $K_\beta$ be the union of all subsets
$R^+(M(\al'',\ga))\subset\cb^*(Y[0])$ where $\beta\neq\al''\in\fl(Y)$ and
\[\cs(\al'',\ga)\le\cs(\beta,\ga),\]
where $\cs(\,\cdot\,,\,\cdot\,)$ is as in \refp{eqn:cs-al-beta}.
Then $K_\beta$ is compact. Choose a closed neighbourhood
$W_\beta$ in $\cb^*(Y[0])$ of the finite set $R^+(M(\beta,\ga))$ such that
$W_\beta$ is disjoint from $K_\beta$, and a smooth function
\[f_\beta:\cb^*(Y[0])\to[0,1]\]
such that the following two conditions hold:
\begin{itemize}
\item $W_\beta$ and $W_{\beta'}$ are disjoint if $\beta\neq\beta'$;
\item $f_\beta=1$ on a neighbourhood of $R^+(M(\beta,\ga))$,
  and $f_\beta=0$ outside $W_{\beta}$.
\end{itemize}

Set $f:=1-\sum_\beta f_\beta$.
Let $\flag$ be the set of all $\beta\in\fl^1_Y$ such that
\[\cs(\al',\ga)>\cs(\beta,\ga)>0.\]
For $\om\in M(\al',\ga)$ and $\beta\in\flag$ we
define $\tau^+_\beta(\om)\in\R$ implicitly by
\[\ce_\om([\tau^+_\beta(\om)-2,\infty))=\cs(\beta,\ga)+\eps,\]
where the constant $\eps$ is as in Subsection~\ref{subsec:enhol}, and set
\[\ttp(\om):=f(R^+(\om))\cdot\tau^+(\om)+
\sum_\beta f_\beta(R^+(\om))\cdot\tau^+_\beta(\om).\]
The function $\ttp$ behaves under translation in the same way as
$\tau^\pm$. Namely, for any real number $s$ one has
\[\ttp(\ct^*_s(\om))=\ttp(\om)-s.\]
For any $\om\in M(\al',\ga)$ let 
$\hrp(\om)$ denote the restriction of $\om$ to the band $\yb{\ttp(\om)}$.
For $i=1,2,3$ let $\trp i$ be the holonomy invariant section of the bundle
$\bbei{\al'}\beta$ over $M(\al',\beta)\times\R$ (as defined in
Subsection~\ref{subsec:enhol}) whose value at $(\om,\ttp(\om))$ is $\rho_i(\hrp(\om))$.

\be{lemma}\label{lemma:ttpinfty}
Let $\om_n$ be a chain-convergent sequence  in $M(\al',\ga)$.
If the last term of the chain-limit of $\om_n$ lies in $\chm(\beta,\ga)$
for some $\beta\in\fl^*(Y)$ of index~$1$ then
\[(\tau^+-\ttp)(\om_n)\to\infty,\]
otherwise the sequence $(\tau^+-\ttp)(\om_n)$ is bounded.
\end{lemma}

\proof Because of the translationary invariance of $\tau^+-\ttp$ we may
assume that $\tau^+(\om_n)=0$. Then $\om_n$ converges over compact subsets of
$\ry$ to some element $\om\in M(\al'',\ga)$ representing the last term
in the chain-limit of $\om_n$. In fact, because no energy
can be lost at $\infty$ by the choice of $\eps$, there are, for any real number
$r$, connections $A_n,A$ representing $\om_n,\om$, respectively, such
that
\begin{equation}\label{eqn:anai}
  \|A_n-A\|_{L^{p,w}_1((r,\infty)\times Y)}\to0,
\end{equation}
as follows from the exponential decay results of \cite[Subsection~4.1]{D5}.
Here, $p,w$ are as in the definition of the space $\ca$ of connections
in Section~\ref{section:mod-sp}.

Suppose first that $\beta:=\al''$ is irreducible of index~$1$. Then
$\ttp(\om_n)=\tau^+_\beta(\om_n)$ for $n\gg0$ and
\[(\tau^+-\tau^+_\beta)(\om_n)=-\tau^+_\beta(\om_n)\to\infty,\]
proving the first assertion of the lemma.

Now suppose the sequence
$(\tau^+-\ttp)(\om_n)$ is not bounded. After passing to a subsequence we may
assume that there exists a $\beta\in\flag$ such that for each $n$ one has
$R^+(\om_n)\in W_\beta$. Suppose, for contradiction, that $\al''\neq\beta$.
Since $W_\beta$ is closed we must have $R^+(\om)\in W_\beta$
as well, hence
\[\cs(\al'',\ga)>\cs(\beta,\ga).\]
From \refp{eqn:anai} we deduce that
\[\tau^+_\beta(\om_n)\to\tau^+_\beta(\om),\]
so
$(\ttp-\tau^+)(\om_n)=\tau^+_\beta(\om_n)$ is bounded. This contradiction shows
that $\al''=\beta$.\square

\begin{lemma}\label{lemma:ttpconv}
  If $\om_n$ is a sequence in $M(\al',\ga)$ which converges over compacta
  to $\om\in M(\al'',\ga)$, where $\al''\in\fl(Y)$ and
  $\ind(\al'')\neq1$, then
  \[\ttp(\om_n)\to\ttp(\om).\]
  \end{lemma}

\proof Let $\beta\in\fl^1_Y$ with $\cs(\beta,\ga)>0$.
If $\cs(\al'',\ga)\le\cs(\beta,\ga)$ then $R^+(\om)\not\in W_\beta$. Since
$W_\beta$ is closed, we have $R^+(\om_n)\not\in W_\beta$ for $n\gg0$. This means
that $\beta$ contributes neither to $\ttp(\om)$ nor to $\ttp(\om_n)$ for
$n\gg0$. If on the other hand $\cs(\al'',\ga)>\cs(\beta,\ga)$ then
\[\tau^+_\beta(\om_n)\to\tau^+_\beta(\om).\]
From this the lemma follows.\square

Let $\tta$ and $\ttd$ be the real-valued functions on $M(\al,\ga)$ defined by
\[\tta:=\frac12(\ttp+\tau^-),\quad\ttd:=\frac12(\ttp-\tau^-).\]
Let
\[\lla:M(\al,\ga)\to[0,\infty),\quad
  \om\mapsto e_\al(R^-(\om))\cdot\ttd(\om),\]
where $e_\al$ is as in \refp{eqn:eal}. As the following lemma shows,
the quantity $\lla(\om)$ measures the extent to which $\om$
factors through the trivial connection $\theta$ over $Y$.

\be{lemma}\label{lemma:tilla-ch-conv}
Let $\om_n$ be a chain-convergent sequence  in $M(\al,\ga)$.
If the first term of the chain-limit of $\om_n$ lies in $\chm(\al,\theta)$ then
$\lla(\om_n)\to\infty$, 
otherwise the sequence $\lla(\om_n)$ is bounded.
\end{lemma}

\proof Because of the translationary invariance of $\lla$ we may assume
$\tau^-(\om_n)=0$ for all $n$,
so that the sequence $\om_n$ converges over compact subsets of $\ry$ to some
$\om\in M(\al,\beta)$, where $\beta\in\fl(Y)$. Then $\om$ represents
the first term of the chain-limit of $\om_n$.

{\bf Part I.} Suppose first that $\beta=\theta$. We will show that
$\lla(\om_n)\to\infty$.
There are two sequences $\tn1,\tn2$ of real numbers such that
\be{itemize}
\item $\ct^*_{\tn1}(\om_n)$ converges over compact subsets of $\ry$ to an
element of $M(\al,\theta)$.
\item $\ct^*_{\tn2}(\om_n)$ converges over compact subsets of $\ry$ to an
element of $M(\theta,\beta')$, where $\beta'$ is an element of $\fl^*(Y)$ which
is either equal to $\ga$ or has index~$1$.
\item $\tn2-\tn1\to\infty$.
\end{itemize}
Define the sequence $r_n$ of real numbers implictly by
\[\ce_{\om_n}((-\infty,r_n])=\cs(\al,\theta)+\eps.\]
Then $r_n<\tau^+(\om_n)$ and $r_n<\tau^+_\beta(\om_n)$ for all $\beta\in\fl_\ga$,
hence $r_n<\ttp(\om_n)$. For large $n$ one therefore has
\[\lla(\om_n)=\ttp(\om_n)-\tau^-(\om_n)>r_n-\tau^-(\om_n).\]
But
\[\tn1-\tau^-(\om_n),\quad\tn2-r_n\]
are both bounded sequences and $\tn2-\tn1\to\infty$, hence
\[\lla(\om_n)>r_n-\tau^-(\om_n)\to\infty.\]

{\bf Part II.} Now suppose $\beta$ is irreducible. We will show that
the sequence $\lla(\om_n)$ is bounded.

{\em Case 1:} $\beta=\ga$. Then $\om_n$ converges to $\om$ in
$M(\al,\ga)$, hence $\lla(\om_n)$ is bounded.

{\em Case 2:} $\ind(\al,\beta)\le4$. For large $n$ one would then have
$R^-(\om_n)\not\in U_\al$, hence $e_\al(R^-(\om_n))=0$ and therefore
$\lla(\om_n)=0$.

{\em Case 3:} $\ind(\al,\beta)=5$, i.e.\ $\ind(\beta)=1$.
For large $n$ one would then have
$R^+(\om_n)\in W_\beta$ and therefore
\[\lla(\om_n)=e_\al(\om_n[0])\cdot\tau^+_\beta(\om_n)
\to e_\al(\om[0])\cdot\tau^+(\om),\]
so that $\lla(\om_n)$ is bounded in this case, too.\square

Given $\al'\in\fl(Y)$, a real number $d$, and a real $3\times2$ matrix
$A'=(a'_{ij})$ of maximal rank we define two sections $\zeta_1,\zeta_2$ of
$\bbei{\al'}\ga$ by
\[\zeta_j(\om,t):=b^+\trp j+(1-b^+)\sum_{i=1}^3a'_{ij}\rho^+_i,\]
where $b^+:=b(\tau^+-\ttp-d)$.  {\em Here, and in the remainder of this
  section}, $b:\R\to\R$ is a smooth function satisfying \refp{eqn:b-prop1}
and \refp{eqn:b-prop2}.

We will show that for $\al'=\ga$ and generic matrix $A'$ the sections
$\zeta_1,\zeta_2$ are linearly independent at any point
$(\om,t)\in M(\al,\ga)\times\R$ with $\lla(\om)\gg0$. We begin by spelling
out sufficient conditions on $A'$ under which this holds.

For any $\beta\in\fl^1_Y$ the finite set
$\chm(\theta,\beta)\times\chm(\beta,\ga)$ is in $1-1$ correspondence with
the set of points $(\om,\om')\in M(\theta,\beta)\times M(\beta,\ga)$
satisfying
\begin{equation}\label{ttom}
  \tau^+(\om)=0=\tau^+(\om').
\end{equation}
(In other words, this is one way of fixing translation.) For each such pair
$(\om,\om')$, represented by a pair $(A,A')$ of connections, say, the holonomy
of $A$ along the path $[0,\infty)\times\{y_0\}$ composed with the holonomy of
$A'$ along $(-\infty,0]\times\{y_0\}$ defines an isomorphism
\[\hol_{\om,\om'}:\bbe_{\om[0]}\to\bbe_{\om'[0]}.\]
For any real number $r$ and $j=1,2$ let
\[\eta_j(r)=r\cdot\hol_{\om,\om'}(\rho_j(\om[0]))+
(1-r)\sum_{i=1}^3a'_{ij}\rho_i(\om'[0]).\]
Then the set
\[C:=\{r\in[0,1]\st\eta_1(r)\wedge\eta_2(r)=0\}\]
has expected dimension $1-2=-1$ and is empty for generic matrices $A'$.
Since $\fl(Y)$ is finite we conclude that for generic $A'$, the set $C$
is empty for any $\beta\in\fl^1_Y$ and any
$(\om,\om')\in M(\theta,\beta)\times M(\beta,\ga)$ satisfying \refp{ttom}.
From now on we assume $A'$ is chosen so that this holds.

\be{lemma}\label{lemma:tau1-lin-ind}
Let $A'$ be as described above.
If $d>0$ is sufficiently large then
the sections $\zeta_1,\zeta_2$ are linearly independent at every
point in $M(\theta,\ga)\times\R$.
\end{lemma}

\proof If the lemma were false then we could find a sequence $d_n$ of real
numbers converging to $\infty$ and for each $n$ an element
$\om_n\in M(\theta,\ga)$ such that $\zeta_1,\zeta_2$, defined with $d_n$
in place of $d$, are linearly dependent at $(\om_n,t)$ for some (hence any)
$t$. Because $A'$ has maximal rank and the assumptions on $\eps$ ensure that
$\rho_1,\rho_2,\rho_3$ are linearly independent at $R^+(\om_n)$, we must have
$b^+(\om_n)>0$, i.e.
\[(\tau^+-\ttp)(\om_n)>d_n-1,\]
which shows that $(\tau^+-\ttp)(\om_n)\to\infty$. After passing to a subsequence
we can assume that the sequence $\om_n$ is chain-convergent and that
$b^+(\om_n)$ converges to some $r\in[0,1]$. By Lemma~\ref{lemma:ttpinfty}
the chain-limit lies in $\chm(\theta,\beta)\times\chm(\beta,\ga)$ for some
$\beta\in\fl^1_Y$. Then the sequences
\[\ct^*_{\tau^+(\om_n)}(\om_n),\quad\ct^*_{\tau^+_\beta(\om_n)}(\om_n)\]
converge over compact subsets of $\ry$ to some $\om\in M(\theta,\beta)$ and
$\om'\in M(\beta,\ga)$, respectively, and \refp{ttom} holds. But then
$\eta_1(r)$ and $\eta_2(r)$ are linearly dependent, contradicting the assumption
on $A'$.\square

From now on we assume that $d$ is chosen so that the
conclusion of Lemma~\ref{lemma:tau1-lin-ind} holds.

\be{lemma}\label{lemma-zeta}
There is a constant $T_1<\infty$ such that the sections $\zeta_1,\zeta_2$
are linearly independent at every point $(\om,t)\in M(\al,\ga)\times\R$ with
$\lla(\om)>T_1$.
\end{lemma}

\proof Recall that if $\zeta_1,\zeta_2$ are linearly independent at
$(\om,t)$ for some real number $t$ then the same holds at $(\om,t')$ for all
$t'$. Now suppose the lemma
were false. Then we could find a sequence $\omn$ in $M(\al,\ga)$ such that
$\lla(\omn)\to\infty$ and $\zeta_1(\omn,t),\zeta_2(\omn,t)$ are linearly
dependent for every $n$. We may also arrange that $\tau^+(\om_n)=0$.
After passing to a subsequence we may assume that
$\omn$ is chain-convergent. From Lemma~\ref{lemma:tilla-ch-conv} we see that
there are two possibilities for the chain-limit.

{\bf Case 1:} The chain-limit of $\om_n$ lies in
$\chm(\al,\theta)\times\chm(\theta,\beta)\times\chm(\beta,\ga)$ for some
$\beta\in\fl^1_Y$. Then $\ttp(\om_n)=\tau^+_\beta(\om_n)$ for $n\gg0$.
Let $\om\in M(\theta,\beta)$ be a representative for the
middle term of the chain-limit. By Lemma~\ref{lemma:ttpinfty} we have
$(\tau^+-\ttp)(\om_n)\to\infty$, so for $t_n:=\ttp(\omn)$ one has
\[\zeta_j(\om_n,t_n)\to\rho_j(R^+(\om)),\]
contradicting the fact that the $\rho_j$ are linearly independent at
$R^+(\om)$.

{\bf Case 2:} The chain-limit of $\om_n$ lies in
$\chm(\al,\theta)\times\chm(\theta,\ga)$. Then $\om_n$ converges over compact
subsets of $\ry$ to some $\om\in M(\theta,\ga)$ satisfying
$\tau^+(\om)=0$. According to Lemma~\ref{lemma:ttpconv} we have
$\ttp(\om_n)\to\ttp(\om)$, so
\[\zeta_j(\om_n,t)\to\zeta_j(\om,t)\]
for any $t$. Hence, $\zeta_1,\zeta_2$ must be linearly dependent at
$(\om,t)$. But $d$ was chosen so that the conclusion of
Lemma~\ref{lemma:tau1-lin-ind} holds, so we have a contradiction.\square

At any point $(\om,t)\in M(\al',\ga)\times\R$ where $\zeta_1,\zeta_2$
are linearly independent let $\xi_1(\om,t),\xi_2(\om,t)$ be the
orthonormal pair of vectors in $\bbe_{\om[t]}$ obtained by applying the
Gram-Schmidt process to $\zeta_1(\om,t)$ and $\zeta_2(\om,t)$, and let
$\xi_3=\xi_1\times\xi_2$ be the fibrewise cross-product of $\xi_1$ and $\xi_2$.
Then $\{\xi_j(\om,t)\}_{j=1,2,3}$ is a positive orthonormal basis for
$\bbe_{\om[t]}$.

We now have the necessary ingredients to define the cut-down moduli space
$\tmx33(\al,\ga)$. Set
\[c:M(\al,\ga)\times\R\to[0,1],\quad(\om,t)\mapsto b(t-\tta(\om))\]
and for $j=1,2,3$ define a section $\si_j$ of the bundle $\bbe_{\al\ga}$ over
$M(\al,\ga)\times\R$ by
\[\si_j:=(1-c)\sum_ia^-_{ij}\rho^-_i+c\sum_ia^+_{ij}\xi_i.\]
Choose a constant $T_1$ for which the conclusion of Lemma~\ref{lemma-zeta}
holds and define a function $g:M(\al,\ga)\times\R\to[0,1]$ by
\begin{equation*}
  g(\om,t):=b(\lla(\om)-T_1)\cdot b(\ttp(\om)-t)\cdot b(t-\tau^-(\om)).
  \end{equation*}
For $j=1,2,3$ we now define a section $\ti s_j$ of $\bbe_{\al\ga}$ by
\[\ti s_j(\om,t):=(1-g(\om,t))\cdot s_j(\om[t])+g(\om,t)\cdot\si_j(\om,t).\]


Now set
\[\tmx33(\al,\ga):=\{(\om,t)\in\mab\times\R\st
\ti s_1(\om,-t)=0,\;\ti s_2(\om,t)=0\}.\]

In the study of the ends of $\tmx33(\al,\ga)$ we will encounter certain
subspaces of $M(\theta,\ga)$ which we now define. For $\om\in M(\theta,\ga)$
and $j=1,2$ set
\[\ti s_j(\om):=(1-b(\ttp(\om)))\cdot s_j(\om[0])
+b(\ttp(\om))\sum_{i=1}^3a^+_{ij}\xi_i(\om,0)\]
and define
\[\ti M_{3;j}(\theta,\ga):=\{\om\in M(\theta,\ga)\st \ti s_j(\om)=0\}.\]
This space has expected dimension $2-3=-1$ and is empty for ``generic''
choices of sections $s_j$ and matrix $A^+$.

\begin{lemma}\label{lemma:ttauest3}
  There is a constant $C_0<\infty$ such that for all
  $(\om,t)\in\tmx33(\al,\ga)$ one has
  \[|t|\le\min(-\tau^-(\om),\ttp(\om))+C_0.\] 
\end{lemma}

\proof That $|t|+\tau^-(\om)$ is uniformly bounded above for
$(\om,t)\in\tmx33(\al,\ga)$ is proved in the same way as the corresponding part
of Lemma~\ref{lemma:ttauest1}. To prove the same for $|t|-\ttp(\om)$,
suppose there were a sequence $(\om_n,t_n)\in\tmx33(\al,\ga)$ with
\[|t_n|-\ttp(\om_n)\to\infty.\]
After passing to a subsequence we may assume the following.
\begin{itemize}
\item The sequence $\om_n$ is chain-convergent;
\item There is a constant $e=\pm1$ such that $|t_n|=et_n$ for all $n$;
\item The sequence $et_n-\tau^+(\om_n)$ converges in $[-\infty,\infty]$ to
  some point $t$.
\end{itemize}
Let $j:=\frac12(3+e)$. Then for $n\gg0$ we have
\[0=\ti s_j(\om_n,et_n)=s_j(\om_n[et_n]).\]
According to Lemma~\ref{lemma:ttpinfty} one of the following two cases
must occur.

{\bf Case 1:} The sequence $(\tau^+-\ttp)(\om_n)$ is bounded. Then
$et_n-\tau^+(\om_n)\to\infty$, so $\om_n[et_n]\to\uga$. By continuity of
$s_j$ we must have $s_j(\uga)=0$, which however will not hold for a
``generic'' section $s_j$.

{\bf Case 2:} $(\tau^+-\ttp)(\om_n)\to\infty$. From Lemma~\ref{lemma:ttpinfty}
we deduce that $\ct^*_{\tau^+(\om_n)}(\om_n)$ converges over compact subsets of
$\ry$ to some $\om\in M(\beta,\ga)$, where $\beta\in\fl^1_Y$. Then
$\ttp(\om_n)=\tau^+_\beta(\om_n)$ for $n\gg0$. Furthermore,
$\ct^*_{\tau^+_\beta}(\om_n)$ converges over compacta to an element of some
moduli space $M(\al',\beta)$, where $\beta\neq\al'\in\fl(Y)$.

{\bf Case 2a:} $t=\pm\infty$. Then the exponential decay results of
\cite[Subsection~4.1]{D5} imply that
$\om_n[et_n]$ converges to $\ubeta$ (if $t=-\infty$) or to $\uga$
(if $t=\infty$). This is ruled out in the same way as Case 1.

{\bf Case 2b:} $t$ finite. Then $\ct^*_{et_n}(\om_n)$ converges over compacta
to $\om':=\ct^*_t(\om)\in M(\beta,\ga)$, and $\om_n[et_n]\to\om'[0]$.
But then $s_j(\om'[0])=0$, which will not hold for a ``generic'' section
$s_j$ of the bundle $\bbe$, since $M(\beta,\ga)$ has dimension~$1$
whereas $\bbe$ has rank~$3$.\square

\begin{lemma}\label{lemma:ttauest4}
  For any constant $C_1<\infty$ there is constant $L>0$ such that for
  all $(\om,t)\in\tmx33(\al,\ga)$ satisfying $\lla(\om)\ge L$ one has
  \[|t|\le\min(-\tau^-(\om),\ttp(\om))-C_1.\]
\end{lemma}

\proof If not, then there would be a constant $C_1<\infty$ and
a sequence $(\om_n,t_n)\in\tmx33(\al,\ga)$ with $\lla(\om_n)\to\infty$
such that either
\begin{description}
\item[(i)] $|t_n|>-\tau^-(\om_n)-C_1$ for all $n$, or
  \item[(ii)] $|t_n|>\ttp(\om_n)-C_1$ for all $n$.
\end{description}
Case (i) is rule out as in the proof of Lemma~\ref{lemma:ttauest2}. Now
suppose (ii) holds. Because $\lla(\om_n)\to\infty$ we have
$\ttd(\om_n)\to\infty$.
From Lemma~\ref{lemma:ttauest3} we deduce that $|t_n|-\ttp(\om_n)$ is bounded,
so
\[|t_n|-\tau^-(\om_n)\to\infty.\]
This implies that $c(\om_n,t_n)=1$ for $n\gg0$. After passing to a subsequence
we may assume that the sequence $\om_n$ chain-converges and
$|t_n|=-et_n$ for some constant $e=\pm1$.

{\bf Case 1:} $(\tau^+-\ttp)(\om_n)$ is bounded. By Lemmas~\ref{lemma:ttpinfty}
and \ref{lemma:tilla-ch-conv} the chain-limit of $\om_n$ must lie in
$\chm(\al,\theta)\times\chm(\theta,\ga)$, so after passing to a subsequence
we may assume that $\om'_n:=\ct^*_{et_n}(\om_n)$ converges over compacta to some
$\om\in M(\theta,\ga)$. Using Lemma~\ref{lemma:ttpconv} we obtain
\[g(\om_n,et_n)=b(\ttp(\om_n)-et_n)=b(\ttp(\om'_n))\to b(\ttp(\om)).\]
Let $j:=\frac12(3+e)$. Then
\[0=\ti s_j(\om_n,et_n)\to\ti s_j(\om).\]
But then $\om$ lies in $\ti M_{3;j}(\theta,\ga)$, which is empty by choice of the
matrix $A^+$.

{\bf Case 2:} $(\tau^+-\ttp)(\om_n)\to\infty$. Then the chain-limit of $\om_n$
lies in $\chm(\al,\theta)\times\chm(\theta,\beta)\times\chm(\beta,\ga)$ for
some $\beta\in\fl^1_Y$. For large $n$ we now have
$\ttp(\om_n)=\tau^+_\beta(\om_n)$ and $\xi_j(\om_n,et_n)=\trp j(\om_n,et_n)$,
$j=1,2$. After passing to a subsequence we may assume that
$\om'_n:=\ct^*_{et_n}(\om_n)$ converges over compacta to some
$\om\in M(\theta,\beta)$. For large $n$ we have
\[g(\om_n,et_n)=b(\tau^+_\beta(\om_n)-et_n)=b(\tau^+_\beta(\om'_n))\to
b(\tau^+(\om)).\]
Let $j:=\frac12(3+e)$. Then
\[0=\ti s_j(\om_n,et_n)\to(1-b(\tau^+(\om)))\cdot s_j(\om[0])
+b(\tau^+(\om))\sum_ia^+_{ij}\rho^+_i(\om,0).\]
Thus, $\om$ lies in $\ti M_3(\theta,\beta)$, which is empty by choice of $A^+$.
\square

\begin{lemma}
  There is a constant $L<\infty$ such that for all $(\om,t)\in\tmx33(\al,\ga)$
  one has $\lla(\om)<L$.
  \end{lemma}
\proof For any $(\om,t)\in\tmx33(\al,\ga)$ with $\lla(\om)>T_1$ let
$h(\om)\in\SO3$ be the matrix whose coefficients $h_{ij}(\om)$ are given by
\[\rho^-_j(\om,t)=\sum_ih_{ij}(\om)\xi_i(\om,t).\]
By Lemma~\ref{lemma:ttauest4} there is an $L\ge T_1+1$ such that for all
$(\om,t)\in\tmx33(\al,\ga)$ with $\lla(\om)\ge L$ one has
\[|t|\le\min(-\tau^-(\om),\ttp(\om))-1,\]
which implies that $\ti s_j(\om,t)=\si_j(\om,t)$. Given such a $(\om,t)$,
the triple
\[(H,v_1,v_2):=(h(\om),-t-\ti\tau_a(\om),t-\ti\tau_a(\om))\in
\SO3\times\R\times\R\]
satisfies the equation
\[(1-b(v_j))HA^-_j+b(v_j)A^+_j=0.\]
for $j=1,2$. However, as observed in the proof of
Proposition~\ref{prop:tmx-finite-set}, these equations have no solution
for generic matrices $A^\pm$.\square

We will now prove Proposition~\ref{prop:psichainmap} in degree $4$
by counting the number of ends of $\tmx33(\al,\ga)$ modulo $2$.

{\bf Ends of $\tmx33(\al,\ga)$:} Let $(\om_n,t_n)$ be a sequence in
$\tmx33(\al,\ga)$. After passing to a subsequence we may assume that
the following hold:
\begin{description}
\item[(i)] The sequences $\ct^*_{-t_n}(\om_n)$ and $\ct^*_{t_n}(\om_n)$
converge over compact subsets of $\ry$.
\item[(ii)] The sequence $\ct^*_{\tau^-(\om_n)}(\om_n)$ converges over compacta
to some $\om\in M(\al,\beta)$, where $\beta\in\fl(Y)$.
\item[(iii)] The sequences $t_n$ and $\tau^-(\om_n)$ converge in
$[-\infty,\infty]$.
\end{description}
Suppose $(\om_n,t_n)$ does not converge in $\tmx33(\al,\ga)$.

{\bf Case 1:} $\beta=\ga$. We show this cannot happen. First observe that
the sequence $\ttd(\om_n)$
converges in $\R$. Since Lemma~\ref{lemma:ttauest3} provides an upper bound
on $\tau^-(\om_n)$ and a lower bound on $\ttp(\om_n)$ it follows that
both sequences must be bounded. Applying the same lemma again we see that
$|t_n|$ is bounded. But then assumptions (ii) and (iii) imply that
$(\om_n,t_n)$ converges in $\tmx33(\al,\ga)$, which we assumed was not the
case.

{\bf Case 2:} $\beta$ irreducible, $\dim M(\al,\beta)\le4$. Then
$\lla(\om_n)=0$ for $n\gg0$. As in the proof of
Proposition~\ref{prop:v3v3} we find that the corresponding number of ends
of $\tmx33(\al,\ga)$ is $\la\psi d\al,\ga\ra$.

{\bf Case 3:} $\beta$ irreducible, $\dim M(\al,\beta)=5$. Then
$\ttp(\om_n)=\tau^+_\beta(\om_n)$ for $n\gg0$, and
\[\ttd(\om_n)\to\tau_d(\om).\]
As in Case~1 we see that the sequences $\tau^-(\om_n)$
and $t_n$ must be bounded, hence they both converge in $\R$ by assumption~(iii).
From (ii) we deduce that $\om_n$ converges over compacta to some
$\om'\in M(\al,\beta)$ (related to $\om$ by a translation).
By Lemma~\ref{lemma:ttpinfty} we have $\xi_j(\om_n,t)=\trp j(\om_n,t)$ for
$n\gg0$ and any $t$, so
\[\si_j(\om_n,t)\to\si_j(\om',t).\]
Setting $t':=\lim t_n$ we conclude that
$(\om',t')\in\tmx33(\al,\beta)$. The corresponding number of ends of
$\tmx33(\al,\ga)$ is $\la d\psi\al,\ga\ra$.\square

\subsection{Calculation of $\psi$}\label{subsec:calc-psi}

\begin{prop}\label{prop:psiv3v2}
  There are constants $\ka^\pm\in\z/2$ independent of $Y$ and satisfying
  $\ka^++\ka^-=1$ such that if $\psi$ is defined in terms of ``generic''
  sections $s_1,s_2$
  that are sufficiently close and $e$ is the sign of
  $\nu(A^+)-\nu(A^-)$ then there is a homomorphism
  $\Xi:C^*(Y)\to C^{*+4}(Y)$ such that
\begin{equation}\label{eqn:psi-v3v2}
  \psi=v_3v_2+\ka^e\del'\del+d\Xi+\Xi d,
\end{equation}
where the cup products $v_2,v_3$ are defined by three ``generic'' sections
of $\bbe$.
\end{prop}

To be precise, if $s'\in\Ga(\bbe)$ satisfies Property~$\pt4$ and
$\frP\subset\Ga(\bbe)$ is any sufficiently large finite-dimensional linear
subspace then for any sufficiently small generic
$(\frp_0,\frp_1)\in\frP\times\frP$
the conclusion of the proposition holds with $s_j=s'+\frp_j$.

The above proposition completes the proof of Proposition~\ref{prop:psixi}
except for the order of $v_2,v_3$, which is insignificant in vue of
Proposition~\ref{prop:v2v3phi}. (The order could be reversed by a small
change in the proof given below.)

\proof Let $\al,\beta\in\fl^*(Y)$ with $\ind(\al,\beta)=5$.
The proof is divided into two parts. The first part
deals with the case $\ind(\al)\not\equiv4\mod8$ in which
no factorization through
  the trivial connection can occur in the moduli space $M(\al,\beta)$. The
  second part handles the case $\ind(\al)\equiv4\mod8$.
  
  {\bf Part (I)} Suppose $\ind(\al)\not\equiv4\mod8$. The proof will consist of
  counting modulo $2$ the ends and boundary points of a $1$--manifold
  $\cm$ obtained by gluing together two $1$--manifolds
  $\cmsi$ and $\cmcyl$ along their common boundary.
  To define these $1$--manifolds, let
\[\Si:=\{z\in\bbc:|\Im(z)|\le3,\;|z|\ge1\}\]
and let $\Si':=\Si/\pm1$ be the surface-with-boundary obtained by identifying
each $z\in\Si$ with $-z$. The image of a point $z\in\Si$ in $\Si'$ will be
denoted by $[z]$.

\begin{figure}[h]
  \caption{A portion of the surface $\Si$}
\centerline{\includegraphics*[width=\textwidth]{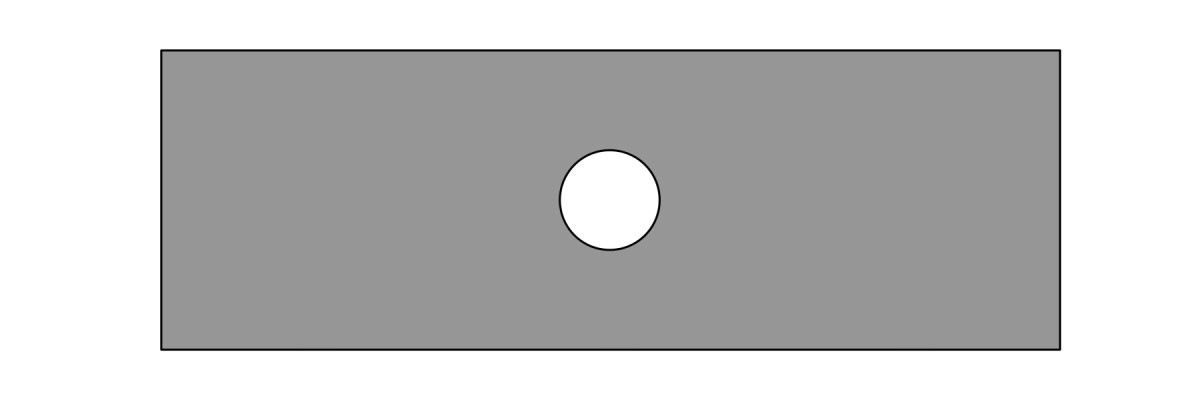}}
\end{figure}
For $-3\le y\le3$ we define a section $\chi_y$ of $\bbe$ by
\[6\chi_y:=(3-y)s_1+(3+y)s_2.\]
In particular,
\[\chi_{-3}=s_1,\quad\chi_3=s_2.\]
Let $\bar\xi\in\Ga(\bbe)$, and let $\hat\xi$ be a section of the bundle
$\bbe\times S^1$ over $\cb^*(Y[0])\times S^1$ satisfying
\[\hat\xi(\ups,-z)=-\hat\xi(\ups,z),\]
so that  $\hat\xi\in\Ga_a(\uline{\bbe})$ in the notation of
Section~\ref{sec:spaces-lin-dep}. We then define a section $\xi$ of the
bundle $\bbe\times\Si$ over $\cb^*(Y[0])\times\Si$ as follows. Let
\begin{equation}\label{eqn:b1def}
  b_1(z):=b(|z|-2).
  \end{equation}
For $\ups\in\cb^*(Y[0])$ and $z=(x,y)\in\Si$ let
\[\xi(\ups,z):=(1-b_1(z))\cdot(\bar\xi(\ups)+\hat\xi(\ups,z/|z|))
  +b_1(z)\chi_y(\ups).\]
  Let $f:\Si\to\R$ be the smooth function given by
  \[f(z):=b_1(z)\Re(z).\]
  Note that $f(z)=\Re(z)$ for $|z|\ge3$, and $f(z)=0$ for $|z|=1$.
  Moreover, $f(-z)=-z$.

\begin{defn}\label{defn:cm12}
  \begin{description}
  \item[(i)] Let $\cmsi=\cmsi(\al,\beta)$
    be the subspace of $M(\al,\beta)\times\Si'$
    consisting of those points $(\om,[z])$ such that
    \[\xi(\om[f(z)],z)=0,\quad\xi(\om[f(-z)],-z)=0.\]
  \item[(ii)] Let $\cmcyl=\cmcyl(\al,\beta)$ be the subspace of
    $M(\al,\beta)\times S^1\times[0,\infty)$
    consisting of those points $(\om,z^2,r)$ such that $z\in S^1$ and
    \[\hat\xi(\om[-r],z)=0,\quad\bar\xi(\om[r])=0.\]
   \end{description}
\end{defn}

If $\bar\xi$ is ``generic'' and $\hat\xi$ is given by a ``generic'' section
of $\ubbe\otimes\uell$ (see Lemma~\ref{lemma:Ga-a})
then $\cmcyl$ will be a smooth
$1$--manifold-with-boundary. Now choose a section $s'\in\Ga(\bbe)$ satisfying
Property~$\pt4$. If $\frP$ is a sufficiently large finite-dimensional
linear subspace of $\Ga(\bbe)$ and $(\frp_0,\frp_1)$ a generic element of
$\frP\times\frP$ then taking $s_j=s'+\frp_j$, $j=1,2$ the space $\cmsi$ will be
a smooth $1$--manifold-with-boundary. If in addition $\frp_0,\frp_1$ are
sufficiently small then for $-3\le y\le3$ the section $\chi_y$
will satisfy Property~$\pt4$ and define the same cup product
$v_3:C^*(Y)\to C^{*+3}(Y)$ as $s'$, by Lemma~\ref{lemma:prop-T}.

The part of the boundary of $\cmsi$ given by $|z|=1$ can be identified with
the boundary of $\cmcyl$ (defined by $r=0$). To see this, let
$(\om,z)\in M(\al,\beta)\times\Si$ with $|z|=1$, and set $\om_0:=\om[0]$. Then
$(\om,[z])\in\cmsi$ if and only if
\[\bar\xi(\om_0)+\hat\xi(\om_0,z)=0=\bar\xi(\om_0)-\hat\xi(\om_0,z),\]
which in turn is equivalent to $(\om,z^2,0)\in\cmcyl$.

This allows us to define a topological $1$--manifold-with-boundary
$\cm=\cm(\al,\beta)$ as a quotient
of the disjoint union $\cmsi\coprod\cmcyl$ by identifying each boundary point
of $\cmcyl$ with the corresponding boundary point of $\cmsi$.

The proposition will be proved by counting the ends and boundary points
of $\cm$ modulo~$2$. Before doing this, we pause to define the homomorphism
$\Xi$. Let $\al',\beta'\in\fl^*(Y)$ with $\ind(\al',\beta')=4$.
Replacing $(\al,\beta)$ by $(\al',\beta')$ in Definition~\ref{defn:cm12}
yields zero-dimensional manifolds $\cm_j(\al',\beta')$, $j=1,2$.
The argument that we will give below to determine the ends of $\cm_j(\al,\beta)$
can also be applied to show that $\cm_j(\al',\beta')$ is compact.
Granted this, we define $\Xi:=\Xi_1+\Xi_2$, where $\Xi_j$ has matrix coefficient
\[\la\Xi_j\al',\beta'\ra:=\#\cm_j(\al',\beta').\]

{\bf Ends of $\cmsi(\al,\beta)$:} Let $(\om_n,[z_n])$ be a sequence in
$\cmsi(\al,\beta)$, where $z_n=(x_n,y_n)\in\R^2$. After passing to a subsequence
we may assume that

\be{description}
\item[(i)] The sequence $\ct^*_{-x_n}(\om_n)$
  converges over compact subsets of $\ry$ to some
  $\om^-\in M(\al^-,\beta^-)$.
\item[(ii)] The sequence $\ct^*_{x_n}(\om_n)$ converges over compact subsets of
$\ry$ to some $\om^+\in M(\al^+,\beta^+)$.
\item[(iii)] The sequence $(x_n,y_n)$ converges in
  $[-\infty,\infty]\times[-3,3]$ to some point $(x,y)$.
\end{description}
Suppose $(\om_n,[z_n])$ does not converge in $\cmsi(\al,\beta)$.

{\em Case 1:} $x$ finite. Then $\ind(\al^+,\beta^+)=4$ and either
$\al^+=\al$ or $\beta^+=\beta$. The corresponding number of ends of
$\cmsi(\al,\beta)$ is $\la(d\Xi_1+\Xi_1d)\al,\beta\ra$ modulo $2$.

{\em Case 2:} $x=\pm\infty$. Then for $n\gg0$ one has
\[0=\xi(\om[\pm x_n],\pm z_n)\to\chi_{\pm y}(\om^\pm[0]).\]
Hence $\chi_{\pm y}(\om^\pm[0])=0$. Since
$\chi_{\pm y}$ satisfy Property~$\pt4$ we must have $\ind(\al^\pm,\beta^\pm)\ge3$,
so
\[5=\ind(\al,\beta)\ge\ind(\al^-,\beta^-)+\ind(\al^+,\beta^+)\ge6.\]
This contradiction shows that there are no ends in the case $x=\pm\infty$.

{\bf Ends of $\cmcyl(\al,\beta)$:} We argue as in part~(I) of the proof of
Proposition~\ref{prop:v2v3phi}. Let $(\om_n,z_n^2,r_n)$ be a sequence in
$\cmcyl(\al,\beta)$. 
After passing to a subsequence we may assume that $r_n$ converges in
$[0,\infty]$ to some point $r$. Then the number of ends modulo~$2$
corresponding to $r<\infty$ is $\la(d\Xi_2+\Xi_2d)\al,\beta\ra$.
Using Proposition~\ref{prop:E-lla} and the description of the cup product
$v_2$ in Lemma~\ref{lemma:v_2-defn-alt} we see that
the number of ends corresponding to $r=\infty$ is $\la v_3v_2\al,\beta\ra$.

{\bf Boundary points of $\cm(\al,\beta)$:} These are the points
$(\om,[z])$ in $M(\al,\beta)\times\Si'$ where $\Im(z)=3$ and
\[0=\xi(\om[x],z)=s_2(\om[x]),\quad0=\xi(\om[-x],-z)=s_1(\om[-x]).\]
The number of such points is by definition $\la\psi\al,\beta\ra$.

Since the number of ends plus the number of boundary point of $\cm$ must
be zero modulo~$2$ we obtain the equation~\refp{eqn:psi-v3v2} in the
case $\ind(\al)\not\equiv4\mod2$.

{\bf Part (II)} Suppose $\ind(\al)\equiv4\mod8$. We adapt the approach used
in the proof of Proposition~\ref{prop:v2v3phi} by deforming the equations
defining the $1$--manifold $\cm(\al,\beta)$ to obtain a $1$--manifold
in which we can control factorizations through the trivial connection.
Cutting away part of the latter $1$--manifold yields a
$1$--manifold-with-boundary $\cm^L$ in which such factorizations do not occur.
Counting modulo~$2$ the ends and boundary points of $\cm^L$ will produce
the formula \refp{eqn:psi-v3v2}.

We continue using the notation introduced in Part~(I) and earlier in
Subsections~\ref{subsec:enhol} - \ref{subsec:proof-v2v3phi}.

In order to define the deformed equations we first introduce maps
$V^\pm:[-3,3]\to\R^3$ given by
\[6V^\pm(y):=(3-y)A^\pm_1
+(3+y)A^\pm_2,\]
where $A^\pm_j$ is the $j$th column of the matrix $A^\pm$ entering in the
definition of $\psi$, see Subsection~\ref{subsec:cochain-map-psi} and
Equation~\ref{eqn:rhorho}.

Choose generic elements
$\bar L^\pm\in\R^3$ and functions $\hat L^\pm:S^1\to\R^3$ satisfying
$\hat L^\pm(-z)=-\hat L^\pm(z)$ for $z\in S^1$. We define maps
$L^\pm:\Si\to\R^3$ by
\[L^\pm(z):=(1-b_1(z))\cdot (\bar L^\pm+\hat L^\pm(z/|z|))+
  b_1(z)\cdot V^\pm(\Im(z)),\]
where the function $b_1$ is as in \refp{eqn:b1def}. Let $\bbesi(\al,\beta)$
be the vector bundle over $\Si\times\mab$ obtained by pulling back
the bundle $\bbe\to\cb^*(Y[0])$ by the map
\[\mab\times\Si\to\cb^*(Y[0]),\quad(\om,z)\mapsto\om[f(z)].\]
Let $c$ and $g$ be the functions defined in
\refp{eqn:c23def} and \refp{eqn:gomt}, respectively.
We define sections $\si,s$ of $\bbesi(\al,\beta)$ by
\begin{gather*}
\si(\om,z):=(1-c(\om,f(z)))\sum_{i=1}^3L^-_i(z)\rho^-_i(\om,f(z))\\
{\vphantom{x}}+c(\om,f(z))\sum_{i=1}^3L^+_i(z)\rho^+_i(\om,f(z)),\\
s(\om,z):=(1-g(\om,f(z)))\cdot\xi(\om[f(z)],z)+g(\om,f(z))\cdot\si(\om,z).
\end{gather*}

\begin{defn}\label{def:tcmsi}
Let $\tcmsi=\tcmsi(\al,\beta)$ be the subspace of $\mab\times\Si'$ consisting
of those points $(\om,[z])$ such that
\[s(\om,z)=0,\quad s(\om,-z)=0.\]
\end{defn}

We define sections $\bar\si,\bar s$ of the bundle $\bbe(\al,\beta)$
over $\mab\times\R$ by
\begin{gather*}
\bar\si(\om,r):=(1-c(\om,r))\sum_{i=1}^3\bar L^-_i\rho^-_i(\om,r)
+c(\om,r)\sum_{i=1}^3\bar L^+_i\rho^+_i(\om,r),\\
\bar s(\om,r):=(1-g(\om,r))\cdot\bar\xi(\om[r])+g(\om,r)\cdot\bar\si(\om,r).
\end{gather*}
Let $\hat\bbe(\al,\beta)$ be the vector bundle over
$\mab\times S^1\times\R$ obtained by pulling back the bundle $\bbe$
by the map
\[\mab\times S^1\times\R\to Y[0],\quad(\om,z,r)\mapsto\om[r].\]
We define sections $\hat\si,\hat s$ of $\hat\bbe(\al,\beta)$ by
\begin{gather*}
\hat\si(\om,z,r):=(1-c(\om,r))\sum_{i=1}^3\hat L^-_i(z)\rho^-_i(\om,r)
+c(\om,r)\sum_{i=1}^3\hat L^+_i(z)\rho^+_i(\om,r),\\
\hat s(\om,z,r):=(1-g(\om,r))\cdot\hat\xi(\om[r],z)+g(\om,r)\cdot\hat\si(\om,z).
\end{gather*}
Note that $\hat s(\om,-z,r)=-\hat s(\om,z,r)$.

\begin{defn}\label{def:tcmcyl}
Let $\tcmcyl=\tcmcyl(\al,\beta)$ be the subspace of
$\mab\times S^1\times[0,\infty)$
consisting of those points $(\om,z^2,r)$ such that $z\in S^1$ and
\[\hat s(\om,z,-r)=0,\quad\bar s(\om,r)=0.\]
\end{defn}

By inspection of the formulas involved one finds that for $|z|=1$ one has
\begin{align*}
\bar\si(\om,0)+\hat\si(\om,z,0)&=\si(\om,z),\\
\bar s(\om,0)+\hat s(\om,z,0)&=s(\om,z).
\end{align*}
Therefore, the part of the boundary of $\tcmsi$ given by $|z|=1$ can be
identified with the boundary of $\tcmcyl$ (defined by $r=0$). By gluing
$\tcmsi$ and $\tcmcyl$ correspondingly
we obtain a topological $1$--manifold-with-boundary $\tcm$.

\begin{lemma}\label{lemma:rez-est1}
  There is a constant $C_0<\infty$ such that for all
  $(\om,[z])\in\tcmsi$ one has
  \[|f(z)|\le\min(-\tau^-(\om),\tau^+(\om))+C_0.\] 
\end{lemma}

\proof The proof is similar to that of Lemma~\ref{lemma:ttauest1}.
We must provide upper bounds on both quantities $|f(z)|+\tau^-(\om)$ and
$|f(z)|-\tau^+(\om)$ for $(\om,[z])\in\tcmsi$. 
The proof is essentially the same in both cases, so we will only spell it out
in the second case. Suppose, for contradiction, that $(\om_n,[z_n])$
is a sequence in $\tcmsi$ with $|f(z)|-\tau^+(\om_n)\to\infty$.
By perhaps replacing $z_n$ by $-z_n$ we can arrange that
$\Re(z_n)\ge0$. Then $f(z_n)\ge0$ as well,
and $g(\om_n,f(z_n))=0$ for $n\gg0$. Let $z_n=(x_n,y_n)$.
After passing to a subsequence we may assume
that $z_n$ converges in $[0,\infty]\times[-3,3]$ to some
point $(x,y)$.

{\bf Case 1:} $x$ finite. Let $z:=(x,y)\in\Si$. The sequence $\om_n$
converges to $\ubeta$ over compact subsets of $\ry$, so 
for large $n$ we have
\[0=\xi(\om_n[f(z_n)],z_n)\to\xi(\ubeta,z).\]
However, the space of all $w\in\Si$ for which $\xi(\ubeta,w)=0$ has
expected dimension $2-3=-1$, so this space is empty for ``generic''
sections $s_1,s_2,\bar\xi,\hat\xi$. Hence, $x$ cannot be finite.

{\bf Case 2:} $x=\infty$. Then $f(z_n)=x_n$ for large $n$.
Now, $\ct^*_{x_n}\om_n$ converges over compacta
to $\ubeta$, so for large $n$ we have
\[0=\xi(\om_n[x_n],z_n)=\chi_{y_n}(\om_n[x_n])\to\chi_y(\ubeta).\]
However, the space of all $t\in[-3,3]$ for which $\chi_t(\ubeta)=0$ has
expected dimension $1-3=-2$, so this space is empty for ``generic''
sections $s_1,s_2$. Hence, $x\not=\infty$.

This contradiction proves the lemma.\square

In the proof of Lemma~\ref{lemma:rez-est2} below
we will encounter certain limits
associated to sequences in $\tcmsi$ with chain-limits in
$\chm(\al,\theta)\times\chm(\theta,\beta)$. These limits lie in 
cut down moduli spaces analogous to those introduced in
Definitions~\ref{def:tcmsi} and \ref{def:tcmcyl},
with $M(\al,\theta)$ or $M(\theta,\beta)$
in place of $\mab$. We now define these cut-down spaces in the case of
$M(\theta,\beta)$ and observe that they are ``generically'' empty.
The case of $M(\al,\theta)$ is similar.

For any $(\om,z)\in\mtb\times\Si$ let
\begin{align*}
s(\om,z):=&(1-b(\tau^+(\om)-f(z)))\cdot\xi(\om[f(z)],z)\\
&\vphantom{x}+b(\tau^+(\om)-f(z))\sum_{i=1}^3L^+_i(z)\rho^+_i(\om,f(z)).
\end{align*}

\begin{defn}
Let $\tcmsi(\theta,\beta)$ be the subspace of $M(\theta,\beta)\times\Si'$
consisting
of those points $(\om,[z])$ such that
\[s(\om,z)=0,\quad s(\om,-z)=0.\]
\end{defn}

Then $\tcmsi(\theta,\beta)$ has expected dimension $3-6=-3$ and is
empty for ``generic'' sections $s_1,s_2,\bar\xi,\hat\xi$ and generic choices of
$A^+,\bar L^+,\hat L^+$.

\begin{defn}
Let $\tcmint(\theta,\beta)$ be the subspace of
$M(\theta,\beta)\times[-3,3]$ consisting of those points $(\om,y)$ such that
\[(1-b(\tau^+(\om)))\cdot\chi_y(\om[0])
+b(\tau^+(\om))\sum_iV^+_i(y)\rho^+_i(\om,0)=0.\]
\end{defn}

We observe that the space $\tcmint(\theta,\beta)$ (a parametrized version of
the space $\ti M_3(\theta,\beta)$ defined in
Subsection~\ref{subsec:proof-v2v3phi})
has expected dimension $2-3=-1$ and is
empty for ``generic'' sections $s_1,s_2$ and generic matrix
$A^+$.

\begin{lemma}\label{lemma:rez-est2}
  For any constant $C_1<\infty$ there is constant $L>0$ such that for
  all $(\om,[z])\in\tcmsi$ satisfying $\lla(\om)\ge L$ one has
  \[|f(z)|\le\min(-\tau^-(\om),\tau^+(\om))-C_1.\]
\end{lemma}

\proof The proof is similar to that of Lemma~\ref{lemma:ttauest2}. If the lemma
did not hold there would be a sequence $(\om_n,[z_n])$ in $\tcmsi$ such that
$\lla(\om_n)\to\infty$ and one of the following two conditions hold:
\begin{description}
\item[(i)] $|f(z_n)|>-\tau^-(\om_n)-C_1$ for all $n$,
  \item[(ii)] $|f(z_n)|>\tau^+(\om_n)-C_1$ for all $n$.
\end{description}
Suppose (ii) holds, the other case being similar.
By replacing $z_n$ by $-z_n$, if necessary, we can arrange that $\Re(z_n)\ge0$.
From Lemma~\ref{lemma:rez-est1} we deduce that the sequence
$f(z_n)-\tau^+(\om_n)$ is bounded, whereas
\[f(z_n)-\tau^-(\om_n)\to\infty.\]
For large $n$ we therefore have
\[c(\om_n,f(z_n))=1,\quad g(\om_n,f(z_n))=b(\tau^+(\om_n)-f(z_n)).\]
Let $z_n=(x_n,y_n)$.
After passing to a subsequence we may assume that
\begin{itemize}
\item $\om'_n:=\ct^*_{x_n}\om_n$ converges over compact subsets of $\ry$ to
some $\om'\in M(\theta,\beta)$;
\item $z_n$ converges in $[0,\infty]\times[-3,3]$ to some point $z=(x,y)$.
  \end{itemize}

{\bf Case 1:} $x$ finite. Then $\om_n$ converges over compacta to some
$\om\in\mtb$, and
\[0=s(\om_n,z_n)\to s(\om,z).\]
Beause the sequence $z_n$ is bounded, we also have $c(\om_n,f(-z_n))=1$ for
large $n$, so
\[0=s(\om_n,-z_n)\to s(\om,-z).\]
But then $(\om,[z])$ belongs to $\tcmsi(\theta,\beta)$, contradicting the fact
that that space is empty.

{\bf Case 2:} $x=\infty$. Since
\[\tau^+(\om'_n)=\tau^+(\om_n)-x_n,\]
we obtain
\[g(\om_n,f(z_n))=b(\tau^+(\om'_n))\quad\text{for $n\gg0$.}\]
Therefore,
\[0=s(\om_n,z_n)\to
(1-b(\tau^+(\om')))\cdot\chi_y(\om'[0])
+b(\tau^+(\om'))\sum_iV^+_i(y)\rho^+_i(\om',0).\]
But this means that $(\om',y)$ belongs to $\tcmint(\theta,\beta)$, which is
empty.

This contradiction proves the lemma.\square

\begin{lemma}\label{lemma:r-est1}
  There is a constant $C_0<\infty$ such that for all
  $(\om,z^2,r)\in\tcmcyl$ one has
\[r\le\min(-\tau^-(\om),\tau^+(\om))+C_0.\]
\end{lemma}

\proof This is similar to the proof of Lemma~\ref{lemma:ttauest1}.\square

\begin{lemma}\label{lemma:r-est2}
  For any constant $C_1<\infty$ there is constant $L>0$ such that for
  all $(\om,z^2,r)\in\tcmcyl$ satisfying $\lla(\om)\ge L$ one has
  \[r\le\min(-\tau^-(\om),\tau^+(\om))-C_1.\]
\end{lemma}

\proof This is similar to the proof of Lemma~\ref{lemma:ttauest2}.\square

Choose $L\ge2$ such that the conclusions of Lemmas~\ref{lemma:rez-est2}
and \ref{lemma:r-est2} hold with $C_1=1$. For all $(\om,[z]\in\tcmsi$
with $\lla(\om)\ge L$ we then have
\[s(\om,z)=\si(\om,z),\]
and for all $(\om,z^2,r)\in\tcmcyl$ with $\lla(\om)\ge L$ we have
\[\hat s(\om,z,-r)=\hat\si(\om,z,-r),\quad\bar s(\om,r)=\bar\si(\om,r).\]
From Lemma~\ref{lemma:glthm} it follows that
$L$ is a regular value of the real functions on
$\tcmsi$ and $\tcmcyl$ defined by $\lla$. Therefore,
\begin{gather*}
\cmlsi:=\{(\om,[z])\in\tcmsi\st\lla(\om)\le L\},\\
\cmlcyl:=\{(\om,z^2,r)\in\tcmcyl\st\lla(\om)\le L\}
\end{gather*}
are smooth $1$--manifolds-with-boundary, and
\[\cm^L:=\cmlsi\cup\cmlcyl\]
is a topological $1$--manifold-with-boundary. (As before we identify the
part of $\cmlsi$ given by $|z|=1$ with the part of $\cmlcyl$ given by $r=0$.)

{\bf Ends of $\cm^L$:} From Lemma~\ref{lemma:rez-est1} we deduce that every
sequence $(\om_n,[z_n])$ in $\cmlsi$ which satisfies $\lla(\om_n)>0$ has a
convergent subsequence. Similarly, it follows from Lemma~\ref{lemma:r-est1}
that every sequence $(\om_n,z_n^2,r_n)$ in $\cmlcyl$ with $\lla(\om_n)>0$ has a
convergent subsequence. (See the proof of Proposition~\ref{prop:v2v3phi},
``Ends of $\mxl23(\al,\beta)$'', Case 2.) Therefore, all ends of $\cm^L$
are associated with sequences on which $\lla=0$. 
The number of such ends,
counted modulo~$2$, is given by the same formula as in Part~(I), namely
\[\la(v_3v_2+d\Xi+\Xi d)\al,\beta\ra.\]

{\bf Boundary points of $\cm^L$:} The boundary of $\cm^L$ decomposes as
\[\prtl\cm^L=\wsi\cup\wsi'\cup\wcyl,\]
where $\wsi$ and $\wcyl$ are the parts of the boundaries of
$\cmlsi$ and $\cmlcyl$, respectively, given by $\lla(\om)=L$, and $\wsi'$
is the part of the boundary of $\cmlsi$ given by $\Im(z)=\pm3$.
By choice of matrices $A^\pm$ there are no points $(\om,t)\in\tmx33(\al,\beta)$
with $\lla(\om)\ge L$, hence $W'_\Si=\tmx33(\al,\beta)$ and
\[\#\wsi'=\la\psi\al,\beta\ra.\]
By Lemma~\ref{lemma:glthm} we can identify
\[\wsi=\chm(\al,\theta)\times\chm(\theta,\beta)\times\nsi,\quad
\wcyl=\chm(\al,\theta)\times\chm(\theta,\beta)\times\ncyl,\]
where $\nsi$ is the set of points $(H,\tau,[z])$ in
$\SO3\times\R\times\Si'$ satisfying
\begin{gather*}
(1-b(f(z)-\tau))HL^-(z)+b(f(z)-\tau)L^+(z),\\
(1-b(f(-z)-\tau))HL^-(-z)+b(f(-z)-\tau)L^+(-z),
\end{gather*}
whereas $\ncyl$ is the set of points $(H,\tau,z^2,r)$ in
$\SO3\times\R\times S^1\times[0,\infty)$ satisfying
\begin{gather*}
(1-b(-r-\tau))H\hat L^-(z)+b(-r-\tau)\hat L^+(z)=0,\\
(1-b(r-\tau))H\bar L^-+b(r-\tau)\bar L^+=0.
\end{gather*}
Here, $(H,\tau)$ corresponds to $(h(\om),\tau_a(\om))$.
It follows from these descriptions that
\[\#(\wsi\cup\wcyl)=\ka\la\del'\del\al,\beta\ra,\]
where $\ka=\#(\nsi\cup\ncyl)\in\z/2$ is independent of the manifold $Y$.

To prove the theorem it only remains to understand the dependence of
$\ka$ on the pair of matrices $A=(A^+,A^-)$. To emphasize the dependence on $A$
we write $\ka=\ka(A)$ and $\nsi=\nsi(A)$. The space $\ncyl$ is
independent of $A$. The part of $\nsi$ corresponding to $|z|=1$ is also
independent of $A$ and is empty for generic $\bar L,\hat L$ for dimensional
reasons.

Let $P$ denote the space of all pairs $(B^+,B^-)$ of $3\times2$ real
matrices with non-zero columns $B^\pm_j$. Let
\[P^\pm:=\{(B^+,B^-)\in P\st\pm(\nu(B^+)-\nu(B^-))>0\},\]
where $\nu$ is as in \refp{eqn:nuB}. Note that each of $P^+,P^-$ is homotopy
equivalent to $S^2\times S^2$ and therefore path connected.

For any smooth path $C:[0,1]\to P$ we define
\[\cnc:=\bigcup_{0\le t\le1}\nsi(C(t))\times\{t\}\subset
\SO3\times\R\times\Si'\times[0,1].\]
As observed above there are no points $(H,\tau,[z],t)$ in $\cnc$ with $|z|=1$.
Since $b_1(z)>0$ for $|z|>1$ we can therefore make $\cnc$ regular
(i.e.\ transversely cut out) by varying $C$ alone. If $\cnc$ is regular
then it is a compact $1$--manifold-with-boundary, and
\[\prtl\cnc=\nsi(C(0))\cup\nsi(C(1))\cup X_C,\]
where $X_C$ is the set of points $(H,\tau,x,t)$ in
$\SO3\times\R\times\R\times[0,1]$ satisfying the two equations
\begin{gather*}
  (1-b(x-\tau))HC^-_1(t)+b(x-\tau)C^+_1(t)=0,\\
  (1-b(-x-\tau))HC^-_2(t)+b(-x-\tau)C^+_2(t)=0.
  \end{gather*}
It follows that
\[\ka(C(0))+\ka(C(1))=\#X_C.\]
If $A,B\in P^+$ then we can find a path $C:[0,1]\to P^+$ from $A$ to $B$.
Then $X_C$ is empty. By perturbing $C(t)$ for $0<t<1$ we can arrange that
$\cnc$ is regular. This yields $\ka(A)=\ka(B)$. The same holds if
$A,B\in P^-$.

Let $\ka^\pm$ be the value that $\ka$ takes on $P^\pm$. To compute
$\ka^++\ka^-$, let $(e_1,e_2,e_3)$ be the standard basis for $\R^3$ and define
$C:[0,1]\to P$ by
\begin{align*}
  -C^+_1(t)&=C^-_1(t):=e_1,\\
  -C^+_2(t)&:=(1-t)e_1+te_2,\\
  C^-_2(t)&:=(1-t)e_2+te_1.
\end{align*}
Then $C(0)\in P^+$ and $C(1)\in P^-$. Moreover,
$X_C$ consists of the single point
$(I,0,0,1/2)$, and this point is regular. (Here $I$ is the identity matrix.)
If we perturb $C$ a little in order
to make $\cnc$ regular then $X_C$ will still consist of a single, regular point.
We conclude that
\[\ka^++\ka^-=\#X_C=1.\]
This completes the proof of the proposition.\square

\appendix

\section{Instantons reducible over open subsets}

The following proposition is implicit in \cite[p\,590]{KM3} but we include a
proof for completeness.

\begin{prop}
  Let $X$ be an oriented connected Riemannian $4$--manifold and $E\to X$
  an oriented Euclidean $3$--plane bundle. Suppose $A$ is a non-flat
  ASD connection in $E$ which restricts to a reducible connection over some
  non-empty open set in $X$. Then there exists a rank~$1$ subbundle of $E$
  which is preserved by $A$.
\end{prop}

\proof This is a simple consequence of the unique continuation argument in the
proof of \cite[Lemma~4.3.21]{DK}. The proof has two parts: local existence
and local uniqueness.

(i) Local existence. By unique continuation, every point in $X$ has a connected
open neighbourhood $V$ such that $A|_V$ is reducible, i.e.\ there exists
a non-trivial automorphism $u$ of $E|_V$ such that $\nabla_Au=0$. The
$1$--eigenspace of $u$ is then a line bundle preserved by $A$.

(ii) Local uniqueness. Because $A$ is not flat, it follows from
unique continuation that the set of points in $X$ where $F_A=0$ has empty
interior. Now let $V$ be any non-empty connected open set in $X$ and suppose
$A$ preserves a rank~$1$ subbundle $\lla\subset E|_V$. We show that $\lla$ is
uniquely determined. Let $x\in V$ be a point where $F_A\neq0$. By the
holonomy description of curvature
(see \cite[Theorem~12.47]{Jeffrey-Lee-Manifolds-DG}) we can find a loop $\ga$
in $V$ based at $x$ such that the holonomy $\hol_\ga(A)$ of $A$ along $\ga$
is close to but different from the identity. The $1$--eigenspace of
$\hol_\ga(A)$ is then $1$--dimensional and must agree with the fibre $\lla_x$.
If $x'$ is an arbitrary point in $V$ then there is a similar description of
$\lla_{x'}$ in terms of the holonomy of $A$ along a loop obtained by
conjugating $\ga$ with a path in $V$ from $x$ to $x'$.
\square

\section{Unique continuation on a cylinder}

As in Subsection~\ref{subsec:inst-cohom} let 
$Y$ be a closed oriented connected $3$-manifold and $P\to Y$ an
$\SO3$ bundle. If $Y$ is not an integral homology sphere then we assume
$P$ is admissible.
Let $J\subset\R$ be an open interval.
We consider the perturbed ASD equation for connections
in the bundle $J\times P\to J\times Y$ obtained by
adding a holonomy perturbation to the Chern-Simons function. For a connection
$A$ in temporal gauge the equation takes the form
\[\frac{\prtl A_t}{\prtl t}=-*F(A_t)+V(A_t),\]
where $A_t$ is the restriction of $A$ to the slice $\{t\}\times P$ and
$V$ is the formal gradient of the perturbation.
The following proposition is probably well known among experts, but we
include a proof for completeness.

\begin{prop}\label{prop:unique-continuation-cylinder}
  Suppose $A,A'$ are perturbed ASD connections in the bundle
  $J\times P\to J\times Y$. If $A$ and $A'$ are in temporal gauge and
$A_T=A'_T$ for some $T\in J$, then $A=A'$.
\end{prop}

\proof We will apply (an adaption of)
the abstract
unique continuation theorem in \cite{Ogawa}. To this end, fix an arbitrary
connection $B$ in $P$ and let
\[c_t=A_t-A'_t,\quad a_t=A_t-B,\quad a'_t=A'_t-B.\]
We have
\[F(A_t)=F(B)+d_Ba_t+a_t\wedge a_t\]
and similarly for $A'_t$, so
\[\frac{\prtl c_t}{\prtl t}+*d_Bc_t=-*(a_t\wedge c_t+c_t\wedge a'_t)
+V(A_t)-V(A'_t).\]
By \cite[Prop.\ 3.5\,(v)]{KM7} we have
\[\|V(A_t)-V(A'_t)\|_{L^2}\le\const\|c_t\|_{L^2},\]
hence
\[\|\frac{\prtl c_t}{\prtl t}+*d_Bc_t\|_{L^2}\le\phi(t)\|c_t\|_{L^2}\]
where
\[\phi(t)=\const(\|a_t\|_\infty+\|a'_t\|_\infty+1).\]
Because $*d_B$ is a formally self-adjoint operator on $1$--forms on $Y$ and
$\phi$ is locally square integrable (in fact, continuous), we deduce
from \cite{Ogawa} that for any
compact subinterval $[t_0,t_1]$ of $J$
there are constants $C_0,C_1$ such that for $t_0\le t\le t_1$ one has
\[\|c_t\|_{L^2}\ge\|c_{t_0}\|_{L^2}\cdot\exp(C_0t+C_1).\]
(\cite{Ogawa} considers the case when $c_t$ is defined for $0\le t<\infty$, but
the approach works equally well in our case.)
Taking $t_1=T$ we obtain $c_t=0$ for
$t<T$. Replacing $c_t$ by $c_{-t}$ we get $c_t=0$ for
$t>T$ as well.\square

\noindent University of Oslo, Norway\\
\noindent Email:\ kfroyshov@math.uio.no

\end{document}